\documentclass{article}
 \usepackage{amsmath,graphicx,amsfonts,amssymb,graphicx,epstopdf}

\newtheorem{theorem}{Theorem}[section]
\newtheorem{corollary}[theorem]{Corollary}
\newtheorem{lemma}[theorem]{Lemma}
\newtheorem{proposition}[theorem]{Proposition}

\newtheorem{definition}[theorem]{Definition}

\newcommand{\dsp}{\displaystyle} 
 
\newcommand{\be}{\begin{equation}}
\newcommand{\ee}{\end{equation}}
\newcommand{\beq}{\begin{eqnarray}}
\newcommand{\eeq}{\end{eqnarray}}
\newcommand{\beqst}{\begin{eqnarray*}}
\newcommand{\eeqst}{\end{eqnarray*}}

\title{Integral Transform Approach to Generalized Tricomi Equations 
}

\author{\scshape Karen Yagdjian}

\begin{document}
\date{}
\maketitle

 \centerline{Department of Mathematics, 
University of Texas-Pan American,}
   \centerline{1201 W.~University Drive, 
Edinburg, TX 78539,  
USA }

\bigskip

\begin{abstract} 
We present some integral transform  that allows to obtain   solutions of the generalized Tricomi 
equation     from  solutions of a simpler equation.
 We used in  \cite{Galstian-Kinoshita-Yagdjian,Galstian-Kinoshita},
\cite{YagTricomi}-\cite{Helsinki_2013} the particular version of this transform
in order to investigate in a unified way several equations such as the linear
and semilinear Tricomi  
equations, Gellerstedt
equation, the wave equation in Einstein-de~Sitter  spacetime, the
wave and the Klein-Gordon equations in the de~Sitter and
anti-de~Sitter spacetimes. 
\medskip 

\noindent {\it MSC:} 35C15; 35Q75; 35Q05; 83F05; 76H05
\medskip 

\noindent {\it Keywords:} Generalized Tricomi Equation; Einstein-de~Sitter  spacetime; Representation of solution
\end{abstract}

\section{Introduction}

In this paper we construct solutions of the generalized Tricomi equation. 
The construction is based on the  extension of the  approach suggested in   \cite{YagTricomi}; this extension allows  to write solution of the generalized Tricomi 
equation  via solutions of a simpler equation.
The particular version of this approach  was used in  \cite{Galstian-Kinoshita-Yagdjian,Galstian-Kinoshita},
\cite{YagTricomi}-\cite{Helsinki_2013}
to investigate in a unified way several equations such as the linear
and semilinear Tricomi  
equations, Gellerstedt
equation, the wave equation in Einstein-de~Sitter  (EdeS) spacetime, the
wave and the Klein-Gordon equations in the de~Sitter and
anti-de~Sitter spacetimes.  The listed  equations play an important
role in the gas dynamics, elementary particle physics, quantum field
theory in curved spaces, and cosmology.
\smallskip

\smallskip

Consider for the smooth function $f=f(x,t)$
the solution $w=w_{A,f}(x,t;b)$ to the  problem
\begin{equation}
\label{0.1JDE}
w_{tt}-  A(x,\partial_x) w =0, \,\, \,
w(x,0;b)= f (x,b), \,\, w_t(x,0;b)=0 ,\,\, t \in [ 0,T_1]\subseteq {\mathbb R} , \,\, x \in \widetilde{\Omega}   \subseteq {\mathbb R}^n,
\end{equation}
with the parameter $b\in I =[t_{in},T] \subseteq {\mathbb R}$, $0 \leq t_{in}< T \leq \infty$,  and with $ 0< T_1 \leq \infty $. Here $\widetilde{\Omega}  $ is a domain in ${\mathbb R}^n $, while  $A(x,\partial_x) $ is the partial differential operator $A(x,\partial_x)=\sum_{|\alpha| \leq m } a_\alpha (x)\partial_x^\alpha $ with smooth coefficients, $a_\alpha \in C^\infty(\widetilde{\Omega}) $.
We are going to present   the integral operator
\begin{equation}
\label{KJDE}
{\mathcal K} [w](x,t)
 =
\int_{  t_{in}}^{t} db
  \int_{ 0}^{ |\phi (t)- \phi (b)|}   K (t;r,b)  w(x,r;b ) \, dr
, \quad x \in \widetilde{\Omega} , \,\, t \in I,
\end{equation}
which maps the function $w=w(x,r;b ) $ into the solution $ u =u(x,t)$ of the generalized Tricomi equation
\begin{equation}
\label{0.3JDE}
u_{tt}-  a^2(t) A(x,\partial_x) u  =f, \qquad x \in \widetilde{\Omega} \,, 
 \,\,   t \in I\,,
\end{equation}
where $a^2(t)=t^\ell $, $\ell \in {\mathbb C} $. In fact, the function $u=u(x,t) $ takes initial values as follows
\[
u(x,t_{in})= 0, \,\, u_t(x,t_{in})=0 ,\quad  x \in \widetilde{\Omega} \,.
\]
In (\ref{KJDE}), $\phi =\phi (t) $ is a distance function produced by $a=a(t) $, that is $\ \phi (t) = \int_{t_{in}}^t a(\tau )\,d\tau $.
Moreover, we also introduce the corresponding operators, which generate solutions of the source-free equation and takes non-vanishing initial values. These operators are constructed in
\cite{YagTricomi} in the case of $\ell>0 $, $A(x,\partial_x)=\Delta  $, $\widetilde{\Omega} = {\mathbb R}^n$, where $\Delta   $ is the Laplace operator on ${\mathbb R}^n $, and, consequently,  equation (\ref{0.1JDE}) is the wave equation.
In the present paper we restrict ourselves to the smooth functions, but it is easily seen that  similar formulas, with the corresponding interpretations, are applicable to the distributions  
as well.

In order to motivate our approach, we consider
the solution $w=w(x,t;b)$ to the Cauchy problem
\begin{equation}
\label{1.7new}
w_{tt}-  \Delta  w =0, \quad   (t,x) \in {\mathbb R}^{1+n},  \qquad
w(x,0;b)= \varphi  (x,b), \,\, w_t(x,0;b)=0 ,\quad x \in {\mathbb R}^n,
\end{equation}
with the parameter $b \in I \subseteq {\mathbb R}$. We denote that solution by $w_\varphi =w_\varphi (x,t;b)$;
if $\varphi $ is independent of the second time variable $b$, then
  we  write simply $w_\varphi (x,t)$.
There are well-known explicit representation formulas for the solution of the   problem (\ref{1.7new}). (See, e.g.,  
\cite{Helgason}, \cite{L-P}, \cite{Shatah}.)

The starting point of the approach suggested in \cite{YagTricomi} is the Duhamel's principle  (see, e.g., \cite[Ch.4]{Shatah}), 
  which has been revised in order to prepare the ground for generalization.
Our {\it first observation} is that the function
\begin{equation}
\label{main}
u(x,t)= \int_{t_{in}}^t \,db \int_{ 0}^{  t-b  } w_f(x,r;b )\,dr\,,
\end{equation}
is the solution of the Cauchy problem
$
 u_{tt}-\Delta u =f(x,t)$ in ${\mathbb R}^{n+1}$, and
$u(x,t_{in})=0$, $u_t(x,t_{in})=0$   in  $\,\, {\mathbb R}^{n}\,,
$
if the function $w_f=w_f(x;t;b ) $ is a solution of the problem  (\ref{1.7new}), where $\varphi =f $.
The  {\it second observation} is that in (\ref{main}) the upper limit $t-b $  of the inner integral is a distance function.
Our {\it third observation}  is that the solution operator ${\mathcal G}\,:\,f \longmapsto u $ can be regarded as a composition of two operators. The first one
\[
{\mathcal W}{\mathcal E}: \,\, f  \longmapsto  w
\]
is a Fourier Integral Operator, which is a solution operator of the Cauchy problem
for wave equation.  The second operator
\[
{\mathcal K}:\,\,w\longmapsto   u
\]
is the integral operator given by (\ref{main}). We regard the variable $b$ in  (\ref{main}) as a ``subsidiary time''. Thus, ${\mathcal G}= {\mathcal K}\circ {\mathcal W}{\mathcal E}$ and we arrive at the  diagram of Figure 1.
\begin{center}
\begin{figure}[hi]
\label{Fig1}
\hspace{1.8cm}\includegraphics[width=0.25\textwidth]{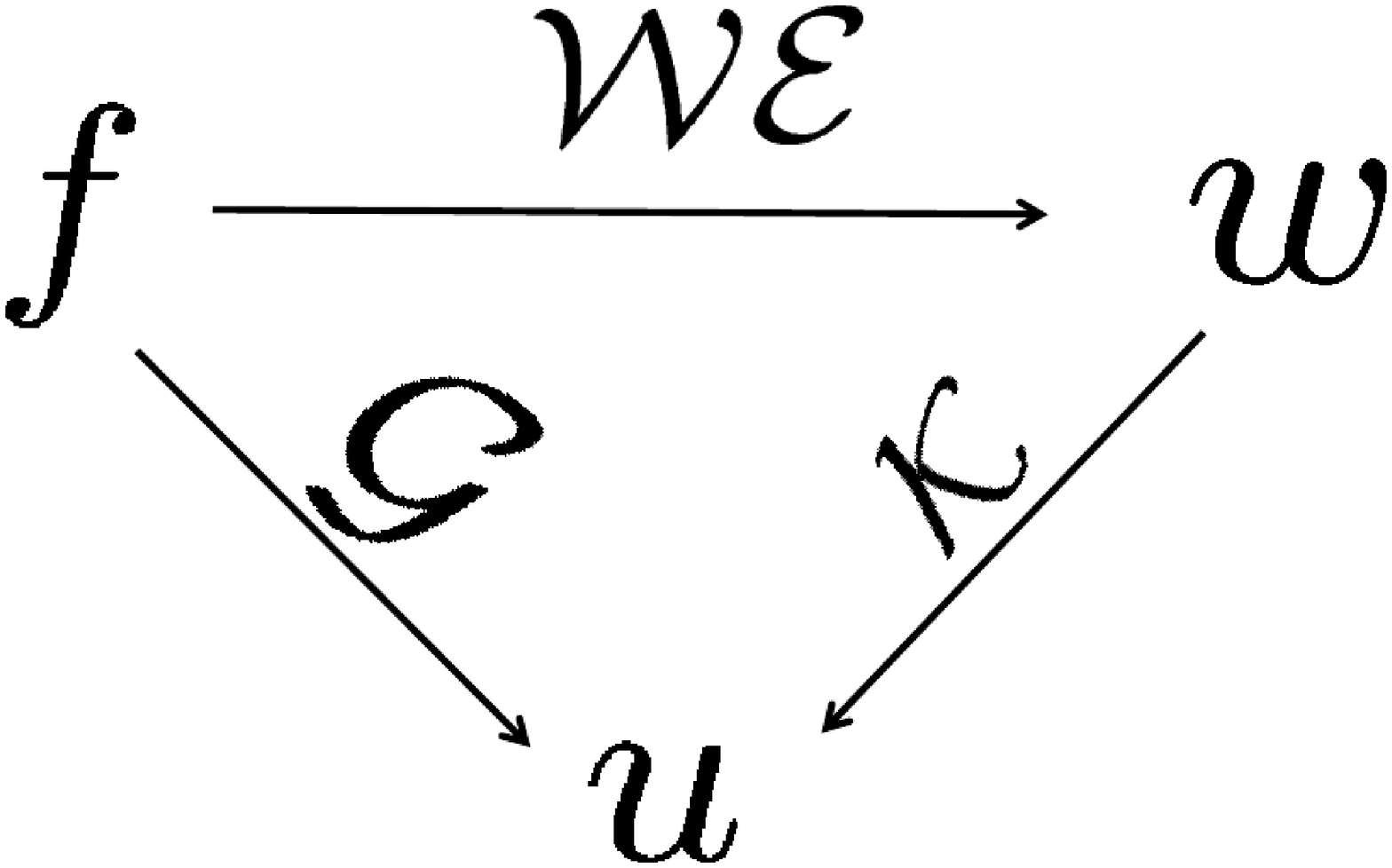}\hspace{2.9cm}
\includegraphics[width=0.25\textwidth]{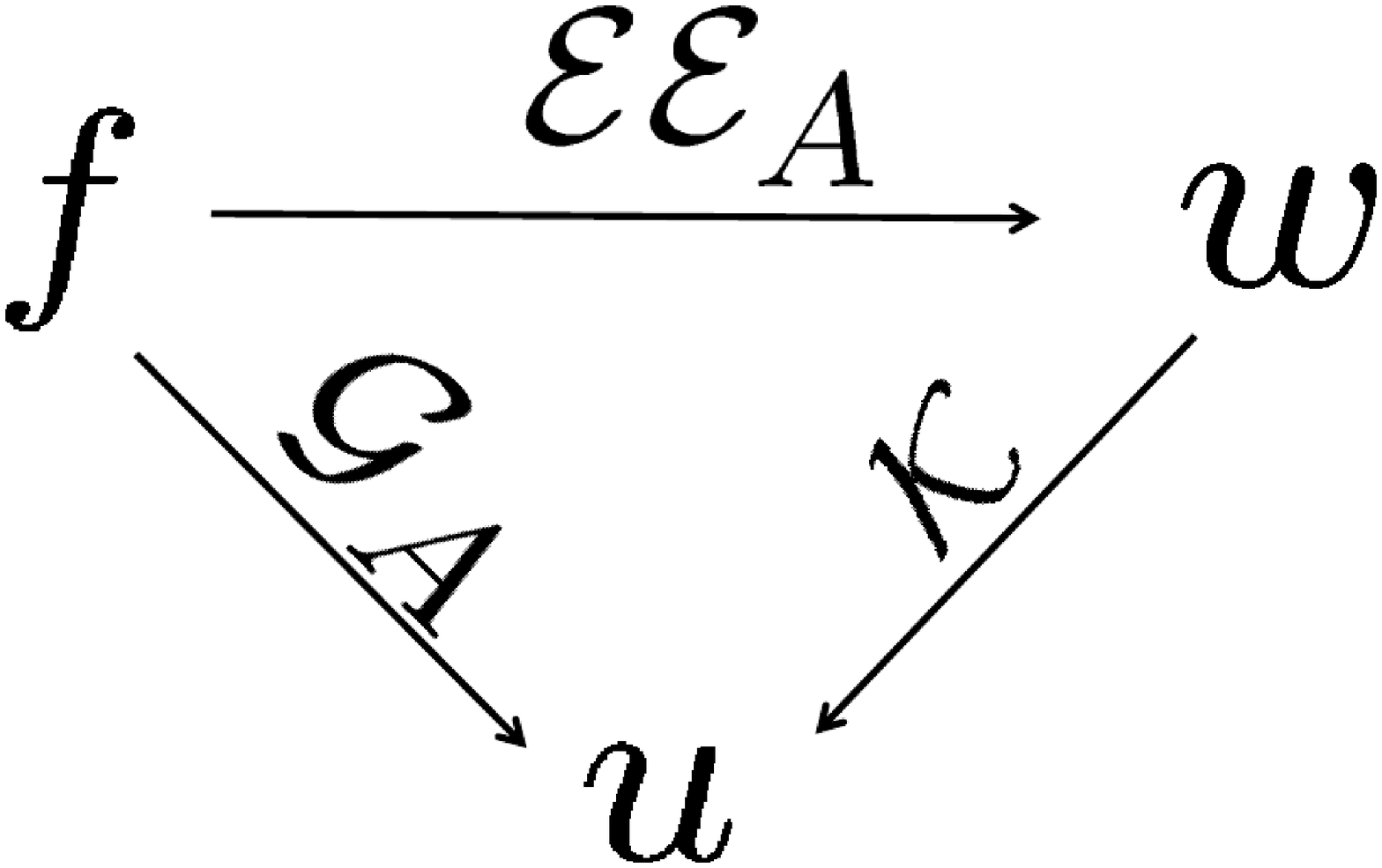}
\caption{\small (a) Case of wave equation $A(x,\partial_x)=\Delta  $ \hspace{0.4cm} { (b)}\, Case of general  $A(x,\partial_x)$}
\end{figure}
\end{center}
\vspace{-0.5cm}
\noindent
Based on the first diagram, we have generated in \cite{yagdjian_Rend_Trieste} a class of operators for which we have  obtained  
 explicit representation formulas for the solutions, and, in particular,  the representations for the fundamental solutions 
of the partial differential operator. Here we emphasize that the integral transform ${\cal K} $ is less singular than 
fundamental solution (Green function) given by $\cal G $ since the operator ${\mathcal W}{\mathcal E} $ takes essential part of singularities.

Our {\it fourth observation}  is that if we plug in (\ref{main}) the solution $w =w (x;t;b ) $   of the Dirichlet problem for the elliptic equation  
$ w_{tt}+ \Delta  w =0$,  $ (t,x) \in {\mathbb R}^{1+n}$,   
$w(x,0;b)= f  (x,b)$,  \, $x \in {\mathbb R}^n$, 
then the integral (\ref{main}) defines the solution $u$ of the equation $
 u_{tt}+\Delta u =f(x,t)+ \int_{t_{in}}^t  w_t (x;0;b )\,db$, such that
$u(x,t_{in})=0$, $u_t(x,t_{in})=0$.   Thus, integral  (\ref{main}), regarded  as an {\sl integral transform}, can be used for elliptic equations as well. 
That integral transform is interesting in its own right. If we want to remove the integral of the right hand side of the equation, then we have to restrict ourselves 
to some particular functions $f$ since the Cauchy problem $w_{tt}+ \Delta  w =0$,    
$w(x,0;b)= f  (x,b)$,   $w_t(x,0;b)= 0$,  
is solvable, even locally, not for every, even smooth, function $f$. 
\smallskip

In the present paper, by varying the first mapping, we extend the class of the equations for which we can obtain  explicit representation formulas for the solutions.  
More precisely, consider the   diagram (b) of Figure~1,
where $w=w_{A,\varphi}(x,t;b)$ is a solution to the  problem (\ref{0.1JDE})
with the parameter $b \in I \subseteq {\mathbb R}$.
If we have a resolving operator of the problem (\ref{0.1JDE}), then,  by applying (\ref{KJDE}),  we can generate   solutions of  another   equation.
Thus, ${\mathcal G}_A = {\mathcal K}\circ {\mathcal E}{\mathcal E}_A $. The new class of equations contains operators with $x$-depending coefficients, and those equations
are not necessarily hyperbolic.
\smallskip

In this paper we restrict ourselves to the   generalized Tricomi  equation,
that is $a(t)=t^{\ell}$, $\ell \in {\mathbb C}$. This class includes, among others, equations of the wave propagating in the so-called Einstein-de~Sitter (EdeS)
universe and in the radiation dominated universe with the spatial slices
of the constant curvature.  We believe that the integral transform and the representation formulas for the solutions that we derive in this article fill up the gap in the literature on that topic.
\smallskip

The transform linking to the generalized Tricomi operator  is generated by the kernel 
\begin{eqnarray}
\label{E}
K (t;r,b)
& = &
E(r,t;b;\gamma  ):= c_\ell \left(  (\phi (t)  + \phi (b))^2  -r^2 \right)^{-\gamma }  F \left(\gamma , \gamma ;1; \frac{(\phi (t)  - \phi
(b))^2 - r^2} {(\phi (t)  + \phi (b))^2 - r^2}  \right) ,
\end{eqnarray}
with   
the distance function $\phi  =\phi (t)  $ and the numbers $\gamma  $, $c_\ell $ defined as follows
\begin{eqnarray}
\label{phi_Tric}
\phi (t)=  \frac{2}{\ell+2}t^{\frac{\ell+2}{2} } ,\quad \gamma := \frac{\ell}{2(\ell+2)}
,\quad  \ell \in {\mathbb C}\setminus \{-2\},\quad c_\ell = \left(\frac{ \ell+2 }{4}\right)^{-\frac{ \ell }{ \ell+2 }} ,
\end{eqnarray}
while $F\big(a, b;c; \zeta \big) $ is the Gauss's hypergeometric function. Here $t_{in}=0$.

\smallskip

 According to Theorem~\ref{Teqforkernel}, the function $  E(r,t;b;\gamma  )$ solves  the following  Tricomi-type equation:
 \begin{eqnarray}
 \label{Eeq}
&  &
E_{tt} (r,t;b;\gamma  ) - t^{\ell} E_{rr}(r,t;b;\gamma  ) =0\,, \qquad 0< b<t
\,.
\end{eqnarray}
The proof of Theorem~\ref{Teqforkernel}, which  is given in Section~\ref{S2}, is straightforward. In fact, that proof is applicable   to the different   distance functions $\phi  =\phi (t)  $,
see, for instance, \cite{Yag_Arx}, where the case of $ a (t) =e^{-t}$ is discussed.
\medskip

There are four important examples of the equations which are amenable to the integral transform approach, when $\ell=3,1, -1,-4/3$; those are the small disturbance equations for the perturbation velocity potential of a two-dimensional 
near sonic uniform flow of dense gases in
a physical plane (see, e.g., \cite{Kluwick}, \cite{T-C}), the Tricomi equation (see, e.g., \cite{Barros-Neto-Gelfand-I,Bers,Cole-Cook,Frankl,Germain,Germain-Bader,Kim,Lau-Price,Lupo-Payne,LupoPayne2002,Lupo_Payne,Morawetz,Morawetz2004,Nocilla,Payne,Smirnov,Tricomi} and bibliography therein), the equation of the waves in  the radiation dominated universe (see, e.g., \cite{Ellis,Gron-Hervik} and bibliography therein) and in the EdeS spacetime (see, e.g., \cite{Ellis,Gron-Hervik,Hawking,Ohanian-Ruffini} and bibliography therein), respectively. 

To introduce the integral transform we need some special geometric structure of the domains of functions. 
\begin{definition}
The set $\Omega \subseteq \overline{{\mathbb R}^{n+1 }_+ }$ is said to be backward time line-connected to $t=0$, if for every point $(x,t) \in \Omega $ the 
line segment 
$\{(x,s)\,|\, s \in (0,t]\, \} $ is also in $\Omega  $; that is $\{(x,s)\,|\, s \in (0,t] \} \subseteq  \Omega$.
\end{definition}
Henceforth we just write "backward time connected"  for such sets. Similarly, if $\Omega \subseteq [0,T]\times {\mathbb R}^{n} $, $T>0$, then one can define a forward time line-connected to $t=T$ set. The union and intersection of the backward time connected sets is also a backward time connected. The interior and the closure of the backward time connected set are also a backward time connected sets. For every set   there exists its minimal backward time connected   covering. 
The domain of the dependence for the wave equation is backward time connected, while domain of influence is forward time connected.

\begin{definition}
Let $\phi  \in C  (\overline{{\mathbb R}_+ }) $ be non-negative increasing function and $\Omega  $ be a  backward time connected set. The  backward time connected set $\Omega_\phi  \subseteq \overline{{\mathbb R}^{n+1 }_+ }$
defined by
\[
\Omega_\phi := \bigcup _{(x,t) \in \Omega } \left\{ (x,\tau ) \,|\, \tau \in (0,\phi (t)]\, \right\}
\]
is said to be a $\phi$-image of $\Omega  $. 
\end{definition}

On Figures 2,3 we illustrate the dependence domains for hyperbolic equations with  $\ell=0,-\frac{4}{3},1,3  $ and $A(x, \partial_x )= \Delta  $. 
\begin{center}
\begin{figure}[hi]
\label{Fig2}
\includegraphics[width=0.24\textwidth]{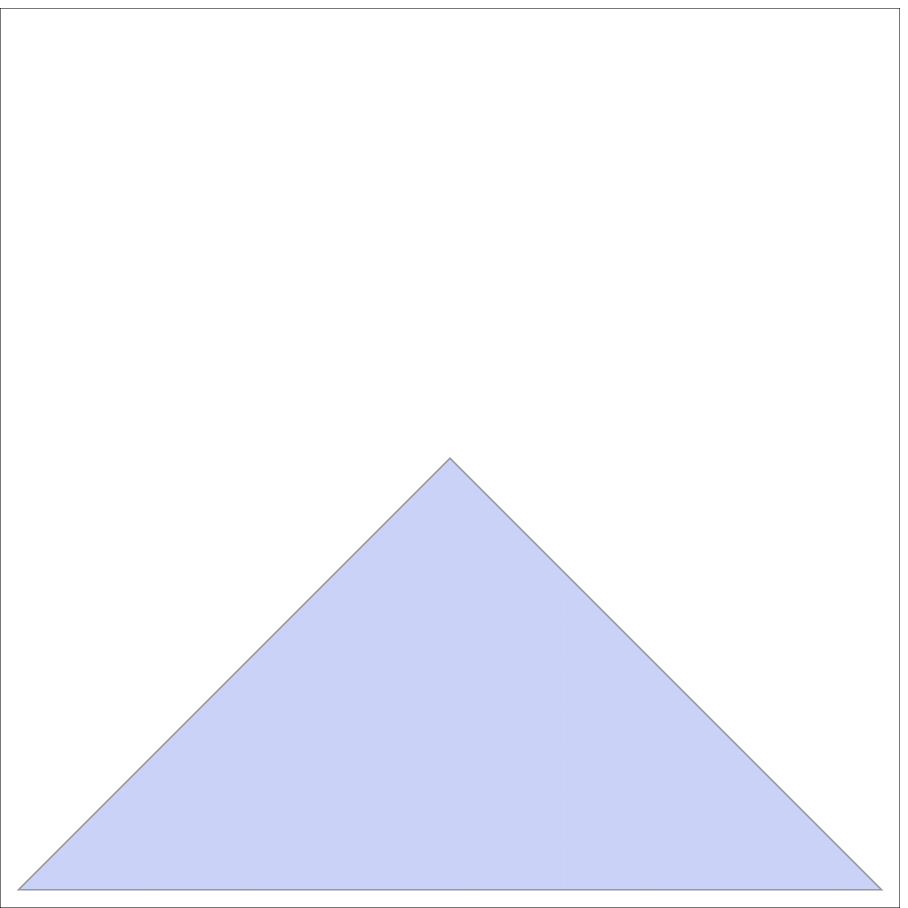}
\includegraphics[width=0.24\textwidth]{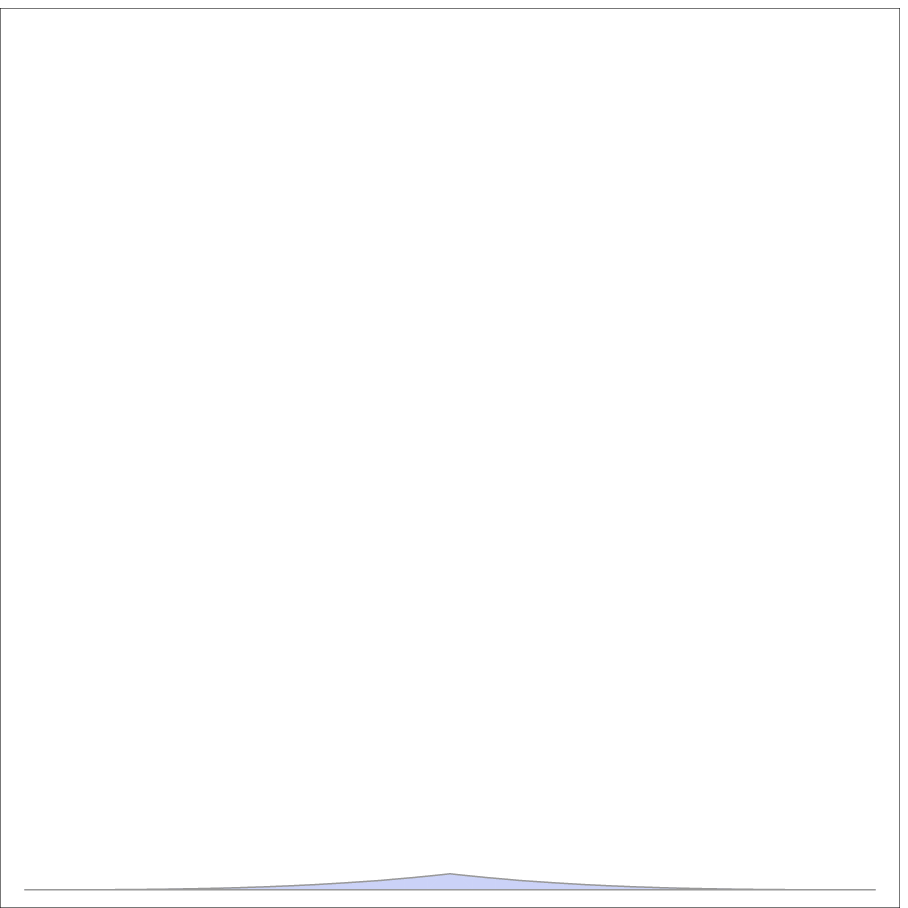}\includegraphics[width=0.24\textwidth]{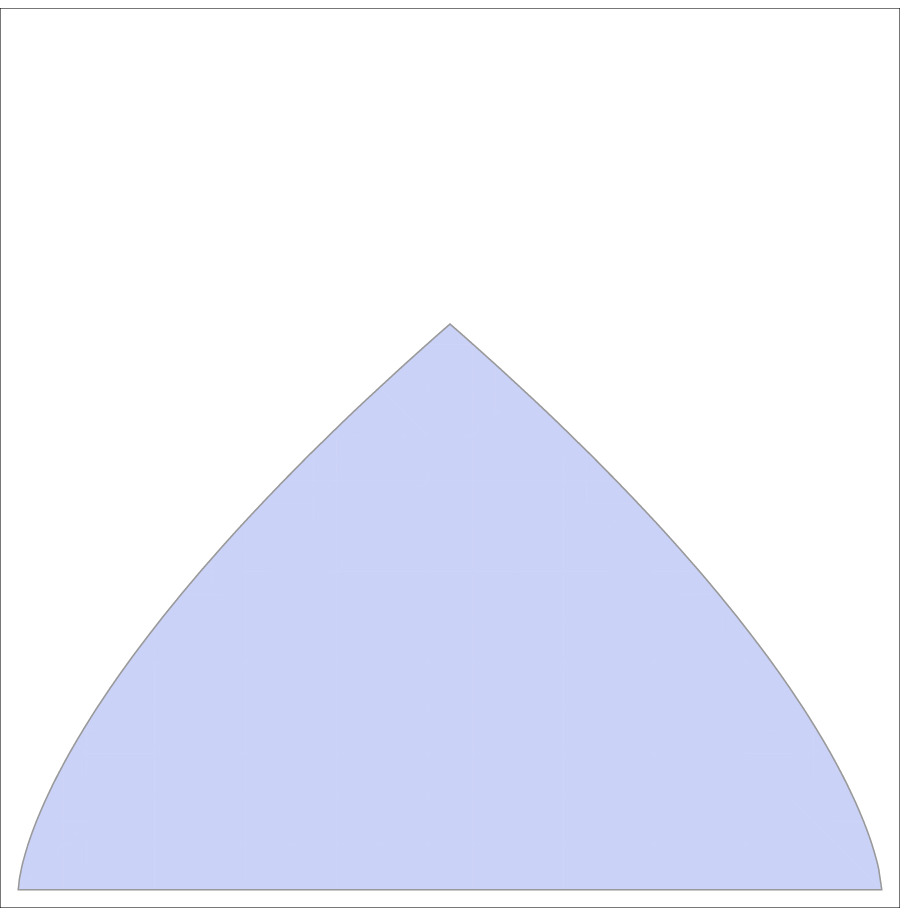}\includegraphics[width=0.24\textwidth]{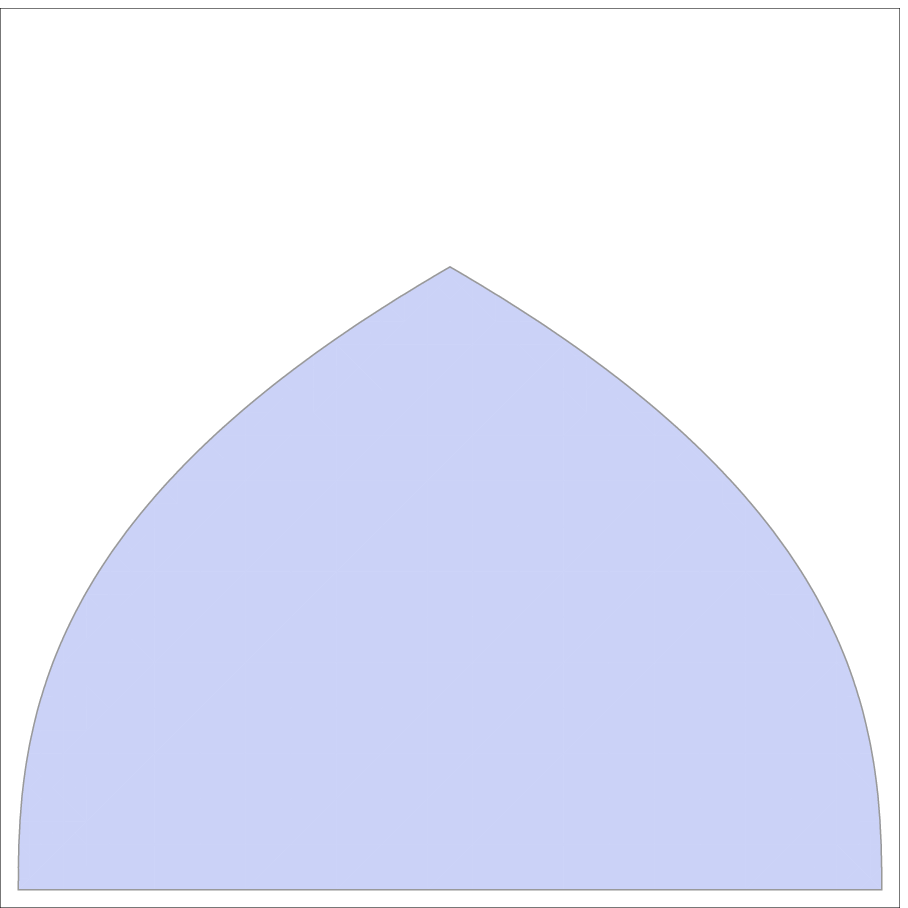}
\caption{$\phi (t)= (2/(\ell+2))t^{(\ell+2)/2}  $,  $ \,\,\phi (t)+ | x | \leq 1$, $x\in [-1,1]$, $\,\,t \in [0,2]$, $\,\, \ell=0,-\frac{4}{3},1,3 $ \hspace{0.4cm}  \,  }
\end{figure}
\end{center}
\begin{center}
\begin{figure}[hi]
\label{Fig3}
\includegraphics[width=0.24\textwidth]{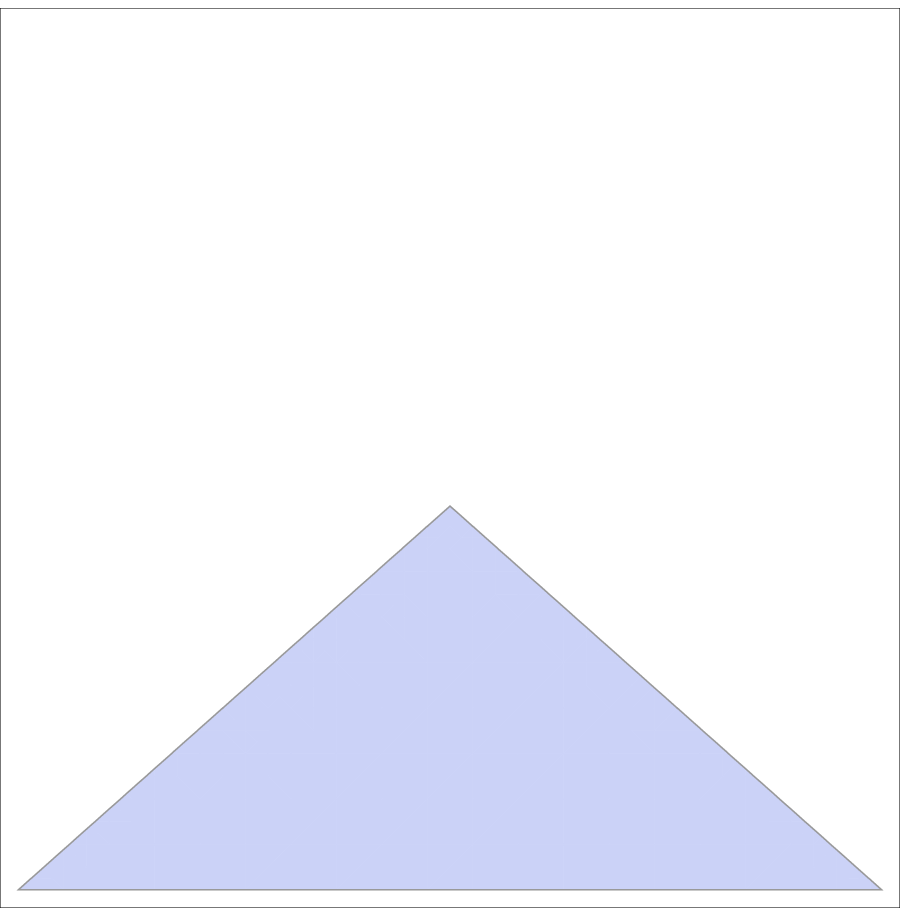}
\includegraphics[width=0.24\textwidth]{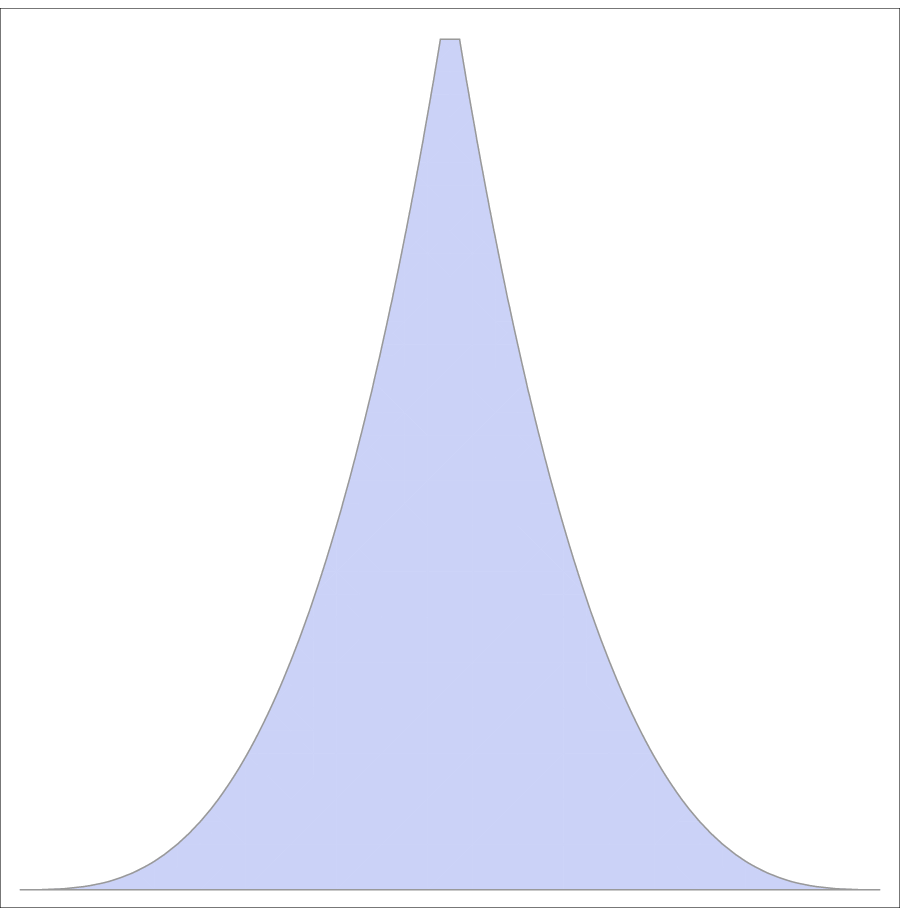}\includegraphics[width=0.24\textwidth]{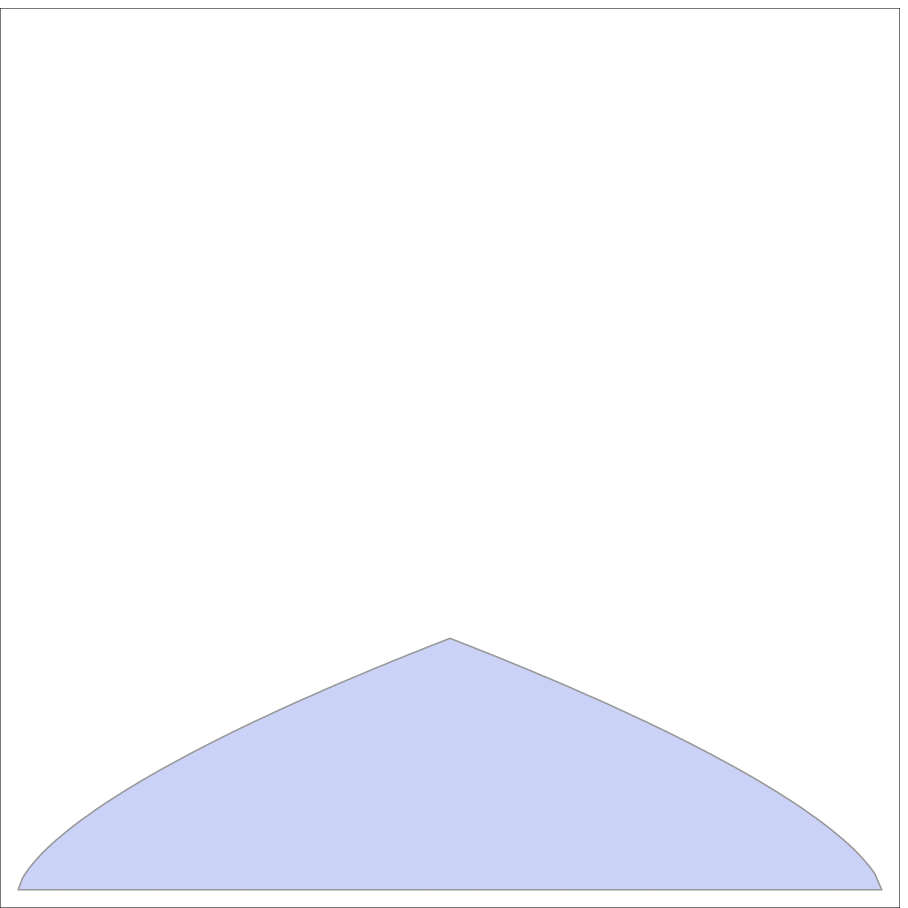}\includegraphics[width=0.24\textwidth]{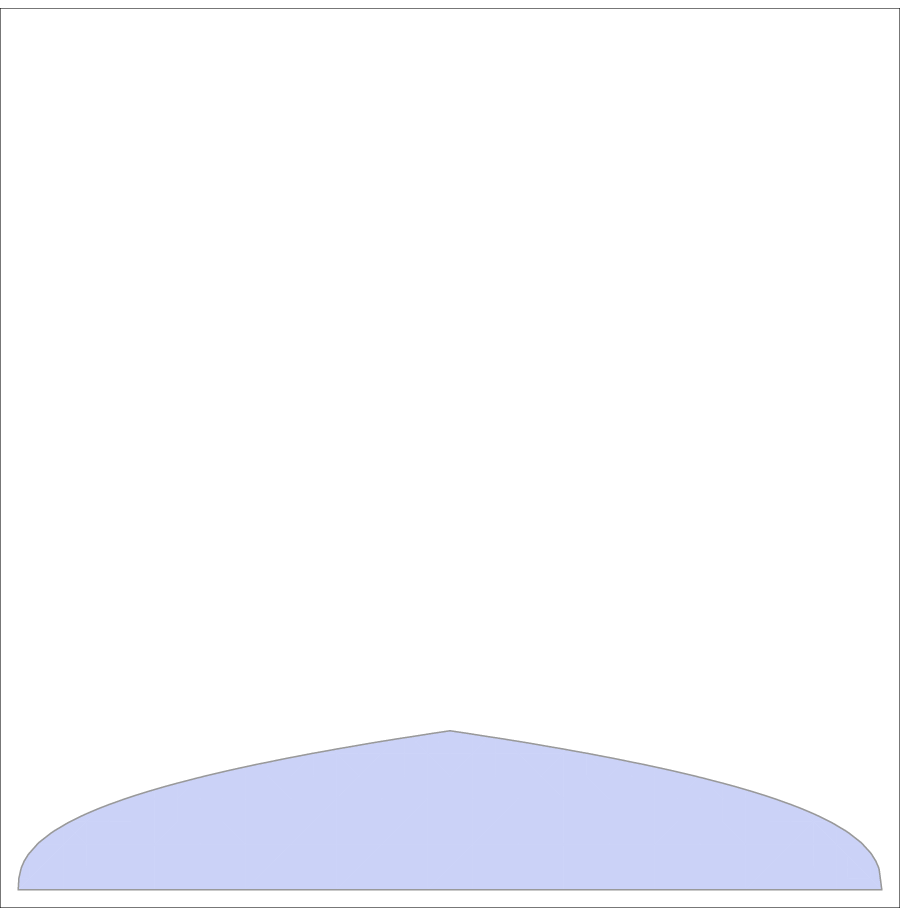}
\caption{$\phi (t)= (2/(\ell+2))t^{(\ell+2)/2}  $,  $\,\,\phi (t)+ | x | \leq 8$, $x\in [-8,8]$, $\,\,t \in [0,18]$, $\,\, \ell=0,-\frac{4}{3},1,3$ \hspace{0.4cm}  \,  }
\end{figure}
\end{center} 
 Figure 4 illustrates the part of the  domain in the hyperbolic region of the Tricomi problem (see, e.g., \cite{Lupo-Payne} and references therein) that has the form 
$ \Omega := \{(x,t) \,|\,|x| < x_0 -\frac{2}{3}t^{\frac{3}{2} },\,\,  -x_0 \leq x\leq x_0,\,\, t>0\}$ for $x_0=1/2$ and $x_0=100$.  Corresponding $ \phi $-images of $\Omega $ are $ \Omega_\phi  := \{(x,t) \,|\,|x| < x_0 -\left(\frac{2}{3} \right)^{5/2} t^{\frac{9}{4}},\,\,  -x_0 \leq x\leq x_0\}$ with $x_0=1/2$ and $x_0=100$, respectively. 

\begin{center}
\begin{figure}[hi] 
\label{Fig4}
\includegraphics[width=0.24\textwidth]{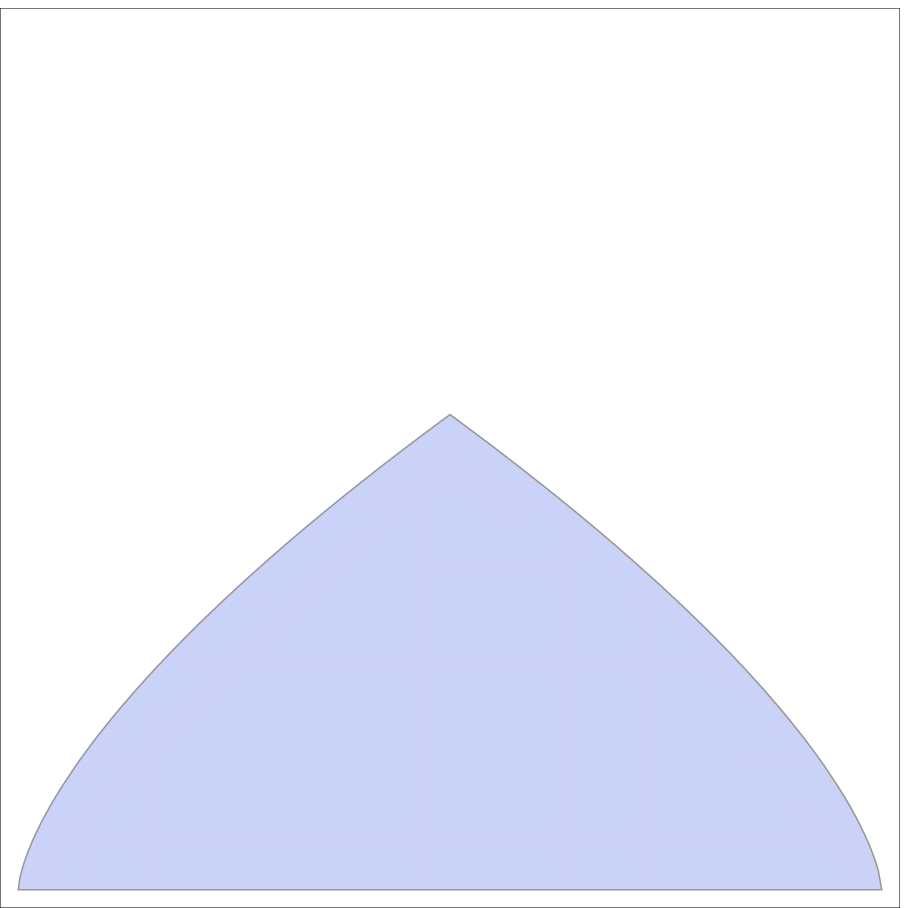}  
\includegraphics[width=0.24\textwidth]{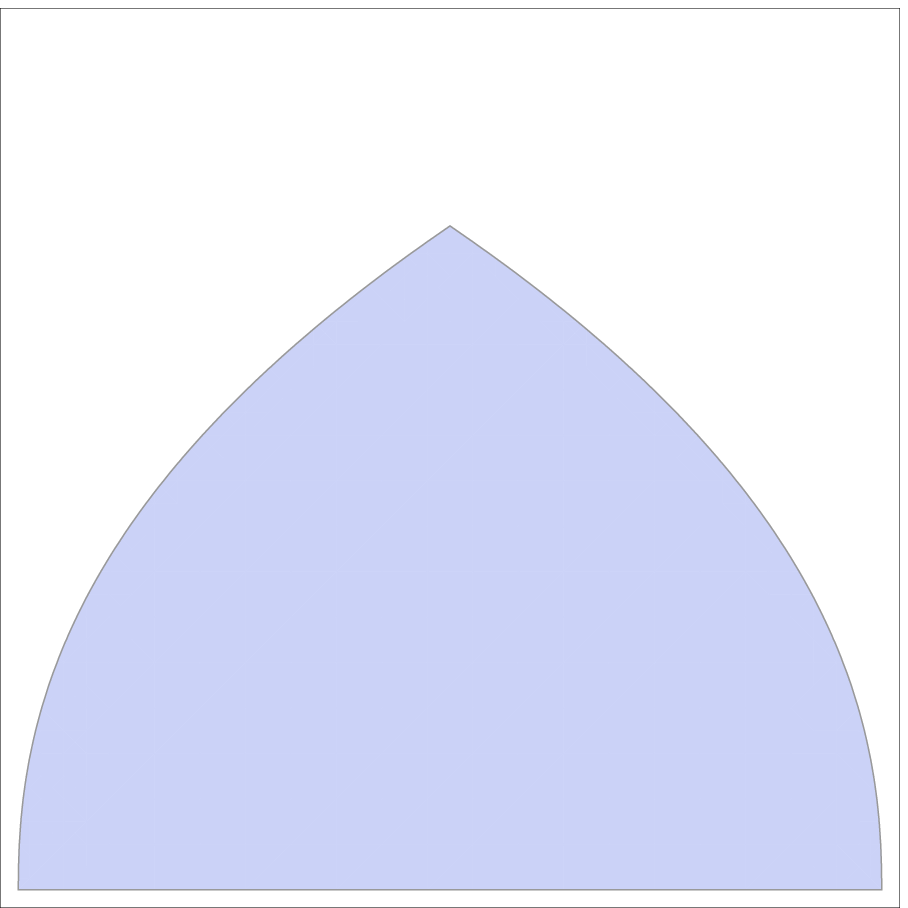} \includegraphics[width=0.24\textwidth]{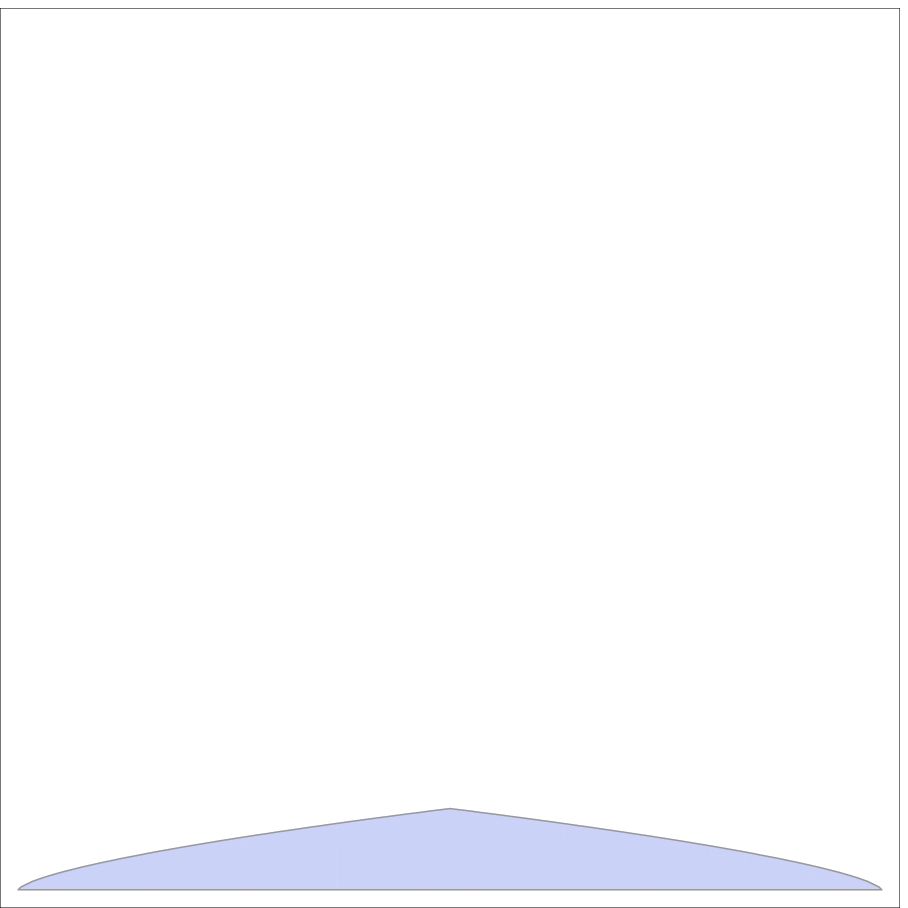}\includegraphics[width=0.24\textwidth]{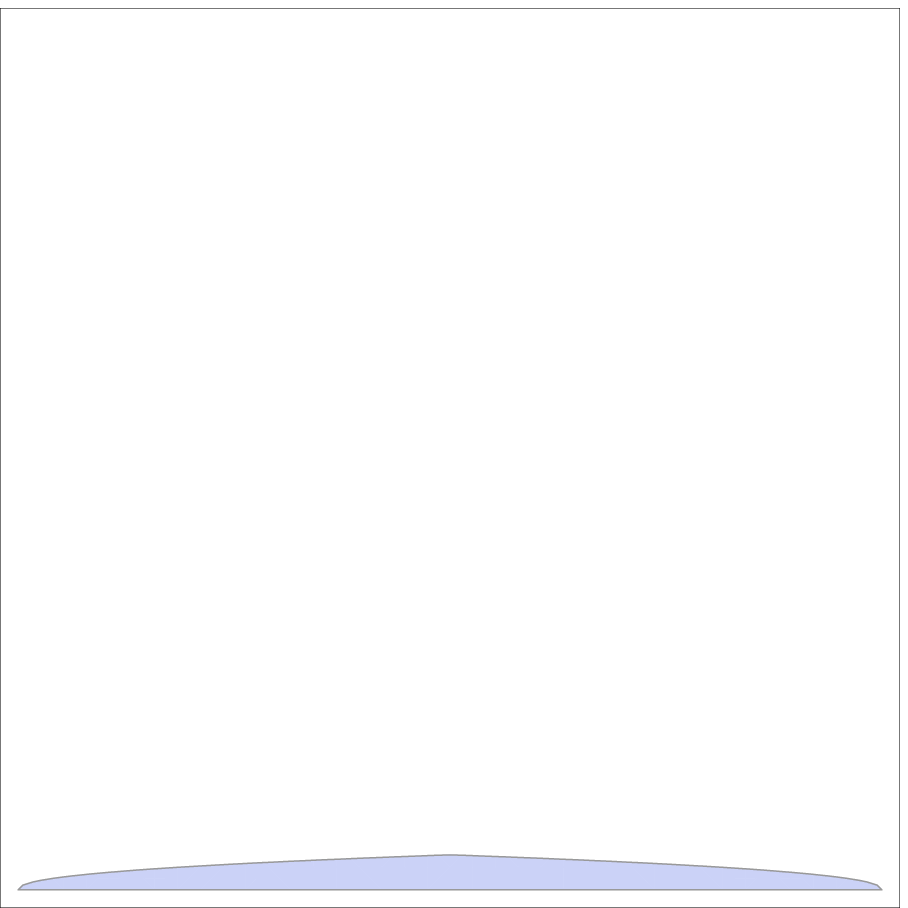} 
\caption{\small $ \Omega $ and $ \Omega_\phi  $ with $x_0=1/2$;  \hspace{2.5cm}   $ \Omega $ and $ \Omega_\phi  $ with $x_0=100$}
\end{figure}
\end{center}

The next theorem describes the main property of the integral transform (\ref{KJDE}). We denote $ U_{x_0}$  a neighborhood of the point $x_0 \in {\mathbb R}^n$.
\begin{theorem}  
\label{T1.8}
Let $f=f(x,t) $  be a  function defined in the backward time connected domain $\Omega  $.
Suppose that for a given $(x_0,t_0) \in \Omega  $ the  function  $w(x,r;b ) \in C_{x,r,t} ^{m,2,0} (U_{x_0} \times [0,  \phi (t_0)]\times [0,t_0])$     
solves the problem 
\begin{eqnarray}
\label{8}
&  &
w_{rr}-  A(x ,\partial_x)  w =0  \quad \mbox{in}\quad U_{x_0}  \quad \mbox{  for all }
\quad  r \in (0, |\phi (t_0)- \phi (b)|)\,,\\
\label{9}
&  &
 w(x ,0;b)= f(x ,b) \quad  \mbox{in}\quad U_{x_0} \quad \mbox{  for all }
\quad  b \in (0,t_0)\,.
\end{eqnarray} 

Then 
for $\ell >-2$   the function
\[
u(x,t)
  =  
c_\ell \int_{ 0}^{t} db
  \int_{ 0}^{ |\phi (t)- \phi (b)|}  \left(  (\phi (t)  + \phi (b))^2  -r^2 \right)^{-\gamma }  
  F \left(\gamma , \gamma ;1; \frac{(\phi (t)  - \phi (b))^2 - r^2} {(\phi (t)  + \phi (b))^2 - r^2}  \right)  w(x,r;b )\,
dr,  
\]
defined  in the past  $ U_{x_0}\times  [0,t_0) $ of   $ U_{x_0}\times \{t_0\} $, is continuous in $ U_{x_0}\times  [0,t_0) $  and it  satisfies  the  equation
\begin{eqnarray}
\label{Tric_eq}
 u_{tt}-t^{\ell}A(x ,\partial_x) u 
& = &
f(x ,t) \nonumber \\
&  &
+ c_\ell(\phi' (t))^2 \int_{ 0}^{t} 
 \left(  \phi (t)  + \phi (b)  \right)^{-2\gamma }  F \left(\gamma , \gamma ;1; \frac{(\phi (t)  - \phi
(b))^2  } {(\phi (t)  + \phi (b))^2  }  \right)  w_{r} (x ,0;b )\, db \,,
\end{eqnarray}
in the sense of distributions   ${\mathcal D}'(U_{x_0}\times  (0,t_0) )  $.  The function $u(x,t) $  takes the vanishing initial value $u(x,0)=0 $ for all $x \in U_{x_0}$.
Moreover, if, additionally, $\ell <4 $, then $u_t $ is continuous in $ U_{x_0}\times  [0,t_0) $  and  $u_t(x,0)=0 $ for all $x \in U_{x_0}$.
\end{theorem}

We stress here that the integral transform $w\longmapsto u $ is point-wise in $x$ and non-local in time. 
Let $\pi_x  $ be a projection $\pi_x: \Omega \longrightarrow {\mathbb R}^n $ of the backward time connected domain $\Omega  $, and denote $ \widetilde{\Omega }:= \pi_x (\Omega)$.

\begin{corollary} 
\label{C1.4}
Let $f=f(x,t) $  be a  function defined in the backward time connected domain $\Omega  $. 
Suppose that   the function $w(x,r;b ) \in C_{x,r,t}^{m,2,0}$    
satisfies   
\begin{eqnarray*}
&  &
w_{rr}-  A(x ,\partial_x)  w =0 \quad  \mbox{  for all }
\quad  (x ,r) \in \Omega_\phi  \quad  \mbox{  and }
\quad (x ,b) \in \Omega  \,,\\
&  &
 w(x ,0;b)= f(x ,b) \quad  \mbox{  for all }
\quad  (x ,b) \in \Omega   \,.
\end{eqnarray*} 

Then 
for $\ell  >-2$   the function
\[
u(x,t)
  =  
c_\ell \int_{ 0}^{t} db
  \int_{ 0}^{ |\phi (t)- \phi (b)|}  \left(  (\phi (t)  + \phi (b))^2  -r^2 \right)^{-\gamma }  F \left(\gamma , \gamma ;1; \frac{(\phi (t)  - \phi
(b))^2 - r^2} {(\phi (t)  + \phi (b))^2 - r^2}  \right)  w(x,r;b )\,
dr, 
\]
defined on    $\Omega $, is continuous and satisfies  the  equation
\begin{eqnarray}
\label{Tric_eqC}
 u_{tt}-t^{\ell}A(x ,\partial_x) u 
& = &
f(x ,t) \nonumber \\
&  &
+ c_\ell(\phi' (t))^2 \int_{ 0}^{t} 
 \left( \phi (t)  + \phi (b) \right)^{-2\gamma }  F \left(\gamma , \gamma ;1; \frac{(\phi (t)  - \phi
(b))^2  } {(\phi (t)  + \phi (b))^2  }  \right)  w_{r} (x ,0;b )\, db \,,
\end{eqnarray}
in the sense of distributions   ${\mathcal D}'(\Omega  )  $.   The function $u$ takes the vanishing initial value
$ u(x ,0)=0 $      for all $x \in \widetilde{\Omega }$. If, additionally, $\ell <4 $, then $u_t $ is continuous in $ \Omega  $  and  $u_t(x,0)=0 $ for all $x \in \widetilde{\Omega }$.
\end{corollary}

If the initial value problem for the operator \, $\partial_t^2- A(x,\partial_x) $ \, admits two initial conditions, 
then we can eliminate the function \,$ w_r$ \, from the right-hand side of equations (\ref{Tric_eq}) and (\ref{Tric_eqC}).

\begin{theorem} 
\label{T1.8b}
Let $f=f(x,t) $  be a smooth function defined in the backward time connected domain $\Omega  $. 
Suppose that   the smooth function $w(x,r;b ) \in C_{x,r,t}^{m,2,0}$    
satisfies   
\begin{eqnarray*}
&  &
w_{rr}-  A(x ,\partial_x)  w =0 \quad  \mbox{  for all }
\quad  (x ,r) \in \Omega_\phi  \quad  \mbox{  and for all }
\quad (x ,b) \in \Omega  \,,\\
&  &
 w(x ,0;b)= f(x ,b),  \quad   w_r(x ,0;b)= 0 \quad \mbox{  for all }
\quad  (x ,b) \in \Omega   \,.
\end{eqnarray*} 

Then 
for $\ell  >-2$   the function
\[
u(x,t)
  =  
c_\ell \int_{ 0}^{t} db
  \int_{ 0}^{ |\phi (t)- \phi (b)|}  \left(  (\phi (t)  + \phi (b))^2  -r^2 \right)^{-\gamma }  F \left(\gamma , \gamma ;1; \frac{(\phi (t)  - \phi
(b))^2 - r^2} {(\phi (t)  + \phi (b))^2 - r^2}  \right)  w(x,r;b )\,
dr, 
\]
defined on    $\Omega $, is continuous and satisfies  the  equation
\begin{eqnarray}
\label{Tric_eqCd}
 u_{tt}-t^{\ell}A(x ,\partial_x) u 
& = &
f(x ,t) \,,
\end{eqnarray}
in the sense of distributions   ${\mathcal D}'(\Omega  )  $.   The function $u$ takes the vanishing initial value
$ u(x ,0)=0 $      for all $x \in \widetilde{\Omega }$. If, additionally, $\ell <4 $, then $u_t $ is continuous in $ \Omega  $  and  $u_t(x,0)=0 $ for all $x \in \widetilde{\Omega }$.
\end{theorem}
\smallskip

For instance, the Cauchy problem for the second order strictly hyperbolic equation admits two initial conditions.
We remind here that for the weakly hyperbolic operators $\partial_{t}^2- \sum_{|\alpha| \leq 2}a_\alpha (x) \partial_x^\alpha $, which satisfy the Levi conditions (see, e.g., \cite{YagBook}), 
the Cauchy problem can be  solved for the smooth initial data. If $m=1$ then the problem with two initial conditions can be solved in Gevery spaces. (See, e.g., \cite{YagBook}.)
The case of $m>2$ covers the beam equation   and hyperbolic in the sense of Petrowski ($p$-evolution ) equations.
On the other hand, the Cauchy-Kowalewski theorem guarantees 
solvability of the   problem in the real analytic functions category for the   partial differential equation  (\ref{Tric_eqCd}) with any positive $\ell$ and $m=2$. 
Furthermore, the operator $A(x,\partial_x)=\sum_{|\alpha| \leq 2}a_\alpha (x) \partial_x^\alpha$ can be replaced with an abstract operator $A$
acting on some linear topological space of functions.\\
\noindent
{\bf Example 1.} Consider equations of the gas dynamics. (a) For the Tricomi equation in the {\it hyperbolic domain}, 
\begin{equation}
\label{5}
 u_{tt}-t\Delta  u=f(x,t),
\end{equation}
$\phi (t)= \frac{2}{3} t^\frac{3}{2}$  and $A(x,\partial_x)= \Delta$. Then for  every $ f \in C({\mathbb R}^n\times [0,T])$ we can solve 
 the Cauchy problem for the 
wave equation  
\[
w_{tt}-\Delta w= 0, \quad w(x,0;b)= f(x,b)\,,\quad w_t(x,0;b)=0\,, \quad x \in {\mathbb R}^n, \quad t \in \left[0, \frac{2}{3} T^\frac{3}{2}\right]
\]
in  ${\mathbb R}^n\times\left[0, \frac{2}{3} T^\frac{3}{2}\right] \times [0,T]$. (For the explicit formula see, e.g., (\ref{vphiodd}),(\ref{vphieven}).)
The solution to the Cauchy problem for (\ref{5}) with vanishing initial data 
is given as follows
\begin{eqnarray*}
u(x,t)
& = &
 3^{-1/3}2^{2/3}\int_{ 0}^{t} db
  \int_{ 0}^{\frac{2}{3} t^\frac{3}{2}- \frac{2}{3} b^\frac{3}{2} }  \left(  \left(\frac{2}{3} t^\frac{3}{2}+ \frac{2}{3} b^\frac{3}{2} \right)^2  -r^2 \right)^{-\frac{1}{6} } \\
&  & \times F \left(\frac{1}{6} , \frac{1}{6} ;1; \frac{\left(\frac{2}{3} t^\frac{3}{2}- \frac{2}{3} b^\frac{3}{2} \right)^2 - r^2} {\left(\frac{2}{3} t^\frac{3}{2}+\frac{2}{3} b^\frac{3}{2} \right)^2 - r^2}  \right)  w(x,r;b )
dr, \quad    t \in [0,T]\,.
\end{eqnarray*}
For the Tricomi equation in the {\it elliptic domain}, 
\begin{equation}
\label{5ell}
 u_{tt}+t\Delta  u=f(x,t), \quad t>0,
\end{equation}
we have $A(x,\partial_x)=-\Delta  $ and, since the Cauchy problem is not well posed,   Theorem~\ref{T1.8} gives representation of the solutions only for some specific functions $f$.

\noindent
(b) The small disturbance equation for the perturbation velocity potential of a two-dimensional 
near sonic uniform flow of dense gases in
a physical plane, has been derived by 
Kluwick \cite{Kluwick}, Tarkenton and Cramer \cite{T-C}. It leads to the equation 
\begin{equation}
\label{Kluwick}
 u_{tt}-t^3\Delta  u=f(x,t),
\end{equation}
with $\ell=3$ and   $\phi (t)= \frac{2}{5} t^\frac{5}{2}$, and $A(x,\partial_x)= \Delta$. 
The solution to the Cauchy problem for (\ref{5}) with vanishing initial data 
is given as follows
\begin{eqnarray*}
u(x,t)
& = &
 \frac{3}{10}\int_{ 0}^{t} db
  \int_{ 0}^{\frac{2}{5} t^\frac{5}{2}- \frac{2}{5} b^\frac{5}{2} }  \left(  \left(\frac{2}{5} t^\frac{5}{2}+ \frac{2}{5}b^\frac{5}{2} \right)^2  -r^2 \right)
^{- \frac{3}{10} } \\
&  & \times F \left(\frac{3}{10} , \frac{3}{10} ;1; 
\frac{\left(\frac{2}{5} t^\frac{5}{2}- \frac{2}{5} b^\frac{5}{2} \right)^2 - r^2} {\left(\frac{2}{5} t^\frac{5}{2}+\frac{2}{5} b^\frac{5}{2} \right)^2 - r^2}  \right)  w(x,r;b )
dr, \quad  t>0\,.
\end{eqnarray*}
\noindent
{\bf Example 2.} Consider the wave equation in the Einstein-de Sitter (EdeS) spacetime with hyperbolic spatial slices. The metric of the Einstein~\&~de~Sitter universe (EdeS universe) is a particular member of the
Friedmann-Robertson-Walker metrics
\begin{eqnarray}
\label{FLRW}
ds^2= -dt^2+ a_{sc}^2(t)\left[ \frac{dr^2}{1-Kr^2}   + r^2  d\Omega ^2 \right]\,,
\end{eqnarray}
where $K=-1,0$, or $+1$, for a hyperbolic, flat or spherical spatial geometry, respectively. 
For the EdeS the scale factor is $a_{sc}(t)=t^{2/3} $.
The covariant d'Alambert's operator,  
\[
\square_g \psi = \frac{1}{\sqrt{|g|}}\frac{\partial }{\partial x^i}\left( \sqrt{|g|} g^{ik} \frac{\partial \psi }{\partial x^k} \right)\,,
\]
in the spherical coordinates is 
\begin{eqnarray*}
\square_{EdeS} \psi
& = &
-  \left( \frac{\partial }{\partial x^0} \right)^2  \psi 
- \frac{2 }{t  } \left(       \frac{\partial \psi }{\partial x^0} \right)\psi 
+  t^{-\frac{4}{3}}\frac{\sqrt{1-Kr^2}}{ r^2 }\frac{\partial }{\partial r}\left(    r^2   \sqrt{1-Kr^2} \frac{\partial \psi }{\partial  r} \right)\\
&  &
+  t^{-\frac{4}{3}}\frac{1 }{ r^2 \sin  \theta  }\frac{\partial }{\partial \theta }\left(   \sin  \theta    \frac{\partial \psi }{\partial \theta } \right)
+ t^{-\frac{4}{3}} \frac{1}{ r^2 \sin^2  \theta } \left( \frac{\partial }{\partial  \phi} \right)^2\psi \,.
\end{eqnarray*}
The change $\psi =t^{-1}u $ of the unknown function leads the equation \, $ \square_{EdeS} \psi =g$ \, to the equation
\begin{eqnarray*}
&  &
  u_{tt} 
-  t^{-4/3}A(x,\partial_x) u =f\,,
\end{eqnarray*}
where
\begin{equation}
A(x,\partial_x) u
 =  
   \frac{\sqrt{1-Kr^2}}{ r^2 }\frac{\partial }{\partial r}\left(    r^2   \sqrt{1-Kr^2} \frac{\partial u }{\partial  r} \right)
+   \frac{1 }{ r^2 \sin  \theta  }\frac{\partial }{\partial \theta }\left(   \sin  \theta    \frac{\partial u }{\partial \theta } \right)
+ \frac{1}{ r^2 \sin^2  \theta } \left( \frac{\partial }{\partial  \phi} \right)^2 u\,.
\end{equation}
The spatial part $X$ of the spacetime (\ref{FLRW}) has a constant curvature $6K$. Operator $ A(x,\partial_x)$ is the Laplace-Beltrami operator  on  $X$.  
The explicit formulas for the solutions of the Cauchy problem for the wave operator on the 
spaces with constant negative curvature  are known, see, for instance, \cite{Helgason,L-P}. Thus, Theorem~\ref{T1.8b}
gives explicit representation for the solution of the Cauchy problem with vanishing initial data for the wave equation in the  EdeS  spacetime with the negative constant curvature $K<0$. In order to keep down the length of this paper, we postpone applications of Theorem~\ref{T1.8b} to the derivation of the Strichartz estimates
and to global well-posedness of the nonlinear generalized Tricomi  equation in the metric (\ref{FLRW}). We note here that $\gamma =-1 $  for the  metric (\ref{FLRW})  (see \cite{Galstian-Kinoshita-Yagdjian})
that makes 
the hypergeometric function polynomial. Moreover, the hypergeometric function is polynomial for $\gamma =-1,-2,\ldots$  that essentially simplifies calculations.

\medskip

The next theorem   represents the integral transforms    for the case of the equation without source term.  In that theorem the transformed function 
has  non-vanishing initial values.
For  $\gamma \in {\mathbb C}$, $Re \,  \gamma >0 $, that is for $\ell \in {\mathbb C} \setminus \overline{D_1(-1,0)} = \{ z \in {\mathbb C}\,|\, |z+1| >1\}$, and, in particular,   
for $ \ell \in (-\infty,-2) \cup (0,\infty)$,  
we define the integral operator
\begin{eqnarray*}
(K_0 v)(x,t)
& := &
 2^{2-2\gamma }
\frac{\Gamma \left( 2\gamma  \right) } {\Gamma^2 \left( \gamma  \right) }
\int_{0}^1   v   (x,  \phi (t)  s)
(1-s^2)^{\gamma - 1   } ds \\
& = &
 \phi (t)  ^{1-2\gamma }   2^{2-2\gamma }
\frac{\Gamma \left( 2\gamma  \right) } {\Gamma^2 \left( \gamma  \right) }
\int_{0}^{ \phi (t) }   v (x,  \tau )
( \phi^2  (t) -\tau ^2)^{\gamma - 1} d\tau  \nonumber \,.
\end{eqnarray*}
 For    $\gamma \in {\mathbb C}$, $Re \, \gamma < 1 $, that is for $\ell \in {\mathbb C}\setminus \overline{D_1(-3,0)}= \{ z \in {\mathbb C}\,|\, |z+3| >1\}$,  and, in particular,  
for $ \ell \in (-\infty,-4) \cup (-2,\infty)$, 
we define the  integral  operator
\begin{eqnarray*}
(K_1 v)(x,t)
& := &
t 2^{2\gamma }
\frac{\Gamma \left(2- 2\gamma  \right) } {\Gamma^2 \left( 1- \gamma  \right) }
\int_{0}^1   v  (x, \phi (t) s)
(1-s^2)^{- \gamma  } ds \\
& = &
 t\phi (t)^{2\gamma -1} 2^{2\gamma }
\frac{\Gamma \left(2- 2\gamma  \right) } {\Gamma^2 \left( 1- \gamma  \right) }
\int_{0}^{\phi (t)}    v (x, \tau )
(\phi ^2 (t)-\tau ^2)^{- \gamma  } d\tau    \nonumber \,. 
\end{eqnarray*}
Thus, both operators are defined simultaneously for $\gamma \in {\mathbb C}$, $0< Re \,  \gamma < 1 $, and,   in particular,  for  $\ell \in (-\infty,-4)\cup (0,\infty)  $.  
Denote 
\[
a_\ell:=  2^{1-2\gamma }
\frac{\ell \,\Gamma \left( 2\gamma  \right) } {2\gamma \Gamma^2 \left( \gamma  \right) }\,,\qquad  
b_\ell:= (\ell+2)2^{2\gamma-1 }
\frac{\Gamma \left(2- 2\gamma  \right) } {\Gamma^2 \left( 1- \gamma  \right) }\,. 
\]
The next theorem describes the properties of the integral transforms $K_0 $ and $K_1 $ in the case when $\ell $  is a positive number.

\begin{theorem}
\label{TK0K1}
Let $\ell$ be a positive number and let $\Omega \subset {\mathbb R}_+^{n+1} $ be a  backward time connected domain.  
Suppose that the function $v \in C_{x,t}^{m,2} ( \overline{\Omega } )$ for given $(x_0,t_0) \in \Omega $  solves  the  equation
\begin{eqnarray}
\label{eqA}
\partial_{t}^2 v -  A(x,\partial_x)  v =0 \qquad \mbox{ at  } \quad x=x_0  \quad \mbox{ and all} \quad \quad t \in (0, \phi (t_0))\, .
\end{eqnarray}

Then  the functions $K_0 v \in C_{x,t}^{m,2} (  \Omega  ) $ and  $K_1 v  \in C_{x,t}^{m,2} (  \Omega  )$ satisfy equations 
 \begin{eqnarray}
\left( \partial_{t}^2-t^{\ell}A(x,\partial_x)\right) K_0v
& = &
a_{\ell} t^{\frac{\ell}{2}-1}
     \partial_t v(x, 0)  \quad \mbox{at} \quad    x=x_0   \quad \mbox{for all} \quad t \in(0,t_0)\,,  
\end{eqnarray}
and
 \begin{eqnarray}
\left( \partial_{t}^2-t^{\ell}A(x,\partial_x)\right) K_1 v
& = &
 b_{\ell}
t^{\frac{\ell}{2}}    \partial_t v   (x, 0) \quad \mbox{at} \quad    x=x_0   \quad \mbox{for all} \quad t \in(0,t_0)\,,
\end{eqnarray}
respectively. They  have at $x=x_0$ the following initial values
\begin{eqnarray*}
(K_0 v)   (x_0,0)
  =     v   (x_0, 0)\,,\qquad (K_0 v)_ t (x_0,0)
  =      0\,,
\end{eqnarray*}
and
\begin{eqnarray*}
(  K_1 v) (x_0,0)
  =   0\,,\qquad 
(  K_1 v)_t(x_0,0)
  =   v   (x_0, 0)\,.
\end{eqnarray*} 
\end{theorem}
Thus, the value   $v(x_0,0)$ of the solutions of (\ref{eqA}) is invariant under operation $K_0 $, while the operator  $K_1 $ acts similarly to the Dirichlet-to-Neumann map.

\begin{corollary}
\label{CPrepTric}
Let $\ell$ be a positive number and $\Omega \subset {\mathbb R}_+^{n+1} $ be a  backward time connected domain.  
Suppose that the function $v \in C_{x,t}^{m,2} (\overline{\Omega _\phi}  )$   solves  the  equation
\[
\partial_t^2 v-  A(x,\partial_x)  v =0     \quad \mbox{ for all} \quad  (x,t) \in \Omega _\phi \, .
\]

Then  the functions $K_0 v \in C_{x,t}^{m,2} (  \Omega  ) $ and  $K_1 v  \in C_{x,t}^{m,2} (  \Omega  )$ satisfy equations 
 \begin{eqnarray}
\left( \partial_{t}^2-t^{\ell}A(x,\partial_x)\right) K_0v
& = &
a_{\ell} t^{\frac{\ell}{2}-1}
     \partial_t v(x, 0)     \quad \mbox{for all} \quad (x,t ) \in \Omega \,,  
\end{eqnarray}
and
 \begin{eqnarray}
\left( \partial_{t}^2-t^{\ell}A(x,\partial_x)\right) K_1 v
& = &
 b_{\ell}
t^{\frac{\ell}{2}}    \partial_t v   (x, 0) \quad \mbox{for all} \quad (x,t ) \in \Omega \,,
\end{eqnarray}
respectively. They  have the following initial values
\begin{eqnarray*}
(K_0 v)   (x ,0)
  =     v   (x , 0)\,,\qquad (K_0 v)_ t (x ,0) =0 \quad \mbox{for all} \quad  x  \in \widetilde{\Omega} 
 \,,
\end{eqnarray*}
and
\begin{eqnarray*}
(  K_1 v) (x ,0)
  =   0\,,\qquad 
(  K_1 v)_t(x ,0)
  =   v   (x , 0) \quad \mbox{for all} \quad  x  \in \widetilde{\Omega}\,.
\end{eqnarray*} 
\end{corollary}
For the Cauchy problem with the full initial data, 
we have the following result
for the generalized Tricomi equation in the hyperbolic domain.

\begin{theorem}
\label{T1.8c} 
Let $\ell$ be a positive number and $\Omega \subset {\mathbb R}_+^{n+1} $ be a  backward time connected domain.  
Suppose that the functions $v_0, v_1  \in C_{x,t}^{m,2} (\overline{\Omega _\phi}  )$   solve  the  problem
\begin{eqnarray*}
&  &
\partial_t^2 v_i-  A(x,\partial_x)  v_i =0     \quad \mbox{ for all} \quad  (x,t) \in \Omega _\phi \,, \\
&  &
 v_i   (x,0)
  =     \sum_{k=0,1} \delta _{ik} \varphi _k  (x  )\,,\qquad  \partial_t v_ {i} (x,0)
  =  0    \,, \qquad  i =0,1, \quad \mbox{for all} \quad (x,t ) \in \widetilde{\Omega_\phi}\,.
\end{eqnarray*}

Then the function $u=K_0v_0 +K_1v_1  \in C_{x,t}^{m,2} (  \Omega  )$ solves the problem   
 \begin{eqnarray*}
&  &
\left( \partial_{t}^2-t^{\ell}A(x,\partial_x)\right) u
= 0      \quad \mbox{for all} \quad (x,t ) \in \Omega \,,  \\
&  &
u   (x,0)
  =    \varphi _0  (x  )\,,\qquad \, \partial_t u  (x,0)
  =   \varphi _1  (x  ) \quad \mbox{for all} \quad  x  \in \widetilde{\Omega} \,.
\end{eqnarray*}
\end{theorem}

In order to make this paper more self-contained, we remind  here that if $A(x,\partial_x)= \Delta  $, then the function $ v_\varphi  (x, t) $ is given by the following formulas (see, e.g., \cite{Shatah}):   for $\varphi \in C_0^\infty ({\mathbb R}^n)$ and for $x \in {\mathbb R}^n$, $n=2m+1$, $m \in {\mathbb N}$,
\begin{eqnarray}
\label{vphiodd}
 v_\varphi  (x, t) :=
 \frac{\partial}{\partial t} \Big( \frac{1}{t} \frac{\partial }{\partial t}\Big)^{\frac{n-3}{2} }
\frac{t^{n-2}}{\omega_{n-1} c_0^{(n)} } \int_{S^{n-1}  }
\varphi (x+ty)\, dS_y\,,
\end{eqnarray}
while for $x \in {\mathbb R}^n$, $n=2m$,  $m \in {\mathbb N}$ ,
\begin{eqnarray}
\label{vphieven}
v_\varphi  (x, t) :=    \frac{\partial }{\partial t}
\Big( \frac{1}{t} \frac{\partial }{\partial t}\Big)^{\frac{n-2}{2} }
\frac{2t^{n-1}}{\omega_{n-1} c_0^{(n)}} \int_{B_1^{n}(0)}  \frac{1}{\sqrt{1-|y|^2}}\varphi (x+ty)\, dV_y
   \,.
\end{eqnarray}
The last formulas can be also written in terms of  Radon transform; for details, see  \cite{Helgason,L-P}. 
\medskip

The case of negative $\ell$ requires some modifications in the setting of the initial conditions at $t=0$. For the  EdeS spacetime these modifications are suggested 
in \cite{Galstian-Kinoshita-Yagdjian}; they are the so-called   weighted initial conditions. For the negative $\ell$ operators  $ K_0 $ and $K_1  $ and corresponding weighted initial conditions will be discussed in the forthcoming paper.

One can consider the Cauchy problem for the equations with negative $\ell$  and with  the initial conditions prescribed at $t=t_{in}>0$.  For $\ell <-2 $, the hyperbolic equations 
in such spacetime   have permanently bounded domain of  influence. Nonlinear equations with a permanently bounded domain of  influence
were studied, in particular, in \cite{Choquet-Bruhat_Sp_2000,Choquet-Bruhat_ND_2000}. In  particular, Choquet-Bruhat \cite{Choquet-Bruhat_ND_2000} proved for small initial data
the global existence and uniqueness of wave maps on the 
FLRW expanding universe with the metric ${\bf g}=-dt^2+R^2(t)   \sigma   $  and    a smooth Riemannian manifold $(S,\sigma )$ of dimension $n \leq 3$, which has a 
  time independent metric $\sigma  $  and   non-zero injectivity radius,
and with $R(t) $ being  a positive increasing function  such that $1/R (t)$ is integrable on $[t_{in}  ,\infty)$.   If the target manifold is flat,  then the wave map equation reduces to a linear system. 
\medskip

In the forthcoming papers we also will apply integral transform approach to the  maximum principle (see, e.g., \cite{Lupo-Payne,Protter}) 
for the generalized Tricomi equation, to the derivation of the $L_p-L_q$ estimates (see \cite{Kim,Yag_Galst_CMP}), to the    
mixed problem for Friedlander model (see \cite{Ivanovici} and references therein), global existence problem for the semilinear 
generalized Tricomi equations on the hyperbolic space (for wave equation see \cite{Anker,Catania-Georgiev,Metcalfe-Taylor}), 
and  Price's law for the corrsponding cosmological models (see, e.g.,  \cite{Metcalfe-Tataru-Tohaneanu} and references therein).
\medskip

This paper is organized as follows. In Section~\ref{S2} we prove the main property  (\ref{Eeq}) of the kernel function $E $. In Section~\ref{S3} we prove Theorem~\ref{T1.8}.
In Section~\ref{S4} we present some properties of operators $K_0$ and $K_1$ and prove Theorem~\ref{T1.8b}. In Appendix we separated some technical lemmas.

\section{The kernel function $E$}
\label{S2}

The function $E$ defined by (\ref{E}) contains the hypergeometric function $F \left(\gamma , \gamma ;1; z  \right) $, where $\gamma =\frac{\ell}{2(\ell+2)} $. The  hypergeometric function is 
defined by the Gauss series on the disk $|z|<1 $, and by analytic continuation over the whole complex $z$-plane and for all complex $\ell \in {\mathbb C} \setminus \{-2\} $ (see, e.g., \cite{Slater}).
 
According to the next theorem, the function $  E(r,t;b;\gamma  )$ solves  the   Tricomi-type equation. 
The similar property (see Theorem~1.12~\cite{Yag_Arx}) possesses another kernel function, which is used to solve  Klein-Gordon equation in the de~Sitter spacetime, that is, 
the equation  with $ a (t) =e^{-t}$ and the mass term.
\begin{theorem}
\label{Teqforkernel}
The kernel function
\begin{eqnarray*}
E(r,t;b;\gamma  )
& = &
c_{\ell} \left(  (\phi (t)  + \phi (b))^2  -r^2 \right)^{-\gamma } F \left(\gamma , \gamma ;1; \frac{(\phi (t)  - \phi
(b))^2 - r^2} {(\phi (t)  + \phi (b))^2 - r^2}  \right) \,,  
\end{eqnarray*} 
where $\ell \in {\mathbb C} \setminus \{-2\} $, $ \phi (t)=\frac{2}{\ell+2}t^{\frac{\ell+2}{2}} $, $\gamma =\frac{\ell}{2(\ell+2)} $, solves the equation
\begin{eqnarray*}
E_{tt}(r,t;b;\gamma  )-(\phi' (t))^2  E_{rr}(r,t;b;\gamma  )=0\,, \quad \mbox{for all} \quad t> b>0,\quad r>0, \,\,  r^2 \not= (\phi (t)  + \phi(b))^2 \,. 
\end{eqnarray*} 
\end{theorem}

Theorem~\ref{Teqforkernel}  generalizes  corresponding   statement from \cite{Barros-Neto-Gelfand-I}, which is due to Darboux~\cite{Darboux}.
In fact, in \cite{Barros-Neto-Gelfand-I} the Tricomi equation is considered, that is, $\ell=1$, $n=1$,  $A(x,\partial_x) = \partial_x^2 $, and the statement is made for the {\sl reduced hyperbolic Tricomi operator},
which is the  Tricomi operator written in the {\sl characteristic coordinates}. For the case of reduced hyperbolic Tricomi operator with $\ell \in {\mathbb R}_+$, see   Lemma~2.1~\cite{YagTricomi}. 
The proof of Theorem~\ref{Teqforkernel}, which  is given  below, is straightforward. It works out   for different   distance functions $\phi  =\phi (t)  $,
see, for instance, \cite{Yag_Arx}, where the case of $ a (t) =e^{-t}$ is discussed. 
\medskip

The case of $\gamma =0 $ is trivial, and thereafter
	 we set $ \gamma \not=0$.
In order to prove Theorem~\ref{Teqforkernel} we need the following technical lemma. Denote
\begin{eqnarray*}
&  &
\alpha  (t,b,r) 
  :=  
\left(  (\phi (t)  + \phi (b))^2  -r^2 \right)^{-\gamma }\,,\qquad
\beta   (t,b,r)
  :=  
\frac{(\phi (t)  - \phi
(b))^2 - r^2} {(\phi (t)  + \phi (b))^2 - r^2} \,,
\end{eqnarray*}
then the next identity is evident 
\begin{equation}
\label{1minusb}
1- \beta   (t,b,r)=\frac{4\phi (t) \phi (b)}{(\phi (t)  + \phi (b))^2  -r^2}\,.
\end{equation}
\begin{lemma}
\label{alphabeta} 
The derivatives of the functions $ \alpha =\alpha   (t,b,r)$ and $\beta =\beta  (t,b,r) $ are as follows  
\begin{eqnarray*} 
\partial_t \alpha   (t,b,r)
& = &
-2\gamma \phi' (t) (\phi (t)  + \phi (b))\alpha   (t,b,r)^{\frac{\gamma +1}{\gamma }} \,,\\
\partial_t^2  \alpha   (t,b,r)
& = &
-2\gamma    \phi' {}'(t) (\phi (t)  + \phi (b))\alpha   (t,b,r)^{\frac{\gamma +1}{\gamma } }  
-2\gamma    ( \phi' (t)  )^2\alpha   (t,b,r)^{\frac{\gamma +1}{\gamma } }    \\
&  &
+ 4\gamma (\gamma +1) (\phi' (t))^2 (\phi (t)  + \phi (b))^2 \alpha^\frac{ \gamma +2}{\gamma }   (t,b,r)\,,\\
\partial_t \beta  (t,b,r)
& = &
4\alpha ^{\frac{2}{\gamma }}  (t,b,r)\phi' (t)  \phi (b) ( \phi^2 (t)  - \phi^2 (b)    + 
   r^2   )  \,,\\
\partial_t^2 \beta  (t,b,r)
& = &
 4  \phi (b) \left(\left(\phi (b)+\phi (t)\right)^2-r^2\right)^{-2} \Bigg[ \phi '{}' (t)   \left(\phi^2 (t)-\phi^2 (b)+r^2\right)    \\
&  &
- 4 (\phi ' (t))^2 \left(\phi (b)+\phi (t)\right) \left(\phi^2 (t)-\phi^2 (b)+r^2\right) 
\left(\left(\phi (b)+\phi (t)\right)^2-r^2\right)^{-1} 
+ 2 \phi (t) (\phi ' (t))^2   \Bigg] 
\end{eqnarray*}
and
\begin{eqnarray*} 
\partial_r \alpha   (t,b,r) 
& = &
 2\gamma r\left( \alpha   (t,b,r)  \right)^{\frac{\gamma+1}{\gamma } } \,,\\
\partial_r^2 \alpha   (t,b,r) 
& = &
2\gamma  \left( \alpha   (t,b,r)  \right)^{\frac{\gamma+1}{\gamma } }
+4\gamma  (\gamma+1 )r^2   \left( \alpha   (t,b,r)  \right)^{\frac{\gamma+2}{\gamma } }  \,,\\
\partial_r\beta  (t,b,r) 
& = &
- 8r\frac{   \phi (t)  \phi (b) } {((\phi (t)  + \phi (b))^2 - r^2)^2}=-2r(1- \beta  (t,b,r) )\left( \alpha   (t,b,r)  \right)^{\frac{1}{\gamma } } \,,\\
\partial_r^2\beta  (t,b,r) 
& = &
- 8 \frac{   \phi (t)  \phi (b) } {((\phi (t)  + \phi (b))^2 - r^2)^3}\left[   (\phi (t)  + \phi (b))^2  +3 r^2 \right]\\
& = &
2( \beta  (t,b,r)-1) \alpha^{\frac{2}{\gamma }}  (t,b,r)\left[   (\phi (t)  + \phi (b))^2  +3 r^2 \right]
\,.
\end{eqnarray*}
\end{lemma}
\medskip

\noindent
{\bf Proof.} These identities can be verified by straightforward calculations. 
Lemma is proved. \hfill $\square$

\medskip

\noindent
{\bf Proof of Theorem~\ref{Teqforkernel}.}  
It is sufficient to consider the function  
\begin{eqnarray*}
E_{c_\ell}(r,t;b;\gamma  )
& := &
\alpha  (t,b,r) F \left(\gamma , \gamma ;1; \beta   (t,b,r) \right) \,,  
\end{eqnarray*}  
since $E(r,t;b;\gamma  ) =c_\ell E_{c_\ell}(r,t;b;\gamma  )$. We have 
\begin{eqnarray*}
\partial_r E_{c_\ell}(r,t;b;\gamma  ) 
& = &
(\partial_r \alpha  (t,b,r)) F \left(\gamma , \gamma ;1; \beta   (t,b,r) \right) + \alpha  (t,b,r) F_z \left(\gamma , \gamma ;1; \beta   (t,b,r) \right)\partial_r \beta   (t,b,r)\,.
\end{eqnarray*} 
Hence,
\begin{eqnarray*}
\partial_r^2 E_{c_\ell}(r,t;b;\gamma  ) 
& = & 
    \alpha_{rr}  (t,b,r)  F \left(\gamma , \gamma ;1; \beta   (t,b,r) \right)  
+   2  \alpha_r  (t,b,r)  \beta_r   (t,b,r) F_z \left(\gamma , \gamma ;1; \beta   (t,b,r) \right)  \\
 &  &
 +   \alpha  (t,b,r) F_{zz} \left(\gamma , \gamma ;1; \beta   (t,b,r) \right)  (\beta_r   (t,b,r))^2  \\
 &  &
 + \alpha  (t,b,r) F_z \left(\gamma , \gamma ;1; \beta   (t,b,r) \right)  \beta_{rr}   (t,b,r) \,,
\end{eqnarray*} 
while
\begin{eqnarray*}
\partial_t^2 E_{c_\ell}(r,t;b;\gamma  ) 
& = & 
    \alpha_{tt}  (t,b,r)  F \left(\gamma , \gamma ;1; \beta   (t,b,r) \right)  
+   2  \alpha_t  (t,b,r)  \beta_t   (t,b,r) F_z \left(\gamma , \gamma ;1; \beta   (t,b,r) \right)  \\
 &  &
 +   \alpha  (t,b,r) F_{zz} \left(\gamma , \gamma ;1; \beta   (t,b,r) \right)  (\beta_t   (t,b,r))^2  \\
 &  &
 + \alpha  (t,b,r) F_z \left(\gamma , \gamma ;1; \beta   (t,b,r) \right)  \beta_{tt}   (t,b,r) \,.
\end{eqnarray*} 
Thus, we obtain
\begin{eqnarray*}
&  &
\partial_t^2 E_{c_\ell}(r,t;b;\gamma  )-t^{\ell}  \partial_t^2E_{c_\ell}(r,t;b;\gamma  )\\
&  = &
    \alpha_{tt}  (t,b,r)  F \left(\gamma , \gamma ;1; \beta   (t,b,r) \right)  
+   2  \alpha_t  (t,b,r)  \beta_t   (t,b,r) F_z \left(\gamma , \gamma ;1; \beta   (t,b,r) \right)  \\
 &  &
 +   \alpha  (t,b,r) F_{zz} \left(\gamma , \gamma ;1; \beta   (t,b,r) \right)  (\beta_t   (t,b,r))^2  \\
 &  &
 + \alpha  (t,b,r) F_z \left(\gamma , \gamma ;1; \beta   (t,b,r) \right)  \beta_{tt}   (t,b,r)  \\
 &  &
-t^{\ell} \Bigg[     \alpha_{rr}  (t,b,r)  F \left(\gamma , \gamma ;1; \beta   (t,b,r) \right)  
+   2  \alpha_r  (t,b,r)  \beta_r   (t,b,r) F_z \left(\gamma , \gamma ;1; \beta   (t,b,r) \right)  \\
 &  &
 +   \alpha  (t,b,r) F_{zz} \left(\gamma , \gamma ;1; \beta   (t,b,r) \right)  (\beta_r   (t,b,r))^2  \\
 &  &
 + \alpha  (t,b,r) F_z \left(\gamma , \gamma ;1; \beta   (t,b,r) \right)  \beta_{rr}   (t,b,r) \Bigg]\,.
\end{eqnarray*}
It follows
\begin{eqnarray}
\label{Ecl}
&  &
\partial_t^2 E_{c_\ell}(r,t;b;\gamma  )-t^{\ell}   \partial_t^2E_{c_\ell}(r,t;b;\gamma  ) \nonumber \\
&  = &
I(t,b,r) F \left(\gamma , \gamma ;1; \beta   (t,b,r) \right)  + J(t,b,r)  F_z \left(\gamma , \gamma ;1; \beta   (t,b,r) \right)  +Y(t,b,r) F_{zz} \left(\gamma , \gamma ;1; \beta   (t,b,r) \right) \,,
\end{eqnarray}
where
\begin{eqnarray*} 
I (t,b,r) 
&  = &
   \alpha_{tt}  (t,b,r) -t^{\ell}   \alpha_{rr}  (t,b,r)  \,, \\
J(t,b,r)
& = &
  2  \alpha_t  (t,b,r)  \beta_t   (t,b,r)   +\alpha  (t,b,r)    \beta_{tt}   (t,b,r)  
-t^{\ell} \left( 2  \alpha_r  (t,b,r)  \beta_r   (t,b,r)  +\alpha  (t,b,r)  \beta_{rr}   (t,b,r)\right)  \,, \\ 
Y(t,b,r)
 & = &
   \alpha  (t,b,r)   (\beta_t   (t,b,r))^2  -  t^{\ell} \alpha  (t,b,r)  (\beta_r   (t,b,r))^2  \,.  
\end{eqnarray*}
We have due to Lemma~\ref{alphabeta} 
\begin{eqnarray*}  
I (t,b,r)  
&  = &
-2\gamma    \phi' {}'(t) (\phi (t)  + \phi (b))\alpha   (t,b,r)^{\frac{\gamma +1}{\gamma } }  
-2\gamma    ( \phi' (t)  )^2\alpha   (t,b,r)^{\frac{\gamma +1}{\gamma } }    \\
&  &
+ 4\gamma (\gamma +1) (\phi' (t))^2 (\phi (t)  + \phi (b))^2 \alpha^\frac{ \gamma +2}{\gamma }   (t,b,r) \\
 &  &
 -t^{\ell}  \left( 2\gamma  \left( \alpha   (t,b,r)  \right)^{\frac{\gamma+1}{\gamma } }
+4\gamma  (\gamma+1 )r^2   \left( \alpha   (t,b,r)  \right)^{\frac{\gamma+2}{\gamma } }  \right)\,.
\end{eqnarray*}
We can write
\begin{eqnarray*} 
I (t,b,r)
&  = &
-2\gamma    \phi' {}'(t) (\phi (t)  + \phi (b))\left(  (\phi (t)  + \phi (b))^2  -r^2 \right)^{-\gamma -1}  
-2\gamma    ( \phi' (t)  )^2 \left(  (\phi (t)  + \phi (b))^2  -r^2 \right)^{-\gamma -1}     \\
&  &
+ 4\gamma (\gamma +1) (\phi' (t))^2 (\phi (t)  + \phi (b))^2 \left(  (\phi (t)  + \phi (b))^2  -r^2 \right)^{-\gamma-2 } \\
 &  &
 -t^{\ell}  \left( 2\gamma    \left(  (\phi (t)  + \phi (b))^2  -r^2 \right)^{-\gamma-1}  
+4\gamma  (\gamma+1 )r^2    \left(  (\phi (t)  + \phi (b))^2  -r^2 \right)^{-\gamma-2} \right)\\
&  = &
\left(  (\phi (t)  + \phi (b))^2  -r^2 \right)^{-\gamma-2}\\
&  &
\times \Bigg\{ -2\gamma    \phi' {}'(t) (\phi (t)  + \phi (b))\left(  (\phi (t)  + \phi (b))^2  -r^2 \right) 
-2\gamma    ( \phi' (t)  )^2 \left(  (\phi (t)  + \phi (b))^2  -r^2 \right)   \\
&  &
+ 4\gamma (\gamma +1) (\phi' (t))^2 (\phi (t)  + \phi (b))^2  
 -t^{\ell}  \left( 2\gamma    \left(  (\phi (t)  + \phi (b))^2  -r^2 \right)  
+4\gamma  (\gamma+1 )r^2   \right)\Bigg\}\,.
\end{eqnarray*}
Then we use the identities 
\begin{eqnarray}
\label{phi}
t^{\ell} = ( \phi' (t)  )^2\,,\qquad \frac{1}{2\gamma }\phi'{}' (t)\phi (t)= \left( \phi' (t)\right)^2 \,,
\end{eqnarray} 
and derive
\begin{eqnarray*}  
I (t,b,r)  
&  = &
\left(  (\phi (t)  + \phi (b))^2  -r^2 \right)^{-\gamma-2}\\
&  &
\times \Bigg\{ -2\gamma    \phi' {}'(t) (\phi (t)  + \phi (b))\left(  (\phi (t)  + \phi (b))^2  -r^2 \right) 
-   \phi' {}' (t)   \phi (t)  \left(  (\phi (t)  + \phi (b))^2  -r^2 \right)   \\
&  &
+ 2 (\gamma +1) \phi' {}'(t)   \phi (t) (\phi (t)  + \phi (b))^2  \\
 &  &
 - \frac{1}{2\gamma } \phi' {}'(t)   \phi (t) \left[ 2\gamma    \left(  (\phi (t)  + \phi (b))^2  -r^2 \right)  
+4\gamma  (\gamma+1 )r^2   \right]\Bigg\}\,.
\end{eqnarray*}
It follows 
\begin{eqnarray*} 
I (t,b,r)
&  = &
\left(  (\phi (t)  + \phi (b))^2  -r^2 \right)^{-\gamma-2}\\
&  &
\times \Bigg\{ -2\gamma    \phi' {}'(t) (\phi (t)  + \phi (b))\left(  (\phi (t)  + \phi (b))^2  -r^2 \right) 
-  \phi' {}' (t)   \phi (t)   \left(  (\phi (t)  + \phi (b))^2  -r^2 \right)  \\
&  &
+ 2 \gamma   \phi' {}'(t)   \phi (t) (\phi (t)  + \phi (b))^2 +  2  \phi' {}'(t)   \phi (t) (\phi (t)  + \phi (b))^2   \\
 &  &
 -   \phi' {}'(t)   \phi (t) \left[     \left(  (\phi (t)  + \phi (b))^2  -r^2 \right)  
+2 (\gamma+1 )r^2   \right]\Bigg\}\\
&  = &
\phi' {}'(t)\left(  (\phi (t)  + \phi (b))^2  -r^2 \right)^{-\gamma-2}\\
&  &
\times \Bigg\{   -2\gamma      \phi (b) \left(  (\phi (t)  + \phi (b))^2  -r^2 \right) 
-      2\phi (t)   \left(  (\phi (t)  + \phi (b))^2  -r^2 \right)  \\
&  &
+  2      \phi (t) (\phi (t)  + \phi (b))^2  -     2   \phi (t) r^2   \Bigg\}   \,.
\end{eqnarray*}
Thus,
\begin{eqnarray} 
\label{I} 
I (t,b,r)  
&  = &
- 2\gamma    \phi' {}' (t)   \phi (b)\left(  (\phi (t)  + \phi (b))^2-r^2\right)^{-\gamma -1}  \,.
\end{eqnarray}

Further, we consider the coefficient $J (t,b,r)  $:  
\begin{eqnarray*} 
J (t,b,r) &  = & 
2  \Bigg[ -2\gamma \phi' (t) (\phi (t)  + \phi (b))\alpha   (t,b,r)^{(\gamma +1) /\gamma }\Bigg]  \Bigg[4\alpha ^{\frac{2}{\gamma }}(t,b,r) \phi' (t)  \phi (b) ( \phi^2 (t)  - \phi^2 (b)    + 
   r^2   )  \Bigg]  \\
&  &
 +\alpha  (t,b,r)       4  \phi (b) \left(\left(\phi (b)+\phi (t)\right)^2-r^2\right)^{-2} \Bigg[ \phi '{}' (t)   \left(\phi^2 (t)-\phi^2 (b)+r^2\right)    \\
&  &
- 4 (\phi ' (t))^2 \left(\phi (b)+\phi (t)\right) \left(\phi^2 (t)-\phi^2 (b)+r^2\right) 
\left(\left(\phi (b)+\phi (t)\right)^2-r^2\right)^{-1} \\
&  &
+ 2 \phi (t) (\phi ' (t))^2   \Bigg]    
-t^{\ell}  2   \Bigg[  2\gamma r  \alpha^{\frac{\gamma+1}{\gamma } }   (t,b,r)   \Bigg]  \Bigg[-2r(1- \beta  (t,b,r) )
  \alpha^{\frac{1}{\gamma } }   (t,b,r)   \Bigg]  \\
&   &
-t^{\ell} \alpha  (t,b,r)  \Bigg[- 8 \frac{   \phi (t)  \phi (b) } {((\phi (t)  + \phi (b))^2 - r^2)^3}
\left[   (\phi (t)  + \phi (b))^2  +3 r^2 \right]\Bigg] \,.  
\end{eqnarray*}
Next we use the identities (\ref{phi}) and definition of $\alpha  (t,b,r)$ in order to rewrite $ J$ as follows:
\begin{eqnarray*}  
J (t,b,r)  
&  = & 
-16    \gamma \frac{1}{2\gamma } \phi '{}'(t) \phi (t) \phi (b)(\phi (t)  + \phi (b))\left(\left(\phi (b)+\phi (t)\right)^2-r^2\right)^{-\gamma -3}  ( \phi^2 (t)  - \phi^2 (b)    + 
   r^2   )    \\
&  &
 +\left(\left(\phi (b)+\phi (t)\right)^2-r^2\right)^{-\gamma }     4  \phi (b) \left(\left(\phi (b)+\phi (t)\right)^2-r^2\right)^{-2} \Bigg[ \phi '{}' (t)   \left(\phi^2 (t)-\phi^2 (b)+r^2\right)    \\
&  &
- 4 \frac{1}{2\gamma } \phi '{}'(t) \phi (t)  \left(\phi (b)+\phi (t)\right) \left(\phi^2 (t)-\phi^2 (b)+r^2\right) 
\left(\left(\phi (b)+\phi (t)\right)^2-r^2\right)^{-1} \\
&  &
+ 2 \phi (t)\frac{1}{2\gamma } \phi '{}'(t) \phi (t)   \Bigg]    \\
 &  &
-\frac{1}{2\gamma } \phi '{}'(t) \phi (t) 2   \Bigg[  2\gamma r \left( \alpha   (t,b,r)  \right)^{\frac{\gamma+1}{\gamma } } \Bigg]  \Bigg[-2r(1- \beta  (t,b,r) )\left( \alpha   (t,b,r)  \right)^{\frac{1}{\gamma } } \Bigg]  \\
&   &
-\frac{1}{2\gamma } \phi '{}'(t) \phi (t)  \alpha  (t,b,r)  \Bigg[- 8 \frac{   \phi (t)  \phi (b) } {((\phi (t)  + \phi (b))^2 - r^2)^3}\left[   (\phi (t)  + \phi (b))^2  +3 r^2 \right]\Bigg] \,.  
\end{eqnarray*}
Hence
\begin{eqnarray*} 
J (t,b,r) 
&  = & 
-8 \phi '{}'(t) \phi (t) \phi (b)(\phi (t)  + \phi (b))\left(\left(\phi (b)+\phi (t)\right)^2-r^2\right)^{-\gamma -3}  ( \phi^2 (t)  - \phi^2 (b)    + 
   r^2   )    \\
&  &
 +\left(\left(\phi (b)+\phi (t)\right)^2-r^2\right)^{-\gamma-2 }     4  \phi (b)   \phi '{}' (t) \Bigg[   \left(\phi^2 (t)-\phi^2 (b)+r^2\right)    \\
&  &
- 2 \frac{1}{ \gamma }  \phi (t)  \left(\phi (b)+\phi (t)\right) \left(\phi^2 (t)-\phi^2 (b)+r^2\right) 
\left(\left(\phi (b)+\phi (t)\right)^2-r^2\right)^{-1} \\
&  &
+   \phi (t)\frac{1}{ \gamma }   \phi (t)   \Bigg]   
+  16\phi '{}'(t) \phi (t)       r^2 \left(\left(\phi (b)+\phi (t)\right)^2-r^2\right)^{-\gamma -3}   
   \phi (t)  \phi (b)   \\
&   &
+4\frac{1}{ \gamma } \phi '{}'(t) \phi (t) \left(\left(\phi (b)+\phi (t)\right)^2-r^2\right)^{-\gamma -3}        \phi (t)  \phi (b)  \left[   (\phi (t)  + \phi (b))^2  +3 r^2 \right]\,.
\end{eqnarray*}
Then
\begin{eqnarray*} 
J (t,b,r) 
&  = & 
\phi '{}'(t)  \phi (b)\Bigg\{ -8 \phi (t)(\phi (t)  + \phi (b))\left(\left(\phi (b)+\phi (t)\right)^2-r^2\right)^{-\gamma -3}  ( \phi^2 (t)  - \phi^2 (b)    + 
   r^2   )    \\
&  &
 +\left(\left(\phi (b)+\phi (t)\right)^2-r^2\right)^{-\gamma-2 }     4    \Bigg[   \left(\phi^2 (t)-\phi^2 (b)+r^2\right)    \\
&  &
- 2 \frac{1}{ \gamma }  \phi (t)  \left(\phi (b)+\phi (t)\right) \left(\phi^2 (t)-\phi^2 (b)+r^2\right) 
\left(\left(\phi (b)+\phi (t)\right)^2-r^2\right)^{-1} 
+   \phi^2 (t)\frac{1}{ \gamma }    \Bigg]    \\
 &  &
+  16 \phi (t)       r^2 \left(\left(\phi (b)+\phi (t)\right)^2-r^2\right)^{-\gamma -3}   
   \phi (t)   \\
&   &
+4\frac{1}{ \gamma } \phi (t) \left(\left(\phi (b)+\phi (t)\right)^2-r^2\right)^{-\gamma -3}        \phi (t)  \left[   (\phi (t)  + \phi (b))^2  +3 r^2 \right]\Bigg\}\,,
\end{eqnarray*}
and, consequently, 
\begin{eqnarray*}  
J (t,b,r)  
&  = & 
\phi '{}'(t)  \phi (b)\left( (\phi (b)+\phi (t) )^2-r^2\right)^{-\gamma -3} \Bigg\{ -8 \phi (t)(\phi (t)  + \phi (b)) ( \phi^2 (t)  - \phi^2 (b)    + 
   r^2   )    \\
&  &
 +       4   \left(\left(\phi (b)+\phi (t)\right)^2-r^2\right)   \left(\phi^2 (t)-\phi^2 (b)+r^2\right)    \\
&  &
-8     \frac{1}{ \gamma }  \phi (t)  \left(\phi (b)+\phi (t)\right) \left(\phi^2 (t)-\phi^2 (b)+r^2\right)   
+   \phi^2 (t)\frac{1}{ \gamma } 4   \left(\left(\phi (b)+\phi (t)\right)^2-r^2\right)       \\
 &  &
+  16 \phi^2 (t)       r^2      
+4\frac{1}{ \gamma } \phi^2 (t)  \left[   (\phi (t)  + \phi (b))^2  +3 r^2 \right]\Bigg\}\,.
\end{eqnarray*}
The terms containing  $  \gamma   $ in the denominators are 
\begin{eqnarray*} 
&   & 
-8     \frac{1}{ \gamma }  \phi (t)  \left(\phi (b)+\phi (t)\right) \left(\phi^2 (t)-\phi^2 (b)+r^2\right)   
+   \phi^2 (t)\frac{1}{ \gamma } 4   \left(\left(\phi (b)+\phi (t)\right)^2-r^2\right)       \\
 &  & 
+4\frac{1}{ \gamma } \phi^2 (t)  \left[   (\phi (t)  + \phi (b))^2  +3 r^2 \right]  
   =   
8\frac{1}{ \gamma } \phi (t)   \phi (b)   \left(\phi^2 (t)-\phi^2 (b)+r^2\right)  \,.
\end{eqnarray*}
The terms without $  \gamma $ are
\begin{eqnarray*} 
&   & 
  -8 \phi (t)(\phi (t)  + \phi (b)) ( \phi^2 (t)  - \phi^2 (b)    +    r^2   )   
 +       4   \left(\left(\phi (b)+\phi (t)\right)^2-r^2\right)   \left(\phi^2 (t)-\phi^2 (b)+r^2\right)     
+  16 \phi^2 (t)       r^2  \\ 
& & 
 = - 4 \left(  (\phi (b)+\phi (t))^2-r^2\right) \left( (\phi (t)-\phi (b))^2-r^2\right)   \,.
\end{eqnarray*}
Thus we have
\begin{eqnarray*}  
J (t,b,r) 
&  = & 
\phi '{}'(t)  \phi (b)\left( (\phi (b)+\phi (t) )^2-r^2\right)^{-\gamma -3}\\
&  &
\times  \Bigg\{8\frac{1}{ \gamma } \phi (t)   \phi (b)   \left(\phi^2 (t)-\phi^2 (b)+r^2\right)- 4 \left(  (\phi (b)+\phi (t))^2-r^2\right) \left( (\phi (t)-\phi (b))^2-r^2\right)  \Bigg\}\,.
\end{eqnarray*}
Finally
\begin{eqnarray}
\label{J}  
J (t,b,r)  
&  = & 
-4   \phi'{}' (t)\phi (b)  ((\phi (t) + \phi (b))^2-r^2)^{-\gamma -2}  \left[- \frac{2}{\gamma}  \phi (b) \phi (t) +  ( \phi (t)-\phi  (b))^2 -   r^2 \right] \,.
\end{eqnarray}

Next we turn to the coefficient $Y (t,b,r) $: 
\begin{eqnarray*} 
Y (t,b,r)  
&  = & 
\alpha  (t,b,r)  \Big\{ (\beta_t   (t,b,r))^2  -  t^{\ell}    (\beta_r   (t,b,r))^2 \Big\}\,.
\end{eqnarray*}
First we consider the second factor of the right-hand side of the last expression and apply Lemma~\ref{alphabeta} 
\begin{eqnarray*} 
&  &
(\beta_t   (t,b,r))^2  -  t^{\ell}    (\beta_r   (t,b,r))^2 \\
& = &
4 \alpha ^{\frac{2}{\gamma }}(t,b,r)(\phi' (t) )^2\Bigg\{  4\alpha ^{\frac{2}{\gamma }}(t,b,r)     \phi^2 (b) ( \phi^2 (t)  - \phi^2 (b)    + 
   r^2   )^2  -    r^2(1- \beta  (t,b,r) )  ^2 \Bigg\}.
\end{eqnarray*}   
Then we use (\ref{1minusb}):   
\begin{eqnarray*} 
&  &
(\beta_t   (t,b,r))^2  -  t^{\ell}    (\beta_r   (t,b,r))^2 \\
& = &
4 \alpha ^{\frac{2}{\gamma }}(t,b,r)(\phi' (t) )^2\Bigg\{  4 ((\phi  (t) +\phi (b))^2-r^2)^{-2}   \phi^2 (b) ( \phi^2 (t)  - \phi^2 (b)    + 
   r^2   )^2  -    r^2 \frac{16\phi^2 (t)\phi^2 (b) }{((\phi  (t) +\phi (b))^2-r^2)^2}     \Bigg\}\\
& = &
16 \alpha ^{\frac{2}{\gamma }}(t,b,r)(\phi' (t) )^2 \phi^2 (b)   \left((\phi (t) +\phi  (b))^2-r^2\right)^{-1} \left((\phi  (t)-\phi (b))^2 -r^2\right)  \,.
\end{eqnarray*}
Hence 
\begin{eqnarray*} 
Y (t,b,r)  
&  = &  
\alpha  (t,b,r)16 \alpha ^{\frac{2}{\gamma }}(t,b,r)(\phi' (t) )^2 \phi^2 (b) \left((\phi (t) +\phi  (b))^2-r^2\right)^{-1}  \left((\phi  (t)-\phi (b))^2 -r^2\right)\\
& = &
 16 \alpha ^{\frac{\gamma +3}{\gamma }}(t,b,r)(\phi' (t) )^2 \phi^2 (b)  \left((\phi  (t)-\phi (b))^2 -r^2\right)\,.
\end{eqnarray*}
Finally
\begin{eqnarray} 
\label{Y} 
Y (t,b,r)  
&  = & 
  16   \phi^2 (b)  (\phi' (t) )^2  \left( \left(\phi (b) -\phi (t)\right)^2 -r^2\right) \left( \left(\phi (b) +\phi (t)\right)^2 -r^2\right)^{-\gamma -3}\,.
\end{eqnarray}
We denote 
\[
G(t,b,r): =  2\gamma^{-1}    \phi' {}' (t)   \phi (b)\left(  (\phi (t)  + \phi (b))^2-r^2\right)^{-\gamma -1} \,.
\]
\begin{lemma} 
\label{L2.3}
Let $z:= \beta  (t,b,r) $, then 
\begin{equation} 
I(t,b,r)= -\gamma ^2 G(t,b,r)\,,\quad J(t,b,r)=\left( 1-(2  \gamma +1)z \right)G(t,b,r) \,,\quad   Y(t,b,r) = z(1-z)G(t,b,r)\,.
\end{equation}
\end{lemma}
\medskip

\noindent
{\bf Proof.} The first equation is evident. For the second one we calculate the ratio   
\begin{eqnarray*}  
 \frac{J(t,b,r)}{ G(t,b,r)} 
& = &
-\frac{ 2      \left[-  2  \phi (b) \phi (t) + \gamma ( \phi (t)-\phi  (b))^2 -\gamma r^2 \right]  }{     \left(  (\phi (t)  + \phi (b))^2-r^2\right) }\,.
\end{eqnarray*}
On the other hand, the following identity is easily seen 
\begin{eqnarray*} 
  1-(2\gamma +1)z  
& = &
- \frac{2}{\left(\phi (b)+\phi (t)\right)^2-r^2}\left[-2\phi (b) \phi (t) +  \gamma    \left(\phi (t)-\phi (b)\right)^2 + \gamma r^2 \right]\,.
\end{eqnarray*}
Hence the second equation is proved. For the third one we consider the ratio $Y/G $ and reduce it to the following
\begin{eqnarray*}  
 \frac{Y(t,b,r)}{ G(t,b,r)} 
& = &
\frac{ 4  \phi  (b)      \phi  (t)  \left( \left(\phi (b) -\phi (t)\right)^2 -r^2\right)  }
{    \left(  (\phi (t)  + \phi (b))^2-r^2\right)^{2}}\,.
\end{eqnarray*}
On the other hand, the following identity is evident
\begin{eqnarray*} 
 z(1-z) 
& = &
 \frac{(\phi (t)  - \phi
(b))^2 - r^2} {(\phi (t)  + \phi (b))^2 - r^2}  \frac{4 \phi (t)   \phi
(b)  } {(\phi (t)  + \phi (b))^2 - r^2} \,.
\end{eqnarray*}
Lemma is proved. \hfill $ \square$

\noindent
{\bf Completion of the proof of Theorem~\ref{Teqforkernel}.}
We use Lemma~\ref{L2.3} in  equation (\ref{Ecl}):
\begin{eqnarray*} 
&  &
\partial_t^2 E_{c_\ell}(r,t;b;\gamma  )-t^{\ell}   \partial_t^2E_{c_\ell}(r,t;b;\gamma  ) \nonumber \\
&  = &
G(t,b,r) \left( - \gamma ^2 F \left(\gamma , \gamma ;1; z\right)  +  \left( 1-(2  \gamma +1)z \right)  F_z \left(\gamma , \gamma ;1; z \right)  + z(1-z) F_{zz} \left(\gamma , \gamma ;1; z \right) \right) \,.
\end{eqnarray*}
The second factor of the right-hand side vanishes due to  the equation for the hypergeometric function $F \left(\gamma , \gamma ;1; z \right)  $.
  Theorem~\ref{Teqforkernel} is proved. \hfill $ \square$\\ 

In the remaining sections Theorem~\ref{Teqforkernel} is used for the real $\ell$.

\section{The problem with vanishing initial data. Proof of Theorem~\ref{T1.8} and Theorem~\ref{T1.8b}}
\label{S3}

The proof of Theorem~\ref{Teqforkernel} given in the previous section is instructive in its own right. 
More exactly, the proofs which are used in \cite{Barros-Neto-Gelfand-I},\cite{Darboux},\cite{YagTricomi} are based on  the transformation   of  equation (\ref{0.3JDE}) 
with $A(x,\partial_x) = \Delta $ and $\ell=1$, due to the  choice of the characteristic coordinates. For   Theorem~\ref{T1.8}, where $x=x_0 $ is frozen,  and even for  Theorem~\ref{T1.8b},  with $x$ running over  ${\mathbb R}^n$  and with the general $A(x,\partial_x) $, that approach does not work, while, as it will be shown in 
this section, the straightforward check of the solution formula can be carried out. It is important that the proof of Theorem~\ref{T1.8}
 presented in this section contains the key calculations inherited from the proof of Theorem~\ref{Teqforkernel}. The proof of Theorem~\ref{T1.8} also uses Lemma~\ref{L2.3}, 
which summarizes those calculations.
\medskip

In this section we consider the case of the  vanishing initial values of $u$.
Assume that $\gamma <1 $ and consider   the function $u=u(x,t)$ written in terms of the auxiliary functions  $\alpha  (t,b,r)$ and $\beta   (t,b,r)$ as follows
\begin{eqnarray}
\label{33}
u(x,t)
& = &
c_\ell \int_{ 0}^{t} db
  \int_{ 0}^{ \phi (t)- \phi (b)} \alpha  (t,b,r)  F \left(\gamma , \gamma ;1; \beta  (t,b,r) \right)  w(x,r;b )
dr, \quad  t>0\,, 
\end{eqnarray}
where $w(x,r;b ) $ is defined by (\ref{8}) and (\ref{9}), and  where we skip  a subindex $0$ of $x_0 $.  In the integral $0< \beta  (t,b,r)<1 $. 

The integrand 
\begin{eqnarray*}
\left(  (\phi (t)  + \phi (b))^2  -r^2 \right)^{-\gamma }   F \left(\gamma , \gamma ;1; \frac{(\phi (t)  - \phi
(b))^2 - r^2} {(\phi (t)  + \phi (b))^2 - r^2} \right)  w(x,r;b ) 
  \quad  t>0,\,\, b \in (0,t) ,\,\, r \in (0, \phi (t)- \phi (b))\,, 
\end{eqnarray*}
of the last iterated integral  is a smooth function except possibly the endpoint $\beta  (t,b,r)=1 $, since  the hypergeometric function $F  \left( a,b;c;z  \right) $
may have singularity at $z=1$  depending on the parameters $a, b, c$. 
  That singularity may affect convergence of the improper integral.   

\begin{proposition}
The integral of  (\ref{33}) defines a continuous function $u   $ such that $A(x,\partial_x)u$ is also continuous and
\begin{eqnarray*} 
A(x,\partial_x)u(x,t)
& = &
c_\ell \int_{ 0}^{t} db
  \int_{ 0}^{ \phi (t)- \phi (b)} \alpha  (t,b,r)  F \left(\gamma , \gamma ;1; \beta  (t,b,r) \right) A(x,\partial_x) w(x,r;b )
dr, \quad  t>0\,, 
\end{eqnarray*}
Moreover, for every compact set $C $ we have
\begin{eqnarray*} 
\max_{x \in C}|u(x,t)|
& \leq  &
c t^{2}   \quad  \mbox{ for small}  \quad  t > 0 \,. 
\end{eqnarray*} 
\end{proposition}
\medskip

\noindent
{\bf Proof.} From now on we use $A \lesssim  B$ to denote the statement that $A\leq CB $ for some absolute constant $C>0$.
According to Lemma~\ref{L5.1},
  for $\gamma <1 $
the integrand  can be estimated as follows
\begin{eqnarray*}
&  &
\left| \left(  (\phi (t)  + \phi (b))^2  -r^2 \right)^{-\gamma }   
F \left(\gamma , \gamma ;1; \frac{(\phi (t)  - \phi (b))^2 - r^2} {(\phi (t)  + \phi (b))^2 - r^2} \right)  w(x,r;b ) \right|\\
& \lesssim   &
   1+ \left(  (\phi (t)  + \phi (b))^2  -r^2 \right)^{-\gamma-1 }\phi (t)  \phi (b) +\left(  (\phi (t)  + \phi (b))^2  -r^2 \right)^{\gamma -1} \phi (t)^{1-2\gamma }  \phi (b) ^{1-2\gamma }   \,.  
\end{eqnarray*}
On the other hand
\begin{eqnarray*}
&  &
\left| \int_{ 0}^{t} db
  \int_{ 0}^{ \phi (t)- \phi (b)} \alpha  (t,b,r)  F \left(\gamma , \gamma ;1; \beta  (t,b,r) \right)  w(x,r;b )
dr \right| \\
&  \lesssim   &
    \int_{ 0}^{t} db \,  
  \int_{ 0}^{ \phi (t)- \phi (b)} \left(  (\phi (t)  + \phi (b))^2  -r^2 \right)^{-\gamma }  dr 
+    \phi (t) \int_{ 0}^{t} db \, \phi (b) 
  \int_{ 0}^{ \phi (t)- \phi (b)} \left(  (\phi (t)  + \phi (b))^2  -r^2 \right)^{-\gamma-1 }  dr \\
&   &
+ \,  \phi (t)^{1-2\gamma }  \int_{ 0}^{t} db \,\phi (b) ^{1-2\gamma }
  \int_{ 0}^{ \phi (t)- \phi (b)} \left(  (\phi (t)  + \phi (b))^2  -r^2 \right)^{ \gamma -1 }     dr  \,.
\end{eqnarray*}
Then we apply Lemma~\ref{L5.2} and obtain
\begin{eqnarray*}
&  &
\left| \int_{ 0}^{t} db
  \int_{ 0}^{ \phi (t)- \phi (b)} \alpha  (t,b,r)  F \left(\gamma , \gamma ;1; \beta  (t,b,r) \right)  w(x,r;b )
dr \right| \\
& \lesssim   &
    \int_{ 0}^{t} \left(\phi (t)  - \phi (b) \right)  \left(\phi (t)  + \phi (b)\right)^{-2   \gamma   } F \left(\frac{1}{2},  \gamma ;\frac{3}{2};\frac{\left(\phi (t)- \phi (b)\right)^2}{\left(\phi (t)  + \phi (b)\right)^2}\right) db\\
&   & 
+   \, \phi (t) \int_{ 0}^{t} \phi (b) \left( \phi (t)  + \phi (b) \right)^{- 2\gamma  -2 } ( \phi (t)- \phi (b)) \,  
F \left(\frac{1}{2},\gamma +1;\frac{3}{2};\frac{(\phi (t)  - \phi (b))^2  } {(\phi (t)  + \phi (b))^2  } \right) db\\
&   & 
+  \,  \phi (t)^{1-2\gamma }  \int_{ 0}^{t} db \,\phi (b) ^{1-2\gamma }
 \left(\phi (t)  - \phi (b) \right)  \left(\phi (t)  + \phi (b)\right)^{ 2  \gamma -2  } F \left(\frac{1}{2},1-\gamma ;\frac{3}{2};\frac{\left(\phi (t)- \phi (b)\right)^2}{\left(\phi (t)  + \phi (b)\right)^2}\right) .
 \end{eqnarray*}
 Due to Lemmas~\ref{L5.4},\ref{L5.5},\ref{L5.6}, the last inequality implies  
 \begin{eqnarray*}
&  &
\left| \int_{ 0}^{t} db
  \int_{ 0}^{ \phi (t)- \phi (b)} \alpha  (t,b,r)  F \left(\gamma , \gamma ;1; \beta  (t,b,r) \right)  w(x,r;b )
dr \right| \\
&  \lesssim   &
   \int_{ 0}^{t} \left(\phi (t)  - \phi (b) \right)  \left(\phi (t)  + \phi (b)\right)^{-2   \gamma   } \\
 & &
\Bigg\{\left[   1 +\left( \frac{  \phi (t)   \phi (b)  } {(\phi (t)  + \phi (b))^2  } \right)\right]  
 + \left( \frac{  \phi (t)   \phi (b)  } {(\phi (t)  + \phi (b))^2  } \right)^{1-\gamma }
\left[ 1+ \left( \frac{  \phi (t)   \phi (b)  } {(\phi (t)  + \phi (b))^2  } \right) \right]  \Bigg\}  db\\
&   & 
+    \, \phi (t) \int_{ 0}^{t} \phi (b) \left( \phi (t)  + \phi (b) \right)^{- 2\gamma  -2 } ( \phi (t)- \phi (b)) \\ 
&  &
 \times  \Bigg\{
\left[ 1+ \left( \frac{ \phi (t)   \phi (b)  } {(\phi (t)  + \phi (b))^2  } \right) \right]   
 + \left( \frac{  \phi (t)   \phi (b)  } {(\phi (t)  + \phi (b))^2  } \right)^{-\gamma }
\left[ 1+  \left( \frac{  \phi (t)   \phi (b)  } {(\phi (t)  + \phi (b))^2  } \right)\right] \Bigg\}  db\\
&   & 
+  \,  \phi (t)^{1-2\gamma }  \int_{ 0}^{t}  \phi (b) ^{1-2\gamma }
 \left(\phi (t)  - \phi (b) \right)  \left(\phi (t)  + \phi (b)\right)^{ 2  \gamma -2  } \\
 &  &
 \times \Bigg\{
\left[ 1+\left( \frac{  \phi (t)   \phi (b)  } {(\phi (t)  + \phi (b))^2  } \right)\right] 
 + \left( \frac{  \phi (t)   \phi (b)  } {(\phi (t)  + \phi (b))^2  } \right)^{ \gamma }\left[  1+  \left( \frac{  \phi (t)   \phi (b)  } {(\phi (t)  + \phi (b))^2  } \right)\right]\Bigg\} \, db   \,.
\end{eqnarray*}
This estimate
 proves that the   integral is a proper integral since $ \ell>-2$. Similarly we verify convergence of the integral with $A(x,\partial_x) $	in the integrand.	
 Proposition is proved. \hfill $\square$

\begin{proposition}
The integral of  the right hand side of (\ref{Tric_eq}), 
\begin{eqnarray*} 
 (\phi' (t))^2 \int_{ 0}^{t} 
 \left( \phi (t)  + \phi (b) \right)^{-2\gamma }  F \left(\gamma , \gamma ;1; \frac{(\phi (t)  - \phi
(b))^2  } {(\phi (t)  + \phi (b))^2  }  \right)  w_{r} (x ,0;b )\, db \,,
\end{eqnarray*}
defines a continuous function. 
\end{proposition}
\medskip

\noindent
{\bf Proof.} The proof is similar to the proof of the previous proposition, therefore we skip it.
\hfill $\square$\\

For $\varepsilon >0 $ consider  the regularization of (\ref{33}), that is the function
\begin{eqnarray}
\label{33epsilon}
u^{(\varepsilon)} (x,t)
& = &
c_\ell \int_{ \varepsilon }^{t} db
  \int_{ 0}^{ \phi (t)- \phi (b)} \alpha  (t,b,r)  F \left(\gamma , \gamma ;1; \beta  (t,b,r) \right)  w(x,r;b )
dr, \quad  t>\varepsilon \,. 
\end{eqnarray}
We set $u^{(\varepsilon)}(x,t) =0 $ if $t \in [0,\varepsilon ] $.  It is evident  that $ u^{(\varepsilon)} $ is continuous   and moreover, $ u^{(\varepsilon)}  \to u  $ in ${\mathcal D}'$ as $\varepsilon \to 0 $. 
It will be shown that $ u^{(\varepsilon)}  $ satisfies equation  
\begin{eqnarray*} 
&  &
 u^{(\varepsilon)}_{tt}-t^{\ell}A(x ,\partial_x) u^{(\varepsilon)} 
  =  
f(x ,t) \nonumber \\
&  &
\hspace{3cm} + c_\ell(\phi' (t))^2 \int_{ \varepsilon }^{t} 
 \left( \phi (t)  + \phi (b) \right)^{-2\gamma }  F \left(\gamma , \gamma ;1; \frac{(\phi (t)  - \phi
(b))^2  } {(\phi (t)  + \phi (b))^2  }  \right)  w_{r} (x ,0;b )\, db , \quad  t>\varepsilon \,.
\end{eqnarray*}
The right hand side of the  last equation converges to the right hand side of (\ref{Tric_eq}); consequently $u $  solves (\ref{Tric_eq}).
\smallskip

Calculations for $(\partial_t^2-t^{\ell}A(x ,\partial_x) ) u$ and for $ (\partial_t^2-t^{\ell}A(x ,\partial_x) )u^{(\varepsilon)} $
are almost identical; the distinction may be only in the terms, containing single integration in $b$. 
For the $(\partial_t^2-t^{\ell}A(x ,\partial_x) ) u$
the differentiation is formal because of the  singularity  may caused by negative powers of $b$ at $b=0$. 
For  $ (\partial_t^2-t^{\ell}A(x ,\partial_x) )u^{(\varepsilon)} $
those calculations are  just differentiation of the integrals smoothly depending on parameters, where $b> \varepsilon $. 
 Since calculations are identical, we give here only formal calculations for
$(\partial_t^2-t^{\ell}A(x ,\partial_x) ) u$. In fact, later on, it will be clear that the total contribution of all possibly divergent integrals is zero.
 
 \smallskip

Then, for the derivative of the function $u= u(x,t) $ we obtain
\begin{eqnarray*}
\partial_t u(x,t) 
& = &
c_\ell   \int_{ 0}^{t}
 \phi' (t)   \left(  4\phi (t)   \phi (b)     \right)^{-\gamma }  w(x,  \phi (t)- \phi (b) ;b ) \,db\\
&   &
+ c_\ell  \int_{ 0}^{t} db
  \int_{ 0}^{ \phi (t)- \phi (b)} 2\phi' (t)  (-\gamma ) (\phi (t)  + \phi (b))\left(  (\phi (t)  + \phi (b))^2  -r^2 \right)^{-\gamma -1} \\
&  & 
\hspace{2.5cm} \times F \left(\gamma , \gamma ;1; \beta  (t,b,r) \right)  w(x,r;b )
dr \\
&  &
+ c_\ell \int_{ 0}^{t} db
  \int_{ 0}^{ \phi (t)- \phi (b)} \left(  (\phi (t)  + \phi (b))^2  -r^2 \right)^{-\gamma } \left( \partial_t \beta  (t,b,r)\right)  F_z \left(\gamma , \gamma ;1; \beta  (t,b,r) \right)   w(x,r;b )\,
dr  \,.
\end{eqnarray*}

The second order derivative can be represented as follows:
\begin{eqnarray*} 
\partial_t^2 u(x,t)  
& = &
A+B+C \,,
\end{eqnarray*}
where we denoted
\begin{eqnarray*}
A
& := &
\partial_t \Bigg\{ c_\ell   \int_{ 0}^{t}
 \phi' (t)   \left(  4\phi (t)   \phi (b)     \right)^{-\gamma }  w(x,( \phi (t)- \phi (b));b ) \,db \Bigg\}\,,\\
B
& :=  &
\partial_t \Bigg\{  c_\ell  \int_{ 0}^{t} db
  \int_{ 0}^{ \phi (t)- \phi (b)} 2\phi' (t)  (-\gamma ) (\phi (t)  + \phi (b))\left(  (\phi (t)  + \phi (b))^2  -r^2 \right)^{-\gamma -1} \\
&  & 
\hspace{2.5cm} \times F \left(\gamma , \gamma ;1; \beta  (t,b,r) \right)  w(x,r;b )
dr  \Bigg\}\,,\\
C
& := &
\partial_t \Bigg\{ c_\ell \int_{ 0}^{t} db
  \int_{ 0}^{ \phi (t)- \phi (b)} \left(  (\phi (t)  + \phi (b))^2  -r^2 \right)^{-\gamma } \\
&  & 
\hspace{2.5cm} \times\left( \partial_t  \beta  (t,b,r) \right)  F_z \left(\gamma , \gamma ;1; \beta  (t,b,r) \right)  w(x,r;b )\,
dr  \Bigg\}\,.
\end{eqnarray*}
Then
\begin{eqnarray*}
\frac{A}{ c_\ell }
& = &
4 ^{-\gamma }   \phi'{}' (t)     \phi (t)^{-\gamma }  \int_{ 0}^{t}
   \phi (b) ^{-\gamma }  w(x,( \phi (t)- \phi (b));b ) \,db  \\ 
&  &
+ \, 4 ^{-\gamma }  (-\gamma ) (\phi' (t) )^{2} \phi (t)^{-\gamma-1 } \int_{ 0}^{t}
   \phi (b) ^{-\gamma }  w(x,( \phi (t)- \phi (b));b ) \,db \\
&  &
+ \,4 ^{-\gamma }   \phi' (t)     \phi (t)^{-\gamma }  
   \phi (t) ^{-\gamma }  w(x,0;t )  
+ \,4 ^{-\gamma }   (\phi' (t) )^2    \phi (t)^{-\gamma }  \int_{ 0}^{t}
   \phi (b) ^{-\gamma }  w_r(x,( \phi (t)- \phi (b));b ) \,db \,.
\end{eqnarray*}
The choice of $w(x,0;t )=f(x,t)$ implies
\begin{eqnarray*}
\frac{A}{ c_\ell }
& = &
4 ^{-\gamma }   \phi'{}' (t)     \phi (t)^{-\gamma }  \int_{ 0}^{t}
   \phi (b) ^{-\gamma }  w(x, \phi (t)- \phi (b) ;b ) \,db  \\ 
&  &
 -4 ^{-\gamma } \gamma   (\phi' (t) )^{2} \phi (t)^{-\gamma-1 } \int_{ 0}^{t}
   \phi (b) ^{-\gamma }  w(x,  \phi (t)- \phi (b) ;b ) \,db \\
&  &
+ \frac{1}{c_\ell } f(x,t) 
+ 4 ^{-\gamma }   (\phi' (t) )^2    \phi (t)^{-\gamma }  \int_{ 0}^{t}
   \phi (b) ^{-\gamma }  w_r(x,  \phi (t)- \phi (b) ;b ) \,db \,, 
\end{eqnarray*}
since
\begin{eqnarray}
&  &
4 ^{-\gamma } c_\ell   \phi' (t)     \phi (t)^{-2\gamma } = 1\,.
\end{eqnarray}

Further, 
\begin{eqnarray*}
\frac{B}{ c_\ell }
& =  &
 2 (-\gamma )  \phi'{}' (t)  \int_{ 0}^{t} db
  \int_{ 0}^{ \phi (t)- \phi (b)} (\phi (t)  + \phi (b))\left(  (\phi (t)  + \phi (b))^2  -r^2 \right)^{-\gamma -1} \\
&  & 
\times F \left(\gamma , \gamma ;1; \beta  (t,b,r) \right)   w(x,r;b )
dr  \\
&  &
+  2 (-\gamma )  ( \phi'  (t))^2  \int_{ 0}^{t} 
  (\phi (t)  + \phi (b))\left(  4\phi (t) \phi (b)  \right)^{-\gamma -1}  
  w(x,\phi (t)- \phi (b);b )\, db  \\
&  &
+ 2 (-\gamma )  (\phi'  (t))^2  \int_{ 0}^{t} db
  \int_{ 0}^{ \phi (t)- \phi (b)} \left(  (\phi (t)  + \phi (b))^2  -r^2 \right)^{-\gamma -1} \\&  & 
\times F \left(\gamma , \gamma ;1; \beta  (t,b,r) \right)   w(x,r;b )
dr  \\
&  &
+ 2 (-\gamma )(-\gamma -1)  (\phi'  (t))^2  \int_{ 0}^{t} db
  \int_{ 0}^{ \phi (t)- \phi (b)} 2(\phi (t)  + \phi (b))^2\left(  (\phi (t)  + \phi (b))^2  -r^2 \right)^{-\gamma -2} \\
&  & 
\times F \left(\gamma , \gamma ;1; \beta  (t,b,r) \right)   w(x,r;b )
dr  \\ 
&  &
+ 2 (-\gamma ) \phi'  (t)  \int_{ 0}^{t} db
  \int_{ 0}^{ \phi (t)- \phi (b)} (\phi (t)  + \phi (b))\left(  (\phi (t)  + \phi (b))^2  -r^2 \right)^{-\gamma -1} \\
&  & 
\times \left( \partial_t \beta  (t,b,r)  \right) F_z \left(\gamma , \gamma ;1; \beta  (t,b,r) \right)   w(x,r;b )
dr  \,.
\end{eqnarray*}
Application of Lemma~\ref{alphabeta} yields
\begin{eqnarray*}
\frac{B}{ c_\ell }
& =  &
 -2  \gamma   \phi'{}' (t)  \int_{ 0}^{t} db
  \int_{ 0}^{ \phi (t)- \phi (b)} (\phi (t)  + \phi (b))\left(  (\phi (t)  + \phi (b))^2  -r^2 \right)^{-\gamma -1} \\
&  & 
\times F \left(\gamma , \gamma ;1; \beta  (t,b,r) \right)  w(x,r;b )
dr  \\
&  &
  - 2  \gamma    ( \phi'  (t))^2  \int_{ 0}^{t} 
  (\phi (t)  + \phi (b))\left(   4\phi (t) \phi (b)  \right)^{-\gamma -1}  
  w(x,\phi (t)- \phi (b);b )\, db  \\
&  &
  - 2  \gamma    (\phi'  (t))^2  \int_{ 0}^{t} db
  \int_{ 0}^{ \phi (t)- \phi (b)} \left(  (\phi (t)  + \phi (b))^2  -r^2 \right)^{-\gamma -1} F \left(\gamma , \gamma ;1; \beta  (t,b,r) \right)  w(x,r;b )
dr  \\
&  &
+ 2   \gamma  ( \gamma +1)   (\phi'  (t))^2  \int_{ 0}^{t} db
  \int_{ 0}^{ \phi (t)- \phi (b)} 2(\phi (t)  + \phi (b))^2\left(  (\phi (t)  + \phi (b))^2  -r^2 \right)^{-\gamma -2} \\
&  & 
\times F \left(\gamma , \gamma ;1; \beta  (t,b,r) \right)  w(x,r;b )
dr  \\ 
&  &
- 2  \gamma    \phi'  (t)  \int_{ 0}^{t} db
  \int_{ 0}^{ \phi (t)- \phi (b)} (\phi (t)  + \phi (b))\left(  (\phi (t)  + \phi (b))^2  -r^2 \right)^{-\gamma -1} \\
&  & 
\times \left( 4\phi' (t)   \phi (b) \frac{ \phi^2 (t)  - \phi^2 (b)    + 
   r^2 } {((\phi (t) + \phi (b))^2 - r^2)^2}   \right) F_z \left(\gamma , \gamma ;1; \beta  (t,b,r) \right)  w(x,r;b )
dr  \,.
\end{eqnarray*}
Consequently,
\begin{eqnarray*}
\frac{B}{ c_\ell }
& =  &
- 2  \gamma   \phi'{}' (t)  \int_{ 0}^{t} db
  \int_{ 0}^{ \phi (t)- \phi (b)} (\phi (t)  + \phi (b))\left(  (\phi (t)  + \phi (b))^2  -r^2 \right)^{-\gamma -1} \\
&  & 
\times F \left(\gamma , \gamma ;1; \beta  (t,b,r) \right)   w(x,r;b )
dr  \\
&  &
 - 2  \gamma  ( \phi'  (t))^2  \int_{ 0}^{t} 
  (\phi (t)  + \phi (b))\left(   4\phi (t) \phi (b)  \right)^{-\gamma -1}  
  w(x,\phi (t)- \phi (b);b )\, db  \\
&  &
 - 2  \gamma  (\phi'  (t))^2  \int_{ 0}^{t} db
  \int_{ 0}^{ \phi (t)- \phi (b)} \left(  (\phi (t)  + \phi (b))^2  -r^2 \right)^{-\gamma -1}  F \left(\gamma , \gamma ;1; \beta  (t,b,r) \right)   w(x,r;b ) \,dr  \\
&  &
+ 2  \gamma (\gamma +1)  (\phi'  (t))^2  \int_{ 0}^{t} db
  \int_{ 0}^{ \phi (t)- \phi (b)} 2(\phi (t)  + \phi (b))^2\left(  (\phi (t)  + \phi (b))^2  -r^2 \right)^{-\gamma -2} \\
&  & 
\times F \left(\gamma , \gamma ;1; \beta  (t,b,r) \right)   w(x,r;b )\,dr  \\ 
&  &
-8 \gamma  (\phi'  (t) )^2 \int_{ 0}^{t} db
  \int_{ 0}^{ \phi (t)- \phi (b)} (\phi (t)  + \phi (b))\left(  (\phi (t)  + \phi (b))^2  -r^2 \right)^{-\gamma -3} \\
&  & 
\times     \phi (b)  ( \phi^2 (t)  - \phi^2 (b)    +    r^2)    
F_z \left(\gamma , \gamma ;1; \beta  (t,b,r) \right)   w(x,r;b ) \, dr  \,.
\end{eqnarray*}
Finally, we derive 
\begin{eqnarray*}
\frac{B}{ c_\ell }
& =  &
 -2  \gamma    \phi'{}' (t)  \int_{ 0}^{t} db
  \int_{ 0}^{ \phi (t)- \phi (b)} (\phi (t)  + \phi (b))\left(  (\phi (t)  + \phi (b))^2  -r^2 \right)^{-\gamma -1} \\
&  & 
\times F \left(\gamma , \gamma ;1; \beta  (t,b,r) \right)   w(x,r;b )
dr  \\
&  &
   -2  \gamma  ( \phi'  (t))^2  \int_{ 0}^{t} 
  (\phi (t)  + \phi (b))\left( 4\phi (t)  \phi (b)  \right)^{-\gamma -1}  
  w(x,\phi (t)- \phi (b);b )\, db  \\
&  &
-2  \gamma   (\phi'  (t))^2  \int_{ 0}^{t} db
  \int_{ 0}^{ \phi (t)- \phi (b)} \left(  (\phi (t)  + \phi (b))^2  -r^2 \right)^{-\gamma -1}  F \left(\gamma , \gamma ;1; \beta  (t,b,r) \right)   w(x,r;b ) \,dr   \\
&  &
+ 2  \gamma ( \gamma +1)  ( \phi'  (t))^2  \int_{ 0}^{t} db
  \int_{ 0}^{ \phi (t)- \phi (b)} 2(\phi (t)  + \phi (b))^2\left(  (\phi (t)  + \phi (b))^2  -r^2 \right)^{-\gamma -2} \\
&  & 
\times F \left(\gamma , \gamma ;1; \beta  (t,b,r) \right)   w(x,r;b )\,dr   \\ 
&  &
-8 \gamma  (\phi'  (t) )^2 \int_{ 0}^{t} db
  \int_{ 0}^{ \phi (t)- \phi (b)} (\phi (t)  + \phi (b))\left(  (\phi (t)  + \phi (b))^2  -r^2 \right)^{-\gamma -3} \\
&  & 
\times     \phi (b)  ( \phi^2 (t)  - \phi^2 (b)    +    r^2)    
F_z \left(\gamma , \gamma ;1; \beta  (t,b,r) \right)   w(x,r;b ) \, dr   \,.
\end{eqnarray*}
Next we consider $C/c_\ell $ and by means of Lemma~\ref{alphabeta} we obtain
\begin{eqnarray*}
\frac{C}{c_\ell}
& = &
  \phi' (t) \int_{ 0}^{t} db
   \left(  4\phi (t) \phi (b) \right)^{-\gamma } \\
&  & 
\times\left( 4\phi' (t)   \phi (b) 
\frac{ 2\phi^2 (t)  - 2\phi (t)   \phi (b)  } {(4\phi (t) \phi (b) )^2}  \right)  F_z \left(\gamma , \gamma ;1; 0 \right)  w(x,r;b )\,dr  \\
&  &
+ (-\gamma )\phi' (t)\int_{ 0}^{t} db
  \int_{ 0}^{ \phi (t)- \phi (b)} 2(\phi (t)  + \phi (b))\left(  (\phi (t)  + \phi (b))^2  -r^2 \right)^{-\gamma -1} \\
&  & 
\times\left( 4\phi' (t)   \phi (b) \frac{ \phi^2 (t)  - \phi^2 (b)    + 
   r^2 } {((\phi (t) + \phi (b))^2 - r^2)^2}   \right)  F_z \left(\gamma , \gamma ;1; \beta  (t,b,r) \right)  w(x,r;b )\, dr \\
&  &
+ \int_{ 0}^{t} db
  \int_{ 0}^{ \phi (t)- \phi (b)} \left(  (\phi (t)  + \phi (b))^2  -r^2 \right)^{-\gamma } \\
&  & 
\times \left[ \partial_t \left( 4\phi' (t)   \phi (b) \frac{ \phi^2 (t)  - \phi^2 (b)    + 
   r^2 } {((\phi (t) + \phi (b))^2 - r^2)^2}   \right)\right]  F_z \left(\gamma , \gamma ;1; \beta  (t,b,r) \right)  w(x,r;b )\, dr \\
&  &
+ \int_{ 0}^{t} db
  \int_{ 0}^{ \phi (t)- \phi (b)} \left(  (\phi (t)  + \phi (b))^2  -r^2 \right)^{-\gamma } \\
&  & 
\times\left( 4\phi' (t)   \phi (b) \frac{ \phi^2 (t)  - \phi^2 (b)    + 
   r^2 } {((\phi (t) + \phi (b))^2 - r^2)^2}   \right)^2 F_{zz}  \left(\gamma , \gamma ;1; \beta  (t,b,r) \right)   w(x,r;b )\, dr \,.
\end{eqnarray*}
If we take into account the following easily verified identity
\begin{eqnarray*}
&  &
\partial_t \left( 4\phi' (t)   \phi (b) \frac{ \phi^2 (t)  - \phi^2 (b)    +     r^2 } {((\phi (t) + \phi (b))^2 - r^2)^2}   \right) \\
& = &
  4\phi'{}' (t)   \phi (b) \frac{ \phi^2 (t)  - \phi^2 (b)    + 
   r^2 } {((\phi (t) + \phi (b))^2 - r^2)^2}   \\
   & &
   +\, 8(\phi' (t))^2   \phi (b) \frac{  - \phi^3  (t)  +3\phi (t) \phi (b)^2 +2 \phi^3  (b) -3\phi (t)r^2-2 \phi (b)r^2 } {((\phi (t) + \phi (b))^2 - r^2)^3} \,,
\end{eqnarray*}
then we can write
\begin{eqnarray*}
\frac{C}{c_\ell}
&  = &
  \phi' (t) \int_{ 0}^{t} db
   \left( 4\phi (t)  \phi (b)  \right)^{-\gamma } \left(  \phi' (t)    
\frac{  \phi  (t)  -    \phi (b)  } {2 \phi (t)\phi (b)}  \right)  F_z \left(\gamma , \gamma ;1; 0 \right)  w(x,(\phi (t)- \phi (b));b )  \\
&  &
 -\gamma \phi' (t)\int_{ 0}^{t} db
  \int_{ 0}^{ \phi (t)- \phi (b)} 2(\phi (t)  + \phi (b))\left(  (\phi (t)  + \phi (b))^2  -r^2 \right)^{-\gamma -1} \\
&  & 
\times\left( 4\phi' (t)   \phi (b) \frac{ \phi^2 (t)  - \phi^2 (b)    + 
   r^2 } {((\phi (t) + \phi (b))^2 - r^2)^2}   \right)  F_z \left(\gamma , \gamma ;1; \beta  (t,b,r) \right)  w(x,r;b )\, dr \\
&  &
+ \int_{ 0}^{t} db
  \int_{ 0}^{ \phi (t)- \phi (b)} \left(  (\phi (t)  + \phi (b))^2  -r^2 \right)^{-\gamma } \Bigg[  4\phi'{}' (t)   \phi (b) \frac{ \phi^2 (t)  - \phi^2 (b)    + 
   r^2 } {((\phi (t) + \phi (b))^2 - r^2)^2} \\
   & &
   + \,8(\phi' (t))^2   \phi (b) \frac{  - \phi^3  (t)  +3\phi (t) \phi (b)^2 +2 \phi^3  (b) -3\phi (t)r^2-2 \phi (b)r^2 } {((\phi (t) + \phi (b))^2 - r^2)^3} \Bigg]  \\
&  &
\times F_z \left(\gamma , \gamma ;1; \beta  (t,b,r) \right)  w(x,r;b )\, dr \\
&  &
+ \int_{ 0}^{t} db
  \int_{ 0}^{ \phi (t)- \phi (b)} \left(  (\phi (t)  + \phi (b))^2  -r^2 \right)^{-\gamma } \\
&  & 
\times\left( 4\phi' (t)   \phi (b) \frac{ \phi^2 (t)  - \phi^2 (b)    + 
   r^2 } {((\phi (t) + \phi (b))^2 - r^2)^2}   \right)^2  F_{zz} \left(\gamma , \gamma ;1; \beta  (t,b,r) \right)  w(x,r;b )\, dr \,.
\end{eqnarray*}
Next we consider 
\begin{eqnarray*} 
\frac{1}{c_\ell} A(x,\partial_x) u 
& =  &
 \int_{ 0}^{t} db
  \int_{ 0}^{ \phi (t)- \phi (b)}  \left(  (\phi (t)  + \phi (b))^2  -r^2 \right)^{-\gamma }  F \left(\gamma , \gamma ;1; \beta  (t,b,r) \right)  A(x,\partial_x)w(x,r;b )\,dr \\
& =  &
\int_{ 0}^{t} db
  \int_{ 0}^{ \phi (t)- \phi (b)}  \left(  (\phi (t)  + \phi (b))^2  -r^2 \right)^{-\gamma }  F \left(\gamma , \gamma ;1; \beta  (t,b,r) \right)  w_{rr} (x,r;b )\,dr \,,
\end{eqnarray*}
where we have used the relation (\ref{8}). Integrating by parts we obtain
\begin{eqnarray*} 
\frac{1}{c_\ell} A(x,\partial_x) u & =  &
 \int_{ 0}^{t}
 \left( 4\phi (t) \phi (b) \right)^{-\gamma }   w_{r} (x,\phi (t)  - \phi (b);b )  \, db\\
 &  &
 - \int_{ 0}^{t} 
 \left(  (\phi (t)  + \phi (b))^2   \right)^{-\gamma }  F \left(\gamma , \gamma ;1; \beta  (t,b,r) \right)   w_{r} (x,0;b ) \, db\\
&  &
-  \int_{ 0}^{t} db
  \int_{ 0}^{ \phi (t)- \phi (b)}  \Big[ \partial_r \left(  (\phi (t)  + \phi (b))^2  -r^2 \right)^{-\gamma }  F \left(\gamma , \gamma ;1; \beta  (t,b,r) \right)\Big]  w_{ r} (x,r;b )\,dr \,.
\end{eqnarray*}
The first term of the last equation will be canceled with the similar term of $C/c_\ell$. We consider  now the last term and transform it as follows
\begin{eqnarray*} 
II & := &
-  \int_{ 0}^{t} db
  \int_{ 0}^{ \phi (t)- \phi (b)}  \Big[ \partial_r \left(  (\phi (t)  + \phi (b))^2  -r^2 \right)^{-\gamma }  F \left(\gamma , \gamma ;1; \beta  (t,b,r) \right)\Big]  w_{ r} (x,r;b )\,dr \\
& = &
-   \int_{ 0}^{t} db
  \int_{ 0}^{ \phi (t)- \phi (b)} (2r\gamma ) \left(  (\phi (t)  + \phi (b))^2  -r^2 \right)^{-\gamma-1 }  F \left(\gamma , \gamma ;1; \beta  (t,b,r) \right)  w_{ r} (x,r;b )\,dr \\
& &
-  \int_{ 0}^{t} db
  \int_{ 0}^{ \phi (t)- \phi (b)}  \left(  (\phi (t)  + \phi (b))^2  -r^2 \right)^{-\gamma }  \Big[ \partial_r \beta  (t,b,r)  \Big] F_z \left(\gamma , \gamma ;1; \beta  (t,b,r) \right)  w_{ r} (x,r;b )\,dr\,.
\end{eqnarray*}
Then 
we use Lemma~\ref{alphabeta} and the integration by parts and obtain
\begin{eqnarray*} 
I I
& = &
- 2\gamma \int_{ 0}^{t} db
  \int_{ 0}^{ \phi (t)- \phi (b)}  r  \left(  (\phi (t)  + \phi (b))^2  -r^2 \right)^{-\gamma-1 }  F \left(\gamma , \gamma ;1; \beta  (t,b,r) \right)  w_{ r} (x,r;b )\,dr \\
& &
+ \, 8  \phi (t)\int_{ 0}^{t} db\, \phi (b) 
  \int_{ 0}^{ \phi (t)- \phi (b)}  \left(  (\phi (t)  + \phi (b))^2  -r^2 \right)^{-\gamma -2}    r  F_z \left(\gamma , \gamma ;1; \beta  (t,b,r) \right)  w_{ r} (x,r;b )\,dr\,.
\end{eqnarray*}
One more integration by parts gives
\begin{eqnarray*}
II & = &
- 2\gamma  \int_{ 0}^{t}
  (\phi (t)- \phi (b))  \left(  4\phi (t)   \phi (b)   \right)^{-\gamma-1 }    w  (x,(\phi (t)- \phi (b)) ;b )\,  db \\
&  &
+ 2\gamma  \int_{ 0}^{t} db
  \int_{ 0}^{ \phi (t)- \phi (b)} \left( \partial_r  r  \left(  (\phi (t)  + \phi (b))^2  -r^2 \right)^{-\gamma-1 }  F \left(\gamma , \gamma ;1; \beta  (t,b,r) \right) \right) w  (x,r;b )\,dr \\
& &
+  8 \phi (t)\int_{ 0}^{t} \, \phi (b) 
   \left( 4\phi (t)   \phi (b) \right)^{-\gamma -2}   (\phi (t)  - \phi (b))  F_z \left(\gamma , \gamma ;1; 0 \right)  w  (x,\phi (t)- \phi (b);b )\, db\\
& &
- 8   \phi (t)\int_{ 0}^{t} db\, \phi (b) 
  \int_{ 0}^{ \phi (t)- \phi (b)} \left( \partial_r  \left(  (\phi (t)  + \phi (b))^2  -r^2 \right)^{-\gamma -2}    r   F_z \left(\gamma , \gamma ;1; \beta  (t,b,r) \right)\right)  w  (x,r;b )\,dr \,.
\end{eqnarray*}  
It is easily seen that 
\begin{eqnarray*} 
II  
& = &
- 2\gamma   \int_{ 0}^{t}
  (\phi (t)- \phi (b))  \left(  4\phi (t)   \phi (b) ) ^2 \right)^{-\gamma-1 }   w  (x,(\phi (t)- \phi (b)) ;b )\,  db \\
&  &
+ 2\gamma   \int_{ 0}^{t} db
  \int_{ 0}^{ \phi (t)- \phi (b)}    \left(  (\phi (t)  + \phi (b))^2  -r^2 \right)^{-\gamma-1 }  F \left(\gamma , \gamma ;1; \beta  (t,b,r) \right)  w  (x,r;b )\,dr \\
&  &
+ 4\gamma (\gamma +1)   \int_{ 0}^{t} db
  \int_{ 0}^{ \phi (t)- \phi (b)}  r^2  \left(  (\phi (t)  + \phi (b))^2  -r^2 \right)^{-\gamma-2 }  F \left(\gamma , \gamma ;1; \beta  (t,b,r) \right) w  (x,r;b )\,dr \\
&  &
+ 2\gamma  \int_{ 0}^{t} db
  \int_{ 0}^{ \phi (t)- \phi (b)}  r  \left(  (\phi (t)  + \phi (b))^2  -r^2 \right)^{-\gamma-1 } \left( - 8r\frac{   \phi (t)  \phi (b) } {((\phi (t)  + \phi (b))^2 - r^2)^2} \right)\\
&  &
\hspace{4cm}   \times F_z \left(\gamma , \gamma ;1; \beta  (t,b,r) \right) w  (x,r;b )\,dr \\
& &
+  8 \phi (t)\int_{ 0}^{t} \, \phi (b) 
   \left( 4\phi (t)   \phi (b) \right)^{-\gamma -2}   (\phi (t)  - \phi (b))  F_z \left(\gamma , \gamma ;1; 0 \right)  w  (x,\phi (t)- \phi (b);b )\, db\\
& &
- 8  \phi (t)\int_{ 0}^{t} db\, \phi (b) 
  \int_{ 0}^{ \phi (t)- \phi (b)}   \left(  (\phi (t)  + \phi (b))^2  -r^2 \right)^{-\gamma -2}   F_z \left(\gamma , \gamma ;1; \beta  (t,b,r) \right)   w  (x,r;b )\,dr \\
& &
- 8  \phi (t)\int_{ 0}^{t} db\, \phi (b) 
  \int_{ 0}^{ \phi (t)- \phi (b)}    2(\gamma +2) \left(  (\phi (t)  + \phi (b))^2  -r^2 \right)^{-\gamma -3}    r^2 \\
  &  &
\hspace{4cm}   \times   F_z \left(\gamma , \gamma ;1; \beta  (t,b,r) \right)   w  (x,r;b )\,dr\\
& &
- 8  \phi (t)\int_{ 0}^{t} db\, \phi (b) 
  \int_{ 0}^{ \phi (t)- \phi (b)} \left(  (\phi (t)  + \phi (b))^2  -r^2 \right)^{-\gamma -2}    r \left( - 8r\frac{   \phi (t)  \phi (b) } {((\phi (t)  + \phi (b))^2 - r^2)^2} \right)   \\
  &  &
\hspace{4cm}  \times 
F_{zz} \left(\gamma , \gamma ;1; \beta  (t,b,r) \right) w  (x,r;b )\,dr\,.
\end{eqnarray*}
Then we take into account the properties
$
F  \left(\gamma , \gamma ;1; 0 \right)=1$, $F_z \left(\gamma , \gamma ;1; 0 \right)=\gamma ^2$ of the hypergeometric function
and obtain

\begin{eqnarray*} 
II 
& = &
\int_{ 0}^{t}
  (\phi (t)- \phi (b))  \left(  4\phi (t)   \phi (b)  \right)^{-\gamma-1 } \Big\{  - 2\gamma   
+    2 \gamma^2   \Big\} w  (x,\phi (t)- \phi (b);b )\, db\\
&  &
+   \int_{ 0}^{t} db
  \int_{ 0}^{ \phi (t)- \phi (b)}    \left(  (\phi (t)  + \phi (b))^2  -r^2 \right)^{-\gamma-2 } \Big\{ 2\gamma \left(  (\phi (t)  + \phi (b))^2  -r^2 \right) 
+    4 \gamma (\gamma +1) r^2 
\Big\} \\
&  &
\hspace{4cm}  \times F \left(\gamma , \gamma ;1; \beta  (t,b,r) \right) w  (x,r;b )\,dr \\
&  &
+  \int_{ 0}^{t} db
  \int_{ 0}^{ \phi (t)- \phi (b)}  \left(  (\phi (t)  + \phi (b))^2  -r^2 \right)^{-\gamma-1 }
\Bigg\{-16\gamma r^2    \frac{   \phi (t)  \phi (b) } {((\phi (t)  + \phi (b))^2 - r^2)^2}  \\
& &
- 8 \phi (t) \phi (b) 
     \left(  (\phi (t)  + \phi (b))^2  -r^2 \right)^{-1}   
- 8  \phi (t)\phi (b) 
   2(\gamma +2) \left(  (\phi (t)  + \phi (b))^2  -r^2 \right)^{-2}    r^2  \Bigg\}\\
  &  &
\hspace{4cm}  \times   F_z \left(\gamma , \gamma ;1; \beta  (t,b,r) \right)   w  (x,r;b )\,dr \\
& &
+\int_{ 0}^{t} db\,   
  \int_{ 0}^{ \phi (t)- \phi (b)} \left(  (\phi (t)  + \phi (b))^2  -r^2 \right)^{-\gamma -2}     \frac{     64 (\phi (t))^2  (\phi (b))^2 r^2} {((\phi (t)  + \phi (b))^2 - r^2)^2}     \\
  &  &
\hspace{4cm}  \times 
F_{zz} \left(\gamma , \gamma ;1; \beta  (t,b,r) \right) w  (x,r;b )\,dr \,.
\end{eqnarray*} 
Consequently,   for \,  $\dsp -\frac{1}{c_\ell}t^{\ell} A(x,\partial_x) u $ \, we have
\begin{eqnarray*} 
&  &
-\frac{1}{c_\ell}t^{\ell} A(x,\partial_x) u \\
& =  &
- (\phi' (t))^2 \int_{ 0}^{t}
 \left(  4\phi (t) \phi (b)  \right)^{-\gamma }     w_{r} (x,\phi (t)  - \phi (b);b )   \, db \\
 &  &
 + (\phi' (t))^2 \int_{ 0}^{t} 
 \left(  \phi (t)  + \phi (b)  \right)^{-2\gamma }  F \left(\gamma , \gamma ;1; \frac{(\phi (t)  - \phi
(b))^2  } {(\phi (t)  + \phi (b))^2  }  \right)   w_{r} (x,0;b ) \, db\\
&  &
+ (\phi' (t))^2 \int_{ 0}^{t}
  (\phi (t)- \phi (b))  \left(  4\phi (t)   \phi (b)  \right)^{-\gamma-1 } 2\gamma (\gamma -1) w  (x,\phi (t)- \phi (b);b )\, db \\
& &
-   \int_{ 0}^{t} db
  \int_{ 0}^{ \phi (t)- \phi (b)}    \left(  (\phi (t)  + \phi (b))^2  -r^2 \right)^{-\gamma-2 } (\phi' (t))^2\Big\{ 2\gamma \left(  (\phi (t)  + \phi (b))^2  -r^2 \right) 
+    4 \gamma (\gamma +1) r^2 
\Big\} \\
&  &
\hspace{4cm}  \times F \left(\gamma , \gamma ;1; \beta  (t,b,r) \right) w  (x,r;b )\,dr \\
&  &
-\int_{ 0}^{t} db
  \int_{ 0}^{ \phi (t)- \phi (b)}  \left(  (\phi (t)  + \phi (b))^2  -r^2 \right)^{-\gamma-1 }(\phi' (t))^2 
\Bigg\{-16\gamma r^2    \frac{   \phi (t)  \phi (b) } {((\phi (t)  + \phi (b))^2 - r^2)^2}  \\
&   &
- 8 \phi (t) \phi (b) 
     \left(  (\phi (t)  + \phi (b))^2  -r^2 \right)^{-1}   - 16 \phi (t)\phi (b) 
 (\gamma +2) \left(  (\phi (t)  + \phi (b))^2  -r^2 \right)^{-2}    r^2  \Bigg\}\\
  &  &
\hspace{4cm}  \times   F_z \left(\gamma , \gamma ;1; \beta  (t,b,r) \right)   w  (x,r;b )\,dr \\
& &
- \int_{ 0}^{t} db\,   
  \int_{ 0}^{ \phi (t)- \phi (b)} \left(  (\phi (t)  + \phi (b))^2  -r^2 \right)^{-\gamma -2}    (\phi' (t))^2 \frac{     64 (\phi (t))^2  (\phi (b))^2 r^2} {((\phi (t)  + \phi (b))^2 - r^2)^2}     \\
  &  &
\hspace{4cm}  \times 
F_{zz} \left(\gamma , \gamma ;1; \beta  (t,b,r) \right) w  (x,r;b )\,dr \,.
\end{eqnarray*}

\noindent
In the double integrals the  terms with \, $F \left(\gamma , \gamma ;1; \beta  (t,b,r) \right)$\, in the expression of $\,\,u_{tt}-t^\ell A(x,\partial_x)u\,\, $  are:
\begin{eqnarray*}
&  & 
 \int_{ 0}^{t} db
  \int_{ 0}^{ \phi (t)- \phi (b)} \Bigg[\left(  (\phi (t)  + \phi (b))^2  -r^2 \right)^{-\gamma -1} 2 (-\gamma )  \phi'{}' (t) (\phi (t)  + \phi (b))\\
&  &
- \left(  (\phi (t)  + \phi (b))^2  -r^2 \right)^{-\gamma -1} 2\gamma  (\phi'  (t))^2 
+ \left(  (\phi (t)  + \phi (b))^2  -r^2 \right)^{-\gamma -2} 4  \gamma (\gamma +1)  ( \phi'  (t))^2 (\phi (t)  + \phi (b))^2\\
&  &
-     \left(  (\phi (t)  + \phi (b))^2  -r^2 \right)^{-\gamma-2 } (\phi' (t))^2\Bigg\{ 2\gamma \left(  (\phi (t)  + \phi (b))^2  -r^2 \right) 
+    4 \gamma (\gamma +1) r^2 
\Bigg\}\Bigg] \\
&  &
\hspace{4cm}  \times F \left(\gamma , \gamma ;1; \beta  (t,b,r) \right) w  (x,r;b )\,dr \,. 
\end{eqnarray*}

\noindent
In the double integrals the  terms with \, $F_{z }\left(\gamma , \gamma ;1; \beta  (t,b,r) \right)$\, in the expression of $\,\, u_{tt}-t^\ell A(x,\partial_x)u \,\, $ are: 
\begin{eqnarray*}
&  &
\int_{ 0}^{t} db
  \int_{ 0}^{ \phi (t)- \phi (b)} \Bigg[ \left(  (\phi (t)  + \phi (b))^2  -r^2 \right)^{-\gamma -3} (\phi (t)  + \phi (b)) 16(-\gamma )   (\phi'  (t) )^2  \phi (b)  ( \phi^2 (t)  - \phi^2 (b)    +    r^2)    
\\
&  &
+ \left(  (\phi (t)  + \phi (b))^2  -r^2 \right)^{-\gamma } \\
&  & 
\times \Bigg[  4\phi'{}' (t)   \phi (b) \frac{ \phi^2 (t)  - \phi^2 (b)    +    r^2 } {((\phi (t) + \phi (b))^2 - r^2)^2}   
   + 8(\phi' (t))^2   \phi (b) \frac{  - \phi^3  (t)  +3\phi (t) \phi (b)^2 +2 \phi^3  (b) -3\phi (t)r^2-2 \phi (b)r^2 } {((\phi (t) + \phi (b))^2 - r^2)^3} \Bigg]  \\
&  &
- \left(  (\phi (t)  + \phi (b))^2  -r^2 \right)^{-\gamma-1 }(\phi' (t))^2 
8\Bigg\{-2\gamma r^2    \frac{   \phi (t)  \phi (b) } {((\phi (t)  + \phi (b))^2 - r^2)^2}  -  \phi (t)\phi (b) 
     \left(  (\phi (t)  + \phi (b))^2  -r^2 \right)^{-1}  \\
&   &
 - 2  \phi (t)\phi (b) 
    (\gamma +2) \left(  (\phi (t)  + \phi (b))^2  -r^2 \right)^{-2}    r^2  \Bigg\}\Bigg]  F_z \left(\gamma , \gamma ;1; \beta  (t,b,r) \right)   w  (x,r;b )\,dr\,. 
\end{eqnarray*}
In the double integrals the terms with \,$F_{zz}\left(\gamma , \gamma ;1; \beta  (t,b,r) \right)$\,  in the expression of $\,\,u_{tt}-t^\ell A(x,\partial_x)u \,\,$  are:
\begin{eqnarray*} 
&  &
\int_{ 0}^{t} db
  \int_{ 0}^{ \phi (t)- \phi (b)} \left(  (\phi (t)  + \phi (b))^2  -r^2 \right)^{-\gamma } \left( 4\phi' (t)   \phi (b) \frac{ \phi^2 (t)  - \phi^2 (b)    + 
   r^2 } {((\phi (t) + \phi (b))^2 - r^2)^2}   \right)^2  \\
&  &
\hspace{4cm} \times F_{zz} \left(\gamma , \gamma ;1; \beta  (t,b,r) \right)  w(x,r;b )\, dr \\
& &
- \int_{ 0}^{t} db\,   
  \int_{ 0}^{ \phi (t)- \phi (b)} \left(  (\phi (t)  + \phi (b))^2  -r^2 \right)^{-\gamma -2}    (\phi' (t))^2 \frac{     64 (\phi (t))^2  (\phi (b))^2 r^2} {((\phi (t)  + \phi (b))^2 - r^2)^2}     \\
  &  &
\hspace{4cm}  \times 
F_{zz} \left(\gamma , \gamma ;1; \beta  (t,b,r) \right) w  (x,r;b )\,dr \\
& = &
\int_{ 0}^{t} db
  \int_{ 0}^{ \phi (t)- \phi (b)} \left(  (\phi (t)  + \phi (b))^2  -r^2 \right)^{-\gamma-4 } \\
&  & 
\times 16 (\phi' (t))^2(\phi (b))^2\Bigg\{    \left( \phi^2 (t)  - \phi^2 (b)    + 
   r^2 \right)^2   
-          4 (\phi (t))^2  r^2  \Bigg\}   
F_{zz} \left(\gamma , \gamma ;1; \beta  (t,b,r) \right) w  (x,r;b )\,dr \,.
\end{eqnarray*}
Now we denote by $I,J,Y$ the coefficients in the integrands in the corresponding to $F,F_z,F_{zz} $ terms. Then
\begin{eqnarray*} 
I
& := &
\left(  (\phi (t)  + \phi (b))^2  -r^2 \right)^{-\gamma -1} 2 (-\gamma )  \phi'{}' (t) (\phi (t)  + \phi (b))\\
&  &
+   \left(  (\phi (t)  + \phi (b))^2  -r^2 \right)^{-\gamma -1} 2 (-\gamma )   (\phi'  (t))^2 \\
&  &
+ \left(  (\phi (t)  + \phi (b))^2  -r^2 \right)^{-\gamma -2} 4  \gamma (\gamma +1)  ( \phi'  (t))^2 (\phi (t)  + \phi (b))^2\\
&  &
-     \left(  (\phi (t)  + \phi (b))^2  -r^2 \right)^{-\gamma-2 } (\phi' (t))^2\Bigg\{ 2\gamma \left(  (\phi (t)  + \phi (b))^2  -r^2 \right) 
+    4 \gamma (\gamma +1) r^2 
\Bigg\} \\ 
& = &
- 2\gamma   \phi'{}' (t) \phi (b) \left(  (\phi (t)  + \phi (b))^2  -r^2 \right)^{-\gamma -1} \,,
\end{eqnarray*}
which coincides with (\ref{I}). Further, after simple calculations we obtain
\begin{eqnarray*} 
J
& := &
 \left(  (\phi (t)  + \phi (b))^2  -r^2 \right)^{-\gamma -3} (\phi (t)  + \phi (b)) 16(-\gamma )   (\phi'  (t) )^2  \phi (b)  ( \phi^2 (t)  - \phi^2 (b)    +    r^2)    
\\
&  &
+ \left(  (\phi (t)  + \phi (b))^2  -r^2 \right)^{-\gamma }  \Bigg[  4\phi'{}' (t)   \phi (b) \frac{ \phi^2 (t)  - \phi^2 (b)    +    r^2 } {((\phi (t) + \phi (b))^2 - r^2)^2} \\  
&  &
   + 8(\phi' (t))^2   \phi (b) \frac{  - \phi^3  (t)  +3\phi (t) \phi (b)^2 +2 \phi^3  (b) -3\phi (t)r^2-2 \phi (b)r^2 } {((\phi (t) + \phi (b))^2 - r^2)^3} \Bigg]  \\
&  &
- \left(  (\phi (t)  + \phi (b))^2  -r^2 \right)^{-\gamma-1 }(\phi' (t))^2 
8\Bigg\{-2\gamma r^2    \frac{   \phi (t)  \phi (b) } {((\phi (t)  + \phi (b))^2 - r^2)^2}  \\
&   &
-  \phi (t)\phi (b) 
     \left(  (\phi (t)  + \phi (b))^2  -r^2 \right)^{-1}   - 2  \phi (t)\phi (b) 
    (\gamma +2) \left(  (\phi (t)  + \phi (b))^2  -r^2 \right)^{-2}    r^2  \Bigg\}\\
& = &
-4  \phi '{}'(t)  \phi (b)  \left(  (\phi (t)  + \phi (b))^2  -r^2 \right)^{-\gamma -2} 
\left[ \frac{ -2   \phi (b)   \phi (t) }{\gamma  }+     (\phi (t)-\phi(b))^2  - r^2 \right]\,,
\end{eqnarray*}
which coincides with (\ref{J}).
Similarly, we derive 
\begin{eqnarray*} 
Y
& := &
 16 \left(  (\phi (t)  + \phi (b))^2  -r^2 \right)^{-\gamma-4 } (\phi' (t))^2(\phi (b))^2\Bigg\{    \left( \phi^2 (t)  - \phi^2 (b)    + 
   r^2 \right)^2   
-          4 (\phi (t))^2  r^2  \Bigg\} \\
& = &
 16   \phi^2 (b) (\phi' (t))^2   \left(  (\phi (t)  + \phi (b))^2  -r^2 \right)^{-\gamma -3} 
\left(  (\phi (t)  - \phi (b))^2  -r^2 \right)    
\end{eqnarray*}
that coincides with (\ref{Y}). Lemma~\ref{L2.3} implies 
\begin{eqnarray*} 
&  &
Y(t,b,r;z)F_{zz} +   J(t,b,r;z)F_z+  I(t,b,r;z )F  \\
& =  &
    G(t,b,r;\gamma ) \Bigg\{ z(1-z)F_{zz} (\gamma , \gamma ;1; z  ) +[ 1-(2\gamma +1)z] F_z  (\gamma , \gamma ;1; z  )
-\gamma ^2 F  (\gamma , \gamma ;1; z  ) \Bigg\}  =0 \,,
\end{eqnarray*}
where $\dsp  z:=\beta (t,b,r)$  \,.

In the  expression of $\,\,u_{tt}-t^\ell A(x,\partial_x)u,\,\, $ we split into two parts the terms  which contain  only single integral   (without factor $c_\ell $):
\begin{eqnarray*} 
&  &
\frac{1}{c_\ell }  f(x,t)
 + (\phi' (t))^2 \int_{ 0}^{t} 
 \left(  \phi (t)  + \phi (b)   \right)^{-2\gamma }  F \left(\gamma , \gamma ;1; \frac{(\phi (t)  - \phi
(b))^2  } {(\phi (t)  + \phi (b))^2  }  \right)   w_{r} (x,0;b )\, db
\end{eqnarray*}
and
\begin{eqnarray*}
&  &
4 ^{-\gamma }  \phi'{}' (t)     \phi (t)^{-\gamma }  \int_{ 0}^{t}
   \phi (b) ^{-\gamma }  w(x, \phi (t)- \phi (b) ;b ) \,db  \\ 
&  &
 -  4 ^{-\gamma }  \gamma   (\phi' (t) )^{2} \phi (t)^{-\gamma-1 } \int_{ 0}^{t}
   \phi (b) ^{-\gamma }  w(x,  \phi (t)- \phi (b) ;b ) \,db \\
&  &
    -  2  \gamma    ( \phi'  (t))^2  \int_{ 0}^{t} 
  (\phi (t)  + \phi (b))\left( 4\phi (t)  \phi (b)  \right)^{-\gamma -1}  
  w(x,\phi (t)- \phi (b);b )\, db \\
  &  &
 +    \phi' (t) \int_{ 0}^{t} db
   \left( 4\phi (t)  \phi (b)  \right)^{-\gamma } \left(  \phi' (t)    
\frac{  \phi  (t)   +      \phi (b)  } {2 \phi (t)\phi (b)}  \right)   \gamma^2  w(x,(\phi (t)- \phi (b));b ) \\
  &  &
+ t^{\ell}\Bigg[ 2\gamma  \int_{ 0}^{t}
  (\phi (t)- \phi (b))  \left(  4\phi (t)   \phi (b)  \right)^{-\gamma-1 }    w  (x,(\phi (t)- \phi (b)) ;b )\,  db\\
  &  &
  - 8 \gamma^2  \phi (t)\int_{ 0}^{t} \, \phi (b) 
   \left( 4\phi (t)   \phi (b) \right)^{-\gamma -2}   (\phi (t)  - \phi (b))      w  (x,\phi (t)- \phi (b);b )\, db \Bigg]\,.
\end{eqnarray*}
In fact,  it is easily seen that the sum of the second group   vanishes:
\begin{eqnarray*} 
&  &
4 ^{-\gamma }     \phi'{}' (t)     \phi (t)^{-\gamma }  
   \phi (b) ^{-\gamma }     
- 4 ^{-\gamma }   \gamma  (\phi' (t) )^{2} \phi (t)^{-\gamma-1 }  
   \phi (b) ^{-\gamma }      \\
&  &
 - 2  \gamma    ( \phi'  (t))^2   
  (\phi (t)  + \phi (b))\left( 4\phi (t)  \phi (b)  \right)^{-\gamma -1} 
 +    \phi' (t) 
   \left( 4\phi (t)  \phi (b)  \right)^{-\gamma } \left(  \phi' (t)    
\frac{  \phi  (t)  -     \phi (b)  } {2 \phi (t)\phi (b)}  \right)   \gamma^2    \\
  &  &
+ t^{\ell}\Bigg[ 2\gamma  
  (\phi (t)- \phi (b))  \left(  4\phi (t)   \phi (b) \right)^{-\gamma-1 }   
  - 8  \gamma^2  \phi (t) \, \phi (b) 
   \left( 4\phi (t)   \phi (b) \right)^{-\gamma -2}   (\phi (t)  - \phi (b)) \Bigg] \\
& = & 
0 \,.
\end{eqnarray*}
Thus, it is proved that the function $u $ solves equation. Ii is easy to check that    $u$ takes the vanishing initial data.

\medskip

\noindent
Now we consider the case of $\ell \in (-2,4) $. We have only to prove that $ u_t$ is continuous up to $t=0$ and then check the second initial value, that is $u_t(x,0)$, of the solution $u $.
We differentiate integral depending on the parameter $t$ and obtain 
\begin{eqnarray*}
\partial_t u(x,t) 
& = &
c_\ell   \int_{ 0}^{t}
 \phi' (t)   \left(  4\phi (t)   \phi (b)     \right)^{-\gamma }  w(x,  \phi (t)- \phi (b) ;b ) \,db\\
&   &
+ c_\ell  \int_{ 0}^{t} db
  \int_{ 0}^{ \phi (t)- \phi (b)} 2\phi' (t)  (-\gamma ) (\phi (t)  + \phi (b))\left(  (\phi (t)  + \phi (b))^2  -r^2 \right)^{-\gamma -1} \\
&  & 
\hspace{2.5cm} \times F \left(\gamma , \gamma ;1; \beta  (t,b,r) \right)  w(x,r;b )
dr \\
&  &
+ c_\ell \int_{ 0}^{t} db
  \int_{ 0}^{ \phi (t)- \phi (b)} \left(  (\phi (t)  + \phi (b))^2  -r^2 \right)^{-\gamma } \left( \partial_t \beta  (t,b,r)\right)  F_z \left(\gamma , \gamma ;1; \beta  (t,b,r) \right)   w(x,r;b )\,
dr   \,.
\end{eqnarray*} 
Next we write it as follows:
\begin{eqnarray*}
\partial_t u(x,t) 
& = &
c_1  \phi' (t)    \phi (t) ^{-\gamma }  \int_{ 0}^{t}
   \phi (b)   ^{-\gamma }  w(x,  \phi (t)- \phi (b) ;b ) \,db\\
&   &
+ c_2  \phi' (t) \int_{ 0}^{t} db
  \int_{ 0}^{ \phi (t)- \phi (b)}   (\phi (t)  + \phi (b))\left(  (\phi (t)  + \phi (b))^2  -r^2 \right)^{-\gamma -1} \\
&  & 
\hspace{2.5cm} \times F \left(\gamma , \gamma ;1; \beta  (t,b,r) \right)  w(x,r;b )
dr \\
&  &
+ c_3  \phi' (t)  \int_{ 0}^{t} db
  \int_{ 0}^{ \phi (t)- \phi (b)} \left(  (\phi (t)  + \phi (b))^2  -r^2 \right)^{-\gamma } \\
&  &
\times \left( \left(  (\phi (t)  + \phi (b))^2  -r^2 \right)^{-2} \phi (b) ( \phi^2 (t)  - \phi^2 (b)    + 
   r^2   ) \right)  F_z \left(\gamma , \gamma ;1; \beta  (t,b,r) \right)   w(x,r;b )\,
dr \\
& = &
c_1  \phi' (t)    \phi (t) ^{-\gamma }  \int_{ 0}^{t}
   \phi (b)   ^{-\gamma }  w(x,  \phi (t)- \phi (b) ;b ) \,db\\
&   &
+ c_2  \phi' (t) \int_{ 0}^{t} db
  \int_{ 0}^{ \phi (t)- \phi (b)}   (\phi (t)  + \phi (b))\left(  (\phi (t)  + \phi (b))^2  -r^2 \right)^{-\gamma -1} \\
&  & 
\hspace{2.5cm} \times F \left(\gamma , \gamma ;1; \beta  (t,b,r) \right)  w(x,r;b )
dr \\
&  &
+ c_3  \phi' (t)  \int_{ 0}^{t} db
  \int_{ 0}^{ \phi (t)- \phi (b)} \left(  (\phi (t)  + \phi (b))^2  -r^2 \right)^{-2-\gamma } \\
&  &
\hspace{2.5cm} \times  \phi (b) \left( ( \phi^2 (t)  - \phi^2 (b)    + 
   r^2   ) \right)  F_z \left(\gamma , \gamma ;1; \beta  (t,b,r) \right)   w(x,r;b )\,dr \,.
\end{eqnarray*} 
If we denote $k=\ell/2 $, then
\begin{eqnarray*}
\partial_t u(x,t) 
& = &
c_1  t^{\frac{k}{2} }  \int_{ 0}^{t}
   b^{-\frac{k}{2}}  w(x,  \phi (t)- \phi (b) ;b ) \,db \quad \mbox{\rm (the first term)}\\
&   &
+ c_2  t^k \int_{ 0}^{t} db \,(\phi (t)  + \phi (b))
  \int_{ 0}^{ \phi (t)- \phi (b)}  \left(  (\phi (t)  + \phi (b))^2  -r^2 \right)^{-\gamma -1} \\
&  & 
\hspace{2.5cm} \times F \left(\gamma , \gamma ;1; \beta  (t,b,r) \right)  w(x,r;b )
dr \quad \mbox{\rm (the second term)}\\
&  &
+ c_3   t^k  \int_{ 0}^{t} db \,b^{k+1}
  \int_{ 0}^{ \phi (t)- \phi (b)} \left(  (\phi (t)  + \phi (b))^2  -r^2 \right)^{-2-\gamma } \\
&  &
\hspace{2.5cm} \times\left(   \phi^2 (t)  - \phi^2 (b)    + 
   r^2  \right)  F_z \left(\gamma , \gamma ;1; \beta  (t,b,r) \right)   w(x,r;b )\,dr \quad \mbox{\rm (the third term)}\,.
\end{eqnarray*}
We have to prove the convergence of all last integrals uniformly  on  compact (in $x$) sets. The first term is convergent if $\ell <4 $ and
\begin{eqnarray*}
&   &
 \left| t^{\frac{k}{2} }  \int_{ 0}^{t}
   b^{-\frac{k}{2}}  w(x,  \phi (t)- \phi (b) ;b ) \,db \right| \lesssim t\,.
\end{eqnarray*}
Consider the second term. According to Lemma~\ref{L5.1}, it can be estimated as follows: 
\begin{eqnarray*}
&   &
\Big|  t^k  \int_{ 0}^{t} db\,(\phi (t)  + \phi (b))
  \int_{ 0}^{ \phi (t)- \phi (b)}   \left(  (\phi (t)  + \phi (b))^2  -r^2 \right)^{-\gamma -1} \\
&  & 
\hspace{2.5cm} \times F \left(\gamma , \gamma ;1; \beta  (t,b,r) \right)  w(x,r;b )dr \Big|\\
&  \lesssim  &
 t^k  \int_{ 0}^{t} db \,(\phi (t)  + \phi (b))
  \int_{ 0}^{ \phi (t)- \phi (b)}  \left(  (\phi (t)  + \phi (b))^2  -r^2 \right)^{-\gamma -1} \\
&  & 
\hspace{1cm} \times \left[ 1+ \left(  (\phi (t)  + \phi (b))^2  -r^2 \right)^{ -1 }\phi (t)  \phi (b) 
+\left(  (\phi (t)  + \phi (b))^2  -r^2 \right)^{2\gamma -1  } \phi (t)^{1-2\gamma }  \phi (b) ^{1-2\gamma }  \right]  
dr \\
&  \lesssim   &
 t^k  \int_{ 0}^{t} db (\phi (t)  + \phi (b))
  \int_{ 0}^{ \phi (t)- \phi (b)}  \left(  (\phi (t)  + \phi (b))^2  -r^2 \right)^{-\gamma -1}   dr \\
&  &
+
 t^k \phi (t)  \int_{ 0}^{t} db \, \phi (b)  (\phi (t)  + \phi (b))
  \int_{ 0}^{ \phi (t)- \phi (b)}  \left(  (\phi (t)  + \phi (b))^2  -r^2 \right)^{-\gamma -2}   dr \\
&  &
+
 t^k \phi (t)^{1-2\gamma }  \int_{ 0}^{t} db \, \phi (b) ^{1-2\gamma }  (\phi (t)  + \phi (b))
  \int_{ 0}^{ \phi (t)- \phi (b)}  \left(  (\phi (t)  + \phi (b))^2  -r^2 \right)^{\gamma -2}    dr \,.
\end{eqnarray*}
Next we apply Lemma~\ref{L5.2} and estimate the last terms as follows 
\begin{eqnarray*} 
& \lesssim   &
 t^k  \int_{ 0}^{t} db (\phi (t)  + \phi (b))
(  \phi (t)- \phi (b)) (  \phi (t)+ \phi (b))^{-2\gamma -2} F\left(\frac{1}{2}, \gamma +1; \frac{3}{2}; \frac{\left(\phi (t)- \phi (b)\right)^2}{\left(\phi (t)  + \phi (b)\right)^2}\right) \\
&  &
+
 t^k \phi (t)  \int_{ 0}^{t} db \, \phi (b)  (\phi (t)  + \phi (b))
 \left(\phi (t)  - \phi (b) \right)  \left(\phi (t)  + \phi (b)\right)^{-2    \gamma -4 } F \left(\frac{1}{2},2+  \gamma ;\frac{3}{2};\frac{\left(\phi (t)- \phi (b)\right)^2}{\left(\phi (t)  + \phi (b)\right)^2}\right) \\
&  &
+
 t^k \phi (t)^{1-2\gamma }  \int_{ 0}^{t} db \, \phi (b) ^{1-2\gamma }  (\phi (t)  + \phi (b))
 \left(\phi (t)  - \phi (b) \right)  \left(\phi (t)  + \phi (b)\right)^{ 2    \gamma -4 } F \left(\frac{1}{2}, 2-\gamma ;\frac{3}{2};\frac{\left(\phi (t)- \phi (b)\right)^2}{\left(\phi (t)  + \phi (b)\right)^2}\right) \\ 
& \lesssim   &
 t^k  \int_{ 0}^{t} db  
(  \phi (t)- \phi (b)) (  \phi (t)+ \phi (b))^{-2\gamma -1} F\left(\frac{1}{2}, \gamma +1; \frac{3}{2}; \frac{\left(\phi (t)- \phi (b)\right)^2}{\left(\phi (t)  + \phi (b)\right)^2}\right) \\
&  &
+
 t^k \phi (t)  \int_{ 0}^{t} db \, \phi (b)  
 \left(\phi (t)  - \phi (b) \right)  \left(\phi (t)  + \phi (b)\right)^{-2    \gamma -3 } F \left(\frac{1}{2},2+  \gamma ;\frac{3}{2};\frac{\left(\phi (t)- \phi (b)\right)^2}{\left(\phi (t)  + \phi (b)\right)^2}\right) \\
&  &
+
 t^k \phi (t)^{1-2\gamma }  \int_{ 0}^{t} db \, \phi (b) ^{1-2\gamma }  
 \left(\phi (t)  - \phi (b) \right)  \left(\phi (t)  + \phi (b)\right)^{ 2    \gamma -3 } F \left(\frac{1}{2}, 2-\gamma ;\frac{3}{2};\frac{\left(\phi (t)- \phi (b)\right)^2}{\left(\phi (t)  + \phi (b)\right)^2}\right)  \\
  & = :&
  A+B+C\,.
\end{eqnarray*}
Then we apply Lemma~\ref{L5.5} to estimate $A$ as follows: 
\begin{eqnarray*} 
A
& =  &
 t^k  \int_{ 0}^{t} db  
(  \phi (t)- \phi (b)) (  \phi (t)+ \phi (b))^{-2\gamma -1} F\left(\frac{1}{2}, \gamma +1; \frac{3}{2}; \frac{\left(\phi (t)- \phi (b)\right)^2}{\left(\phi (t)  + \phi (b)\right)^2}\right) \\
& \lesssim  &
 t^k  \int_{ 0}^{t} db  
(  \phi (t)- \phi (b)) (  \phi (t)+ \phi (b))^{-2\gamma -1}\\
&  &
\times \Bigg[ 1+ \left( \frac{ 4 \phi (t)  \phi (b) }{\left(\phi (t)  + \phi (b) \right)^2 } \right) + \left( \frac{ 4 \phi (t)  \phi (b) }{\left(\phi (t)  + \phi (b) \right)^2 } \right)^{-\gamma  } \left[1+ \left( \frac{ 4 \phi (t)  \phi (b) }{\left(\phi (t)  + \phi (b) \right)^2 } \right)\right]  \Bigg]   \\
& \lesssim  &
 t 
 + \, t^k  \int_{ 0}^{t}   
(  \phi (t)- \phi (b)) (  \phi (t)+ \phi (b))^{-2\gamma -3}   \left(  \phi (t)  \phi (b)   \right)      db \\
&   &
+ \,  t^k  \int_{ 0}^{t}    
(  \phi (t)- \phi (b)) (  \phi (t)+ \phi (b))^{ -1}  \left(  \phi (t)  \phi (b)   \right)^{-\gamma  }   db\\
&    &
+ \, t^k 
\int_{ 0}^{t}    
(  \phi (t)- \phi (b)) (  \phi (t)+ \phi (b))^{-3}  \left(  \phi (t)  \phi (b)   \right)^{1-\gamma  }      db \\
& \lesssim  &
 t\,.
\end{eqnarray*}
Further, we use condition $\ell < 4$ and Lemma~\ref{L5.9} to estimate $B$ as follows:
\begin{eqnarray*} 
B
&  \lesssim & 
 t^k \phi (t)  \int_{ 0}^{t} db \, \phi (b)  
 \left(\phi (t)  - \phi (b) \right)  \left(\phi (t)  + \phi (b)\right)^{-2    \gamma -3 } F \left(\frac{1}{2},2+  \gamma ;\frac{3}{2};\frac{\left(\phi (t)- \phi (b)\right)^2}{\left(\phi (t)  + \phi (b)\right)^2}\right)  \\
 & \lesssim &
 t^k \phi (t)  \int_{ 0}^{t} db \, \phi (b)  
 \left(\phi (t)  - \phi (b) \right)  \left(\phi (t)  + \phi (b)\right)^{-2    \gamma -3 } \\
&  &
\times \Bigg[  \left( \frac{ 4 \phi (t)  \phi (b) }{\left(\phi (t)  + \phi (b) \right)^2 } \right)^{-1}+ \left( \frac{ 4 \phi (t)  \phi (b) }{\left(\phi (t)  + \phi (b) \right)^2 } \right)^{-\gamma -1} \left[1+ \left( \frac{ 4 \phi (t)  \phi (b) }{\left(\phi (t)  + \phi (b) \right)^2 } \right)\right]  \Bigg]   \\ 
 & \lesssim &
 t^k   \int_{ 0}^{t}    
 \left(\phi (t)  - \phi (b) \right)  \left(\phi (t)  + \phi (b)\right)^{-2    \gamma -1 }  \, db \\ 
 &   &
+ t^k  \left(   \phi (t)   \right)^{-\gamma  }\int_{ 0}^{t}   
 \left(\phi (t)  - \phi (b) \right)  \left(\phi (t)  + \phi (b)\right)^{   -1 }   \left(    \phi (b)   \right)^{-\gamma } \, db   \\ 
 &   &
+  t^k \phi (t)^{1-\gamma }   \int_{ 0}^{t}    \phi (b)^{1-\gamma }  
 \left(\phi (t)  - \phi (b) \right)  \left(\phi (t)  + \phi (b)\right)^{  -3 }    \, db  \\   
 & \lesssim &
t \,.
\end{eqnarray*}
Then we apply Lemma~\ref{L5.9} to estimate $C$ 
\begin{eqnarray*}  
C
& = &
 t^k \phi (t)^{1-2\gamma }  \int_{ 0}^{t} db \, \phi (b) ^{1-2\gamma }  
 \left(\phi (t)  - \phi (b) \right)  \left(\phi (t)  + \phi (b)\right)^{ 2    \gamma -3 } F \left(\frac{1}{2}, 2-\gamma ;\frac{3}{2};\frac{\left(\phi (t)- \phi (b)\right)^2}{\left(\phi (t)  + \phi (b)\right)^2}\right)   \\
&  \lesssim  &
 t 
 +  t^k \phi (t)^{1-2\gamma }  \int_{ 0}^{t} db \, \phi (b) ^{1-2\gamma }  
 \left(\phi (t)  - \phi (b) \right)  \left(\phi (t)  + \phi (b)\right)^{ -1}  \left( \phi (t)  \phi (b) \right)^{\gamma -1} \\
 &  &
 + t^k \phi (t)^{1-2\gamma }  \int_{ 0}^{t} db \, \phi (b) ^{1-2\gamma }  
 \left(\phi (t)  - \phi (b) \right)  \left(\phi (t)  + \phi (b)\right)^{ 2    \gamma -3 } \left( \frac{ 4 \phi (t)  \phi (b) }{\left(\phi (t)  + \phi (b) \right)^2 } \right)^{\gamma }   \\
&  \lesssim  &
 t  \,.
\end{eqnarray*} 
Hence, if $  \ell \in (-2,4)$, then the second term is estimated as follows
\begin{eqnarray*}
&   &
 t^k  \int_{ 0}^{t} db
  \int_{ 0}^{ \phi (t)- \phi (b)}   (\phi (t)  + \phi (b))\left(  (\phi (t)  + \phi (b))^2  -r^2 \right)^{-\gamma -1}  
    \lesssim   
t\,.
\end{eqnarray*} 
Then we consider the third term and apply  
 Lemma~\ref{L5.7}.
Theorem~\ref{T1.8} is proved. \hfill $\square$

Theorem~\ref{T1.8b} follows directly from  Theorem~\ref{T1.8}. 

\section{Operators $K_0$ , $K_1$. Problem without source term. Proof of Theorem~\ref{TK0K1}}
\label{S4}

First we check the initial values of the images of smooth functions.
In the following lemma we consider operator $K_0 $, which is defined for all $ \ell \in (-\infty,-2)\cup (0,\infty)$.  
\begin{lemma}
Assume that $\gamma >0 $ and $\phi (0)=\phi' (0)=0 $. Let $\Omega  $ be a backward time connected domain. Then for every   function $v=v(x,t) \in C_{x,t}^{0,1} (\overline{\Omega }) $ 
we have 
\begin{eqnarray*}
(K_0 v)   (x,0)
  =     v   (x, 0)\,,\quad 
(K_0 v)_ t (x,0)
  =      0  \quad \mbox{for all} \quad x \in \widetilde {\Omega }\,.
\end{eqnarray*} 
\end{lemma}
\medskip

\noindent
{\bf Proof.} It is evident that
\begin{eqnarray*}
(K_0 v)   (x,0)
& = &  
2^{2-2\gamma }
\frac{\Gamma \left( 2\gamma  \right) } {\Gamma^2 \left( \gamma  \right) }
\int_{0}^1   v  (x, 0)
(1-s^2)^{\gamma - 1   } ds =v   (x, 0)\,, 
\end{eqnarray*}
and
\begin{eqnarray}
\label{13}
\partial_t K_0 v(x,t)
& = &
\phi' (t)  2^{2-2\gamma }
\frac{\Gamma \left( 2\gamma  \right) } {\Gamma^2 \left( \gamma  \right) }
\int_{0}^1   \partial_r v  (x, r)|_{r=\phi (t) s}
s(1-s^2)^{\gamma - 1   } ds    \,.
\end{eqnarray}
Lemma is proved. \hfill $\square$

Remind that the operator $K_1 $  is defined for all $ \ell \in (-\infty,-4)\cup (-2,\infty)$. 
\begin{lemma}
Assume that $\gamma <1 $ and $\phi (0) =0 $,   $\lim_{t \to 0}t \phi' (t) =0 $. Let $\Omega  $ be a backward time connected domain. For every  
 function $v=v(x,t)  \in C_{x,t}^{0,1} (\overline{\Omega })  $ 
we have 
\begin{eqnarray*}
  (K_1 v) (x,0)
  =  0, \quad (  K_1 v)_t(x,0)
  =   v   (x, 0)\,.
\end{eqnarray*} 
\end{lemma}
\medskip

\noindent
{\bf Proof.} The first  relation is evident. 
Next, we have 
\begin{eqnarray*}
\label{dtK1}
\partial_t K_1 v(x,t)
& = &
2^{2\gamma }
\frac{\Gamma \left(2- 2\gamma  \right) } {\Gamma^2 \left( 1- \gamma  \right) }
\int_{0}^1   v  (x, \phi (t) s)
(1-s^2)^{- \gamma  } ds  \nonumber \\
&  &
+\, t \phi' (t) 2^{2\gamma }
\frac{\Gamma \left(2- 2\gamma  \right) } {\Gamma^2 \left( 1- \gamma  \right) }
\int_{0}^1   \partial_r v  (x, r)|_{r=\phi (t) s}
s(1-s^2)^{- \gamma  } ds  \,.
\end{eqnarray*}
Hence
\begin{eqnarray*}
(  K_1 v)_t(x,0)
& = &
2^{2\gamma }
\frac{\Gamma \left(2- 2\gamma  \right) } {\Gamma^2 \left( 1- \gamma  \right) }
\int_{0}^1   v  (x, 0)
(1-s^2)^{- \gamma  } ds  = v   (x, 0)
\end{eqnarray*} 
since $\phi (0) =0 $ and $\lim_{t \to 0}t \phi' (t) =0 $. Lemma is proved.
\hfill $\square$
\medskip

Now we turn to the equation.
\begin{proposition}
Assume that $\gamma >0$ and the domain $\Omega $ is backward time connected. For $ v \in C_{x,t}^{m,2}(\overline{\Omega_\phi} )$, such that
\begin{eqnarray}
\left( \partial_{t}^2-A(x,\partial_x)\right) v
& = &
0\quad \mbox{at}\,\, x=x_0\quad \mbox{for all }\,\, t \in (0, \phi (T)]\,,
\end{eqnarray}
where $(x_0,T) \in \Omega  $, we have
 \begin{eqnarray}
\left( \partial_{t}^2-t^{\ell}A(x,\partial_x)\right) K_0v
& = &
2^{1-2\gamma }
\frac{ \ell \, \Gamma \left( 2\gamma  \right) } { 2\gamma\Gamma^2 \left( \gamma  \right) }
  t^{\frac{\ell}{2}-1}    \partial_t v   (x, 0) \quad \mbox{at}\,\, x=x_0\quad \mbox{for all }\,\, t \in (0,T ]\,.
\end{eqnarray}
\end{proposition}

\noindent
{\bf Proof. }
Consider  (\ref{13}), it follows:
\begin{eqnarray*}
\partial_t^2 K_0 v(x,t)
& = &
  \phi'{}' (t)  2^{2-2\gamma }
\frac{\Gamma \left( 2\gamma  \right) } {\Gamma^2 \left( \gamma  \right) }
\int_{0}^1   \partial_r v (x, r)|_{r=\phi (t) s}
s(1-s^2)^{\gamma - 1   } ds  \\
&  &
+ (  \phi' (t))^2  2^{2-2\gamma }
\frac{\Gamma \left( 2\gamma  \right) } {\Gamma^2 \left( \gamma  \right) }
\int_{0}^1   \partial_r^2 v (x, r)|_{r=\phi (t) s}
s^2(1-s^2)^{\gamma - 1   } ds  \,.
\end{eqnarray*}
On the other hand, 
\begin{eqnarray*} 
A(x,\partial_x)K_0 v(x,t)
& = &
   2^{2-2\gamma }
\frac{\Gamma \left( 2\gamma  \right) } {\Gamma^2 \left( \gamma  \right) }
\int_{0}^1     A(x,\partial_x) v (x, \phi (t) s) 
(1-s^2)^{\gamma - 1   } ds \\
& = &
   2^{2-2\gamma }
\frac{\Gamma \left( 2\gamma  \right) } {\Gamma^2 \left( \gamma  \right) }
\int_{0}^1    \partial_r^2 v (x, r)|_{r=\phi (t) s} 
(1-s^2)^{\gamma - 1   } ds \,.
\end{eqnarray*}
Consequently,
\begin{eqnarray*} 
(\partial_t^2  -t^\ell A(x,\partial_x) ) K_0 v(x,t)
& = &
  \phi'{}' (t)  2^{2-2\gamma }
\frac{\Gamma \left( 2\gamma  \right) } {\Gamma^2 \left( \gamma  \right) }
\int_{0}^1   \partial_r v (x, r)|_{r=\phi (t) s}
s(1-s^2)^{\gamma - 1   } ds  \\
&  &
+ (  \phi' (t))^2  2^{2-2\gamma }
\frac{\Gamma \left( 2\gamma  \right) } {\Gamma^2 \left( \gamma  \right) }
\int_{0}^1   \partial_r^2 v (x, r)|_{r=\phi (t) s}
s^2(1-s^2)^{\gamma - 1   } ds \\
&  &
-  t^{\ell} 2^{2-2\gamma }
\frac{\Gamma \left( 2\gamma  \right) } {\Gamma^2 \left( \gamma  \right) }
\int_{0}^1   \left[ \partial_r^2 v (x, r)|_{r=\phi (t) s}\right]
(1-s^2)^{\gamma - 1   } ds \,.
\end{eqnarray*}
Integration by parts in the first term leads to 
\begin{eqnarray*}
(\partial_t^2  -t^\ell A(x,\partial_x) ) K_0 v(x,t)
& = &
 - \phi'{}' (t)  2^{2-2\gamma }
\frac{\Gamma \left( 2\gamma  \right) } {\Gamma^2 \left( \gamma  \right) }
\int_{0}^1   \partial_r v (x, r)|_{r=\phi (t) s}
\frac{1}{2\gamma }\frac{d}{ ds}(1-s^2)^{\gamma } ds  \\
&  &
+ (  \phi' (t))^2  2^{2-2\gamma }
\frac{\Gamma \left( 2\gamma  \right) } {\Gamma^2 \left( \gamma  \right) }
\int_{0}^1   \partial_r^2 v (x, r)|_{r=\phi (t) s}
s^2(1-s^2)^{\gamma - 1   } ds \\
&  &
-   t^{\ell}  2^{2-2\gamma }
\frac{\Gamma \left( 2\gamma  \right) } {\Gamma^2 \left( \gamma  \right) }
\int_{0}^1   \left[ \partial_r^2 v_{\varphi _0 } (x, r)|_{r=\phi (t) s}\right]
(1-s^2)^{\gamma - 1   } ds \\
& = &
  \phi'{}' (t)    2^{2-2\gamma }
\frac{\Gamma \left( 2\gamma  \right) } {\Gamma^2 \left( \gamma  \right) }
   \partial_r v (x, 0) 
\frac{1}{2\gamma }   \\
&  &
+\phi'{}' (t)   2^{2-2\gamma }
\frac{\Gamma \left( 2\gamma  \right) } {\Gamma^2 \left( \gamma  \right) }\frac{1}{2\gamma }
\int_{0}^1  \phi (t) \left[ \partial_r^2 v (x, r)|_{r=\phi (t) s}\right]
(1-s^2)^{\gamma } ds  \\
&  &
+  t^{\ell}  2^{2-2\gamma }
\frac{\Gamma \left( 2\gamma  \right) } {\Gamma^2 \left( \gamma  \right) }
\int_{0}^1   \partial_r^2 v (x, r)|_{r=\phi (t) s}
s^2(1-s^2)^{\gamma - 1   } ds \\
&  &
-   t^{\ell}  2^{2-2\gamma }
\frac{\Gamma \left( 2\gamma  \right) } {\Gamma^2 \left( \gamma  \right) }
\int_{0}^1   \left[ \partial_r^2 v (x, r)|_{r=\phi (t) s}\right]
(1-s^2)^{\gamma - 1   } ds \,.
\end{eqnarray*}
Then,  (\ref{phi}) implies 
\begin{eqnarray*}
(\partial_t^2  -t^\ell A(x,\partial_x) ) K_0 v(x,t)
& = &
 \phi'{}' (t)   2^{2-2\gamma }
\frac{\Gamma \left( 2\gamma  \right) } {\Gamma^2 \left( \gamma  \right) }
   \partial_r v (x, 0) 
\frac{1}{2\gamma }   \\
&  &
+ \frac{1}{2\gamma }\phi'{}' (t)   \phi (t)  2^{2-2\gamma }
\frac{\Gamma \left( 2\gamma  \right) } {\Gamma^2 \left( \gamma  \right) }
\int_{0}^1 \left[ \partial_r^2 v (x, r)|_{r=\phi (t) s}\right]
(1-s^2)^{\gamma } ds  \\
&  &
-   ( \phi' (t))^2   2^{2-2\gamma }
\frac{\Gamma \left( 2\gamma  \right) } {\Gamma^2 \left( \gamma  \right) }
\int_{0}^1   \left[ \partial_r^2 v (x, r)|_{r=\phi (t) s}\right]
(1-s^2)^{\gamma  } ds \\
& = &
 \phi'{}' (t)   2^{2-2\gamma }
\frac{\Gamma \left( 2\gamma  \right) } {\Gamma^2 \left( \gamma  \right) }
   \partial_r v (x, 0) 
\frac{1}{2\gamma }  \,.
\end{eqnarray*} 
Proposition is proved. \hfill $\square$

\begin{proposition}
Assume that $\gamma <1$ and $\Omega $ is backward time connected. For $ v \in C_{x,t}^{m,2}(\overline{\Omega_\phi} )$ such that
\begin{eqnarray}
\left( \partial_{t}^2-A(x,\partial_x)\right) v
& = &
0\quad \mbox{at}\,\,\, x= x_0\quad \mbox{for all }\,\, t \in [0, \phi (T)]
\end{eqnarray}
we have
 \begin{eqnarray}
\left( \partial_{t}^2-t^{\ell}A(x,\partial_x)\right) K_1 v
& = &
  ( \ell+2 ) 2^{2\gamma-1}
\frac{\Gamma \left(2- 2\gamma  \right) } {\Gamma^2 \left( 1- \gamma  \right) }
t^{\frac{\ell}{2}}    \partial_t v   (x, 0)\quad \mbox{at}\,\,\, x=x_0\quad \mbox{for all }\,\, t \in [0, T]\,.
\end{eqnarray}
\end{proposition}
\medskip

\noindent
{\bf Proof. }
Consider the derivative $\partial_t^2 K_1 v(x,t)$. According to (\ref{dtK1}) we have:
\begin{eqnarray*} 
\partial_t^2 K_1 v(x,t) 
& =  &
\partial_t \Bigg[ 2^{2\gamma }
\frac{\Gamma \left(2- 2\gamma  \right) } {\Gamma^2 \left( 1- \gamma  \right) }
\int_{0}^1   v  (x, \phi (t) s)
(1-s^2)^{- \gamma  } ds  \nonumber \\
&  & 
+\,  t \phi' (t) 2^{2\gamma }
\frac{\Gamma \left(2- 2\gamma  \right) } {\Gamma^2 \left( 1- \gamma  \right) }
\int_{0}^1   \partial_r v  (x, r)|_{r=\phi (t) s}
s(1-s^2)^{- \gamma  } ds  \Bigg]\,.
\end{eqnarray*}
Then
\begin{eqnarray*} 
\partial_t^2 K_1 v(x,t)  
& = &
\phi' (t) 2^{2\gamma +1}
\frac{\Gamma \left(2- 2\gamma  \right) } {\Gamma^2 \left( 1- \gamma  \right) }
\int_{0}^1   \partial_r v(x, r)|_{r=\phi (t) s}
s (1-s^2)^{- \gamma  } ds  \nonumber \\
&  & 
+\,  t \phi'{}' (t) 2^{2\gamma }
\frac{\Gamma \left(2- 2\gamma  \right) } {\Gamma^2 \left( 1- \gamma  \right) }
\int_{0}^1   \partial_r v (x, r)|_{r=\phi (t) s}
s(1-s^2)^{- \gamma  } ds  \\
&  & 
+\,  t (\phi' (t))^2 2^{2\gamma }
\frac{\Gamma \left(2- 2\gamma  \right) } {\Gamma^2 \left( 1- \gamma  \right) }
\int_{0}^1   \partial_r^2 v (x, r)|_{r=\phi (t) s}
s^2(1-s^2)^{- \gamma  } ds  \,,
\end{eqnarray*}
that is
\begin{eqnarray*} 
\partial_t^2 K_1 v(x,t)  
& = &
-\phi' (t) 2^{2\gamma}
\frac{\Gamma \left(2- 2\gamma  \right) } {\Gamma^2 \left( 1- \gamma  \right) }
\int_{0}^1   \partial_r v (x, r)|_{r=\phi (t) s}
\frac{1}{1-\gamma }\frac{d}{ds}(1-s^2)^{1- \gamma  } ds  \nonumber \\
&  & 
-\,  t \phi'{}' (t) 2^{2\gamma }
\frac{\Gamma \left(2- 2\gamma  \right) } {\Gamma^2 \left( 1- \gamma  \right) }
\int_{0}^1   \partial_r v(x, r)|_{r=\phi (t) s}
\frac{1}{2(1-\gamma )}\frac{d}{ds}(1-s^2)^{1- \gamma  }ds  \\
&  & 
+\,  t (\phi' (t))^2 2^{2\gamma }
\frac{\Gamma \left(2- 2\gamma  \right) } {\Gamma^2 \left( 1- \gamma  \right) }
\int_{0}^1   \partial_r^2 v (x, r)|_{r=\phi (t) s}
s^2(1-s^2)^{- \gamma  } ds   \,.
\end{eqnarray*}
Next we integrate by parts
\begin{eqnarray*} 
\partial_t^2 K_1 v(x,t) 
& =  &
\frac{1}{1-\gamma } 2^{2\gamma}
\frac{\Gamma \left(2- 2\gamma  \right) } {\Gamma^2 \left( 1- \gamma  \right) }
 \phi' (t)   \partial_r v  (x, 0) 
   \\
&  &
+\phi' (t)\phi (t) \frac{1}{1-\gamma }2^{2\gamma}
\frac{\Gamma \left(2- 2\gamma  \right) } {\Gamma^2 \left( 1- \gamma  \right) }
\int_{0}^1   \partial_r^2 v  (x, r)|_{r=\phi (t) s}
(1-s^2)^{1- \gamma  } ds  \nonumber \\
&  & 
+\, \frac{1}{2(1-\gamma )} t \phi'{}' (t) 2^{2\gamma }
\frac{\Gamma \left(2- 2\gamma  \right) } {\Gamma^2 \left( 1- \gamma  \right) }
 \partial_r v (x, 0)  \\
&  & 
+\,  t \phi'{}' (t) \phi (t)\frac{1}{2(1-\gamma )}2^{2\gamma }
\frac{\Gamma \left(2- 2\gamma  \right) } {\Gamma^2 \left( 1- \gamma  \right) }
\int_{0}^1   \partial_r^2 v (x, r)|_{r=\phi (t) s}
(1-s^2)^{1- \gamma  }ds  \\
&  & 
+\,  t (\phi' (t))^2 2^{2\gamma }
\frac{\Gamma \left(2- 2\gamma  \right) } {\Gamma^2 \left( 1- \gamma  \right) }
\int_{0}^1   \partial_r^2 v (x, r)|_{r=\phi (t) s}
s^2(1-s^2)^{- \gamma  } ds \,.
\end{eqnarray*}
It is easily seen that
\begin{eqnarray*}  
\phi' (t)\phi (t) \frac{1}{1-\gamma }+t \phi'{}' (t) \phi (t)\frac{1}{2(1-\gamma )}= t (\phi' (t))^2\,,
\end{eqnarray*}
and, consequently,
\begin{eqnarray*} 
\partial_t^2 K_1 v(x,t) 
& =  &
\Bigg[ \frac{1}{1-\gamma } 2^{2\gamma}
\frac{\Gamma \left(2- 2\gamma  \right) } {\Gamma^2 \left( 1- \gamma  \right) }
 \phi' (t)   
+\, \frac{1}{2(1-\gamma )} t \phi'{}' (t) 2^{2\gamma }
\frac{\Gamma \left(2- 2\gamma  \right) } {\Gamma^2 \left( 1- \gamma  \right) }
\Bigg] \partial_r v (x, 0)  \\
&  & 
+\,  t (\phi' (t))^2 2^{2\gamma }
\frac{\Gamma \left(2- 2\gamma  \right) } {\Gamma^2 \left( 1- \gamma  \right) }
\int_{0}^1   \partial_r^2 v (x, r)|_{r=\phi (t) s}
(1-s^2)^{1- \gamma  } ds \\
&  & 
+\,  t (\phi' (t))^2 2^{2\gamma }
\frac{\Gamma \left(2- 2\gamma  \right) } {\Gamma^2 \left( 1- \gamma  \right) }
\int_{0}^1   \partial_r^2 v (x, r)|_{r=\phi (t) s}
s^2(1-s^2)^{- \gamma  } ds \\
& =  &
 2^{2\gamma}\frac{\Gamma \left(2- 2\gamma  \right) } {\Gamma^2 \left( 1- \gamma  \right) }t (\phi' (t))^2\frac{1}{\phi (t)}\partial_r v (x, 0)  \\
&  & 
+\,  t (\phi' (t))^2 2^{2\gamma }
\frac{\Gamma \left(2- 2\gamma  \right) } {\Gamma^2 \left( 1- \gamma  \right) }
\int_{0}^1   \partial_r^2 v (x, r)|_{r=\phi (t) s}
(1-s^2)^{1- \gamma  } ds \\
&  & 
+\,  t (\phi' (t))^2 2^{2\gamma }
\frac{\Gamma \left(2- 2\gamma  \right) } {\Gamma^2 \left( 1- \gamma  \right) }
\int_{0}^1   \partial_r^2 v (x, r)|_{r=\phi (t) s}
s^2(1-s^2)^{- \gamma  } ds \,.
\end{eqnarray*}
On the other hand
\begin{eqnarray*}    
t^{\ell}  A(x,\partial_x)  K_1 v(x,t)
& =  &
 t (\phi' (t))^2 2^{2\gamma }
\frac{\Gamma \left(2- 2\gamma  \right) } {\Gamma^2 \left( 1- \gamma  \right) }
\int_{0}^1  A(x,\partial_x) v (x, \phi (t) s)
(1-s^2)^{- \gamma  } ds \\
& = &
 t (\phi' (t))^2 2^{2\gamma }
\frac{\Gamma \left(2- 2\gamma  \right) } {\Gamma^2 \left( 1- \gamma  \right) }
\int_{0}^1  \partial_r^2 v (x, r)|_{r=\phi (t) s}
(1-s^2)^{- \gamma  } ds\,.
\end{eqnarray*}
Hence
\begin{eqnarray*}    
\partial_t^2 K_1 v(x,t) - t^{\ell}  A(x,\partial_x)  K_1 v(x,t)
& = &
 2^{2\gamma}\frac{\Gamma \left(2- 2\gamma  \right) } {\Gamma^2 \left( 1- \gamma  \right) }t (\phi' (t))^2\frac{1}{\phi (t)}\partial_r v (x, 0)   \\
&  & 
+\,  t (\phi' (t))^2 2^{2\gamma }
\frac{\Gamma \left(2- 2\gamma  \right) } {\Gamma^2 \left( 1- \gamma  \right) }
\int_{0}^1   \partial_r^2 v (x, r)|_{r=\phi (t) s}
(1-s^2)^{1- \gamma  } ds \\
&  & 
+\,  t (\phi' (t))^2 2^{2\gamma }
\frac{\Gamma \left(2- 2\gamma  \right) } {\Gamma^2 \left( 1- \gamma  \right) }
\int_{0}^1   \partial_r^2 v (x, r)|_{r=\phi (t) s}
s^2(1-s^2)^{- \gamma  } ds \\
&  &
- t (\phi' (t))^2 2^{2\gamma }
\frac{\Gamma \left(2- 2\gamma  \right) } {\Gamma^2 \left( 1- \gamma  \right) }
\int_{0}^1  \partial_r^2 v (x, r)|_{r=\phi (t) s}
(1-s^2)^{- \gamma  } ds \\
& = & 
 2^{2\gamma}\frac{\Gamma \left(2- 2\gamma  \right) } {\Gamma^2 \left( 1- \gamma  \right) }t (\phi' (t))^2\frac{1}{\phi (t)}\partial_r v (x, 0)  \,.
\end{eqnarray*}
Proposition is proved. \hfill $\square$
\medskip

In particular, for $\ell \in   (0,\infty) $ both lemmas and propositions  are applicable  and we arrive at the following result. 
\begin{corollary}
Assume that $\ell \in   (0,\infty) $ and $\Omega  $ is backward time connected domain. For given   $\varphi =\varphi  (x) $ 
assume that the functions $v_{\varphi_0} \in C_{x,t}^{m,2}  (\overline{\Omega _\phi} ) $ and  $v_{\varphi_1}  \in C_{x,t}^{m,2}  (\overline{\Omega _\phi} )$,  solve    the  problem
\[
v_{tt}-  A(x,\partial_x)  v =0, \quad v(x,0)= \varphi_i  (x)\, , \quad v_t(x,0)=0 \quad \mbox{for all} \quad (x,t) \in \Omega_\phi   ,
\]
with $i=0$ and $i=1$, respectively. Then   the function 
\begin{eqnarray}
\label{1.7}
u(x,t)
& = &
 2^{2-2\gamma }
\frac{\Gamma \left( 2\gamma  \right) } {\Gamma^2 \left( \gamma  \right) }
\int_{0}^1   v_{\varphi _0 }  (x, \phi (t) s)
(1-s^2)^{\gamma - 1   } ds \\
&  &
+\, t 2^{2\gamma }
\frac{\Gamma \left(2- 2\gamma  \right) } {\Gamma^2 \left( 1- \gamma  \right) }
\int_{0}^1   v_{\varphi _1 }   (x, \phi (t) s)
(1-s^2)^{- \gamma  } ds , \qquad   (x,t) \in \Omega\,,  \nonumber
\end{eqnarray}
is a solution of the    problem
\begin{eqnarray*} 
&  &
u_{tt}-  t^{\ell} A(x,\partial_x)  u =0 \quad \mbox{for all} \quad (x,t) \in \Omega, \\
&  &
 u(x,0)= \varphi_0 (x), \quad u_t(x,0)= \varphi _1  (x )
\quad \mbox{for all} \quad x \in \widetilde{\Omega }\,.
\end{eqnarray*}
\end{corollary}

Finally, we remark that Theorem~\ref{T1.8}, Corollary~\ref{C1.4}, and Theorem~\ref{T1.8b}  are applicable to the case of  $\ell \in {\mathbb C}$, provided that $\ell \in {\mathbb C}\setminus \overline{ D_1((-3,0))}= \{ z \in {\mathbb C}\,|\, |z+3| > 1 \}$. In order to carry out the continuation of the integral transforms  in the complex plane $\ell\in {\mathbb C}$
we need some analyticity assumption on the solution of (\ref{8}), (\ref{9}). That question is out of the scope of this paper.

\section{Appendix}

In this section $ r \in (0,\phi (t)  - \phi (b))$, where $b\in (0,t) $ and  $\phi (t)  $ is defined in (\ref{phi_Tric}). 
The  
formula  15.3.6 of Ch.15\cite{A-S}  
 ties together points $z=0$ and $z=1$: 
\begin{eqnarray} 
\label{15.3.6}
 F  \left( a,b;c;z  \right) 
& = &
\frac{\Gamma (c)\Gamma (c-a-b)}{\Gamma (c-a)\Gamma (c-b)}F  \left( a,b;a+b-c+1;1-z  \right) \\
 &  &
 + (1-z)^{c-a-b}\frac{\Gamma (c)\Gamma (a+b-c)}{\Gamma (a)\Gamma (b)}F  \left( c-a,c-b;c-a-b+1;1-z  \right) \,, \quad |\arg (1-z)| <\pi \,. \nonumber  
\end{eqnarray}

\begin{lemma} 
\label{L5.1}
We have
\begin{eqnarray*}
&  &
\left|    
F \left(\gamma , \gamma ;1; \frac{(\phi (t)  - \phi (b))^2 - r^2} {(\phi (t)  + \phi (b))^2 - r^2} \right)   \right|\\
& \leq  &
 c \left[ 1+ \left(  (\phi (t)  + \phi (b))^2  -r^2 \right)^{ -1 }\phi (t)  \phi (b) 
+\left(  (\phi (t)  + \phi (b))^2  -r^2 \right)^{2\gamma -1  } \phi (t)^{1-2\gamma }  \phi (b) ^{1-2\gamma }  \right]  \,.
\end{eqnarray*}
\end{lemma}
\medskip

\noindent
{\bf Proof.}  We have due to (\ref{15.3.6}) and to the condition $\gamma <1 $ the following relation
\begin{eqnarray*}  
 F  \left( \gamma ,\gamma ;1;z  \right) 
& = &
\frac{ \Gamma (1-2\gamma )}{\Gamma (1-\gamma )^2  }F  \left( \gamma ,\gamma ;2\gamma  ;1-z  \right) \\
 &  &
 + (1-z)^{1-2\gamma }\frac{ \Gamma (2\gamma -1)}{\Gamma (\gamma )^2  }F  \left( 1-\gamma ,1-\gamma ;2-2\gamma ;1-z  \right) \,, \quad |\arg (1-z)| <\pi \,. \nonumber  
\end{eqnarray*}
Here
\begin{eqnarray*}
&  &
F  \left( \gamma ,\gamma ;2\gamma  ;1-z  \right)
=1+\frac{1}{2} \gamma  (1-z )+O\left((1-z)^2\right)\,,\\
&  &
F  \left( 1-\gamma ,1-\gamma ;2-2\gamma ;1-z  \right)
= 1+\frac{1}{2} (1-\gamma ) (1-z)+O\left((1-z )^2\right)
\end{eqnarray*}
imply
\begin{eqnarray*}
 F  \left( \gamma ,\gamma ;1;z  \right) 
& = &
\frac{ \Gamma (1-2\gamma )}{\Gamma (1-\gamma )^2  }\left[ 1+\frac{1}{2} \gamma  (1-z )+O\left((1-z)^2\right)\right]   \\
 &  &
 + (1-z)^{1-2\gamma }\frac{ \Gamma (2\gamma -1)}{\Gamma (\gamma )^2  }\left[  1+\frac{1}{2} (1-\gamma ) (1-z)+O\left((1-z )^2\right)\right]    \,. \nonumber  
\end{eqnarray*}
It remains to plug  $ z= \frac{(\phi (t)  - \phi (b))^2 - r^2} {(\phi (t)  + \phi (b))^2 - r^2} $ in the last formula.  
Lemma is proved. \hfill $\square$

\begin{lemma}  
\label{L5.2}
For every $ d_1, d_2 \in {\mathbb R}$  we have
\begin{eqnarray*}
&  & 
  \int_{ 0}^{ \phi (t)- \phi (b)}  \left(  (\phi (t)  + \phi (b))^2  -r^2 \right)^{-d_1\gamma -d_2}    dr  \\
&  = &
\left(\phi (t)  - \phi (b) \right)  \left(\phi (t)  + \phi (b)\right)^{-2   d_1\gamma -2d_2 } F \left(\frac{1}{2},d_2+ d_1\gamma ;\frac{3}{2};\frac{\left(\phi (t)- \phi (b)\right)^2}{\left(\phi (t)  + \phi (b)\right)^2}\right) \,.
\end{eqnarray*}
\end{lemma}
\medskip

\noindent
{\bf Proof.}  Indeed for $T>B>0$ we have
\begin{eqnarray*}
&  & 
 \int_{ 0}^{ T- B}  \left(  (T +B)^2  -r^2 \right)^{-d_1\gamma -d_2}    dr =
  (T-B) (B+T)^{-2 (d_2+\gamma  d_1)} F\left(\frac{1}{2},d_2+d_1 \gamma ;\frac{3}{2};\frac{(B-T)^2}{(B+T)^2}\right)\,.
\end{eqnarray*}
Lemma is proved. \hfill $\square$

\begin{lemma} 
\label{L5.4}
We have
\begin{eqnarray*} 
&  &
F \left(\frac{1}{2},  \gamma ;\frac{3}{2};\frac{\left(\phi (t)- \phi (b)\right)^2}{\left(\phi (t)  + \phi (b)\right)^2}\right)\\
& = &
\frac{\sqrt{\pi}\Gamma (1-\gamma )}{2\Gamma (\frac{3}{2}-\gamma )}
\left[   1 +\frac{1}{2} \left( \frac{4 \phi (t)   \phi (b)  } {(\phi (t)  + \phi (b))^2  } \right)+ O\left(\left( \frac{4 \phi (t)   \phi (b)  } {(\phi (t)  + \phi (b))^2  } \right)^2\right)\right]  \\
 &  &
 + \left( \frac{4 \phi (t)   \phi (b)  } {(\phi (t)  + \phi (b))^2  } \right)^{1-\gamma }\frac{ \Gamma (  \gamma -1)}{2\Gamma (\gamma )}
\left[ 1+\frac{\left(\frac{3}{2}- \gamma \right)}{2-\gamma } \left( \frac{4 \phi (t)   \phi (b)  } {(\phi (t)  + \phi (b))^2  } \right)+O\left(\left( \frac{4 \phi (t)   \phi (b)  } {(\phi (t)  + \phi (b))^2  } \right)^2\right) \right] \,.
\end{eqnarray*}
\end{lemma}
\medskip

\noindent
{\bf Proof.}  
Due to (\ref{15.3.6})
we have 
\begin{eqnarray*} 
 F \left(\frac{1}{2},  \gamma ;\frac{3}{2};z\right) 
& = &
\frac{\sqrt{\pi}\Gamma (1-\gamma )}{2\Gamma (\frac{3}{2}-\gamma )}
F  \left( \frac{1}{2},\gamma ; \gamma  ;1-z  \right) \\
 &  &
 + (1-z)^{1-\gamma }\frac{ \Gamma (  \gamma -1)}{2\Gamma (\gamma )}
F  \left( 1,\frac{3}{2}-\gamma ;2-\gamma ;1-z  \right) \,. 
\end{eqnarray*}
On the other hand,
\begin{eqnarray*} 
&   &  
F  \left( \frac{1}{2},\gamma ; \gamma  ;1-z  \right)= 1 +\frac{1}{2} (1-z)+ O\left((z-1)^2\right)\,,\\
 &  &
F  \left( 1,\frac{3}{2}-\gamma ;2-\gamma ;1-z  \right) 
=  1+\frac{\left(\frac{3}{2}- \gamma \right)}{2-\gamma } (1-z)+O\left((z-1)^2\right)\,.
\end{eqnarray*}
Hence, 
\begin{eqnarray*} 
 F \left(\frac{1}{2},  \gamma ;\frac{3}{2};z\right) 
& = &
\frac{\sqrt{\pi}\Gamma (1-\gamma )}{2\Gamma (\frac{3}{2}-\gamma )}
\left[   1 +\frac{1}{2} (1-z)+ O\left((z-1)^2\right)\right]  \\
 &  &
 + (1-z)^{1-\gamma }\frac{ \Gamma (  \gamma -1)}{2\Gamma (\gamma )}
\left[ 1+\frac{\left(\frac{3}{2}- \gamma \right)}{2-\gamma } (1-z)+O\left((z-1)^2\right) \right] \,.
\end{eqnarray*}
It remains to plug $z= \frac{\left(\phi (t)- \phi (b)\right)^2}{\left(\phi (t)  + \phi (b)\right)^2}$ in the last formula.  
Lemma is proved. \hfill $\square$

\begin{lemma} 
\label{L5.5}
We have
\begin{eqnarray*} 
&  &
 F  \left( \frac{1}{2},\gamma +1;\frac{3}{2};\frac{(\phi (t)  - \phi (b))^2  } {(\phi (t)  + \phi (b))^2  } \right) \\
& = &
\frac{\sqrt{\pi}\Gamma ( - \gamma  )}{2\Gamma (\frac{1}{2}- \gamma  )}
\left[ 1+\frac{ 1}{2} \left( \frac{4 \phi (t)   \phi (b)  } {(\phi (t)  + \phi (b))^2  } \right) +O\left( \left( \frac{4 \phi (t)   \phi (b)  } {(\phi (t)  + \phi (b))^2  } \right)^2\right)\right]   \\
 &  &
 + \left( \frac{4 \phi (t)   \phi (b)  } {(\phi (t)  + \phi (b))^2  } \right)^{-\gamma }\frac{ \Gamma ( \gamma  )}{2\Gamma (\gamma +1)}
\left[ 1-\frac{\left(\gamma -\frac{1}{2}\right)}{1-\gamma }  \left( \frac{4 \phi (t)   \phi (b)  } {(\phi (t)  + \phi (b))^2  } \right)+O\left( \left( \frac{4 \phi (t)   \phi (b)  } {(\phi (t)  + \phi (b))^2  } \right)^2\right)\right] \,.
\end{eqnarray*}
\end{lemma}
\medskip

\noindent
{\bf Proof.}
Due to (\ref{15.3.6})
we obtain
\begin{eqnarray*} 
 F  \left( \frac{1}{2},\gamma +1;\frac{3}{2};z  \right) 
& = &
\frac{\sqrt{\pi}\Gamma ( - \gamma  )}{2\Gamma (\frac{1}{2}- \gamma  )}
F  \left( \frac{1}{2},\gamma +1; \gamma  +1;1-z  \right) \\
 &  &
 + (1-z)^{-\gamma }\frac{ \Gamma ( \gamma  )}{2\Gamma (\gamma +1)}
F  \left( 1,\frac{1}{2}- \gamma ;1-\gamma  ;1-z  \right)\,, \quad |\arg (1-z)| <\pi \,, \nonumber  
\end{eqnarray*}
where
\begin{eqnarray*} 
&  &
F  \left( \frac{1}{2},\gamma +1; \gamma  +1;1-z  \right) 
= 1+\frac{ 1}{2}(1-z) +O\left((z-1)^2\right)\,,\\
 &  & 
F  \left( 1,\frac{1}{2}- \gamma ;1-\gamma  ;1-z  \right) 
= 1-\frac{\left(\gamma -\frac{1}{2}\right)}{1-\gamma } (1-z )+O\left((z-1)^2\right)\,.
\end{eqnarray*}
It follows
\begin{eqnarray*} 
 F  \left( \frac{1}{2},\gamma +1;\frac{3}{2};z  \right) 
& = &
\frac{\sqrt{\pi}\Gamma ( - \gamma  )}{2\Gamma (\frac{1}{2}- \gamma  )}
\left[ 1+\frac{ 1}{2}(1-z) +O\left((z-1)^2\right)\right]   \\
 &  &
 + (1-z)^{-\gamma }\frac{ \Gamma ( \gamma  )}{2\Gamma (\gamma +1)}
\left[ 1-\frac{\left(\gamma -\frac{1}{2}\right)}{1-\gamma } (1-z )+O\left((z-1)^2\right)\right] \,.
\end{eqnarray*}
It remains to plug $z= \frac{(\phi (t)  - \phi (b))^2  } {(\phi (t)  + \phi (b))^2  }$ in the last formula.  
Lemma is proved. \hfill $\square$

\begin{lemma} 
\label{L5.6}
We have
\begin{eqnarray*} 
&  & 
F \left(\frac{1}{2},1-\gamma ;\frac{3}{2};\frac{\left(\phi (t)- \phi (b)\right)^2}{\left(\phi (t)  + \phi (b)\right)^2}\right)  \\
& = & 
\frac{\sqrt{\pi}\Gamma ( \gamma )}{2\Gamma (\frac{1}{2}+\gamma )}\left[ 1+\frac{1}{2}\left( \frac{4 \phi (t)   \phi (b)  } {(\phi (t)  + \phi (b))^2  } \right)+O\left(\left( \frac{4 \phi (t)   \phi (b)  } {(\phi (t)  + \phi (b))^2  } \right)^2\right)\right] \\
 &  &
 + \left( \frac{4 \phi (t)   \phi (b)  } {(\phi (t)  + \phi (b))^2  } \right)^{ \gamma }\frac{\Gamma ( -\gamma )}{2\Gamma (1-\gamma)}\left[  1+ \left( \frac{1}{2 \gamma }+1\right) \left( \frac{4 \phi (t)   \phi (b)  } {(\phi (t)  + \phi (b))^2  } \right)+ O\left(\left( \frac{4 \phi (t)   \phi (b)  } {(\phi (t)  + \phi (b))^2  } \right)^2\right)\right] \,.
\end{eqnarray*}  
\end{lemma}
\medskip

\noindent
{\bf Proof.} 
Due to  (\ref{15.3.6})
we obtain  
\[ 
F \left(\frac{1}{2},1-\gamma ;\frac{3}{2};z\right) \\
  =  
\frac{\sqrt{\pi}\Gamma ( \gamma )}{2\Gamma (\frac{1}{2}+\gamma )}F  \left( \frac{1}{2},1-\gamma;  1-\gamma ;1-z  \right)  
 + (1-z)^{ \gamma }\frac{ \Gamma ( -\gamma )}{2\Gamma (1-\gamma)}F  \left( 1,\frac{1}{2}+\gamma ; \gamma ;1-z  \right) \,,    
\]
where 
\begin{eqnarray*}  
&  &
F  \left( \frac{1}{2},1-\gamma;  1-\gamma ;1-z  \right) =  
1+\frac{1-z}{2}+O\left((z-1)^2\right)\,,\\
 &  &
F  \left( 1,\frac{1}{2}+\gamma ; \gamma ;1-z  \right) = 
 1+ \left( \frac{1}{2 \gamma }+1\right) (1-z)+ O\left((z-1)^2\right)\,.
\end{eqnarray*}
It follows
\begin{eqnarray*}  
F \left(\frac{1}{2},1-\gamma ;\frac{3}{2};\frac{\left(\phi (t)- \phi (b)\right)^2}{\left(\phi (t)  + \phi (b)\right)^2}\right)  
& = & 
\frac{\sqrt{\pi}\Gamma ( \gamma )}{ 2\Gamma (\frac{1}{2}+\gamma )}\left[ 1+\frac{1}{2}(1-z)+O\left((z-1)^2\right)\right] \\
 &  &
 + (1-z)^{ \gamma }\frac{ \Gamma ( -\gamma )}{2\Gamma (1-\gamma)}\left[  1+ \left( \frac{1}{2 \gamma }+1\right) (1-z)+ O\left((z-1)^2\right)\right] \,.
\end{eqnarray*} 
It remains to plug $z= \frac{\left(\phi (t)- \phi (b)\right)^2}{\left(\phi (t)  + \phi (b)\right)^2}$ in the last formula. 
Lemma is proved. \hfill $\square$

\begin{lemma} 
\label{L5.7}
We have
\begin{eqnarray*}
&  &    
\Bigg| F \left(\gamma+1 , \gamma +1;2; \frac{(\phi (t)  - \phi (b))^2 - r^2} {(\phi (t)  + \phi (b))^2 - r^2} \right) \Bigg|  \\
& \lesssim  &
\sum_{i=0,1,2}
 \left( \frac{  \phi (t)   \phi (b)  } {(\phi (t)  + \phi (b))^2  -r^2 } \right)^i  
 + \left( \frac{ \phi (t)   \phi (b)  } {(\phi (t)  + \phi (b))^2  -r^2 } \right)^{-2\gamma }\sum_{i=0,1,2}
 \left( \frac{  \phi (t)   \phi (b)  } {(\phi (t)  + \phi (b))^2  -r^2 } \right)^i   \,.
\end{eqnarray*}
\end{lemma}
\medskip

\noindent
{\bf Proof.} Due to (\ref{15.3.6}) we have
 \begin{eqnarray*} 
F \left(\gamma+1 , \gamma +1;2; z \right)   
& = &
\frac{2\Gamma (-2\gamma )}{\Gamma (1-\gamma )^2 }F  \left( \gamma +1 , \gamma +1 ;2\gamma +1;1-z  \right) \\
 &  &
 + (1-z)^{-2\gamma }\frac{2\Gamma (2\gamma  )}{\Gamma ((\gamma +1))^2}
F  \left(1- \gamma ,1- \gamma;1-2\gamma ;1-z  \right)\,.
\end{eqnarray*}
Hence
 \begin{eqnarray*}
&  &    
F  \left( \gamma +1 , \gamma +1 ;2\gamma +1;1-z  \right) 
=1+\frac{(\gamma +1)^2 }{2 \gamma +1}(1-z)+O\left((z-1)^2\right)\,,\\
 &  &
F  \left(1- \gamma ,1- \gamma;1-2\gamma ;1-z  \right) 
=1-\frac{(\gamma -1)^2}{2 \gamma -1} (1-z)+O\left((z-1)^2\right)\,.
\end{eqnarray*}
Thus
\begin{eqnarray*}
F \left(\gamma+1 , \gamma +1;2; z \right)  
& = &
\frac{2\Gamma (-2\gamma )}{\Gamma (1-\gamma )^2 }
\Bigg[1+\frac{(\gamma +1)^2 }{2 \gamma +1}(1-z)+O\left((z-1)^2\right) \Bigg]  \\
 &  &
 + (1-z)^{-2\gamma }\frac{2\Gamma (2\gamma  )}{\Gamma ((\gamma +1))^2}
\Bigg[1-\frac{(\gamma -1)^2}{2 \gamma -1} (1-z)+O\left((z-1)^2\right) \Bigg] \,.
\end{eqnarray*}
It remains to plug $z=\frac{(\phi (t)  - \phi (b))^2 - r^2} {(\phi (t)  + \phi (b))^2 - r^2}  $ in the last formula
\begin{eqnarray*}
&  &    
F \left(\gamma+1 , \gamma +1;2; \frac{(\phi (t)  - \phi (b))^2 - r^2} {(\phi (t)  + \phi (b))^2 - r^2} \right)   \\
& = &
\frac{2\Gamma (-2\gamma )}{\Gamma (1-\gamma )^2 }
\Bigg[1+\frac{(\gamma +1)^2 }{2 \gamma +1}\left( \frac{4 \phi (t)   \phi (b)  } {(\phi (t)  + \phi (b))^2 -r^2 } \right)+O\left(\left( \frac{4 \phi (t)   \phi (b)  } {(\phi (t)  + \phi (b))^2  -r^2 } \right)^2 \right) \Bigg]  \\
 &  &
 + \left( \frac{4 \phi (t)   \phi (b)  } {(\phi (t)  + \phi (b))^2  -r^2 } \right)^{-2\gamma }\frac{2\Gamma (2\gamma  )}{\Gamma ( \gamma +1 )^2}\\
 &  &
 \hspace{2cm} \times 
\Bigg[1-\frac{(\gamma -1)^2}{2 \gamma -1} \left( \frac{4 \phi (t)   \phi (b)  } {(\phi (t)  + \phi (b))^2  -r^2 } \right)+O\left(\left( \frac{4 \phi (t)   \phi (b)  } {(\phi (t)  + \phi (b))^2  -r^2 } \right)^2\right) \Bigg] \,.
\end{eqnarray*}
Lemma is proved. \hfill $\square$

\begin{lemma}
\label{L5.9}
For all $z \in (0,1)$ we have
\begin{eqnarray*}  
\left| F \left(\frac{1}{2},2+  \gamma ;\frac{3}{2}; z \right) \right| 
& \lesssim & 
 1+ (1-z)^{-1 }  +(1-z)^{-\gamma } +(1-z)^{-\gamma -1} \,,\\
\left| F \left(\frac{1}{2}, 2-\gamma ;\frac{3}{2};z\right) \right| 
&  \lesssim  &
1+(1-z)^{\gamma -1 }  +(1-z)^{\gamma   }  \,.
\end{eqnarray*}
\end{lemma}
We left the proof of the last lemma to the reader.
\begin{lemma}
\label{L5.8} If $ \ell \in (- 2,4 )$, then  
\begin{eqnarray}
\label{46}
&  &
\Bigg| t^k  \int_{ 0}^{t} db \,b^{k+1}
  \int_{ 0}^{ \phi (t)- \phi (b)} \left(  (\phi (t)  + \phi (b))^2  -r^2 \right)^{-2-\gamma } \\
&  &
\hspace{2.5cm} \times\left(   \phi^2 (t)  - \phi^2 (b)    + 
   r^2  \right)  F_z \left(\gamma , \gamma ;1; \beta  (t,b,r) \right)   w(x,r;b )\,dr  \Bigg|
 \lesssim  
t \quad \mbox{\rm for all }  \quad t  >0\,. \nonumber
\end{eqnarray} 
\end{lemma}
\medskip

\noindent
{\bf Proof.} Indeed, according to the previous lemma   we have
\begin{eqnarray*}
&  &
  \int_{ 0}^{ \phi (t)- \phi (b)} \left(  (\phi (t)  + \phi (b))^2  -r^2 \right)^{-2-\gamma } \left(   \phi^2 (t)  - \phi^2 (b)    +    r^2  \right)\\
&  &
\times\Bigg[\sum_{i=0 }^2
 \left( \frac{  \phi (t)   \phi (b)  } {(\phi (t)  + \phi (b))^2  -r^2 } \right)^i  
 +  \sum_{i=0 }^2
 \left( \frac{  \phi (t)   \phi (b)  } {(\phi (t)  + \phi (b))^2  -r^2 } \right)^{i-2\gamma }   \Bigg]\,dr\\
 &  \lesssim &
\phi (t) (\phi (t)- \phi (b)) \sum_{i=0 }^2
\phi (t)^i   \phi (b)^i \\
&  &
\times \int_{ 0}^{ \phi (t)- \phi (b)}\Bigg[\left(   (\phi (t)  + \phi (b))^2  -r^2  \right)^{-i-2-\gamma }  
 +  
 ( \phi (t)  + \phi (b))^{-2\gamma }\left( (\phi (t)  + \phi (b))^2  -r^2  \right)^{-i-2+\gamma }   \Bigg] dr\,.
 \end{eqnarray*} 
 Then we apply   Lemma~\ref{L5.2} 
 \begin{eqnarray*} 
&  \lesssim  &
\phi (t) (\phi (t)- \phi (b))^2 \sum_{i=0 }^2 \phi (t)^i\phi (b)^i  \\
&  &
\times \Bigg[     \left(\phi (t)  + \phi (b)\right)^{ -2    \gamma -2i-4 } F \left(\frac{1}{2},i+2+\gamma ;\frac{3}{2};\frac{\left(\phi (t)- \phi (b)\right)^2}{\left(\phi (t)  + \phi (b)\right)^2}\right)\\
&  & 
 +    
( \phi (t)  + \phi (b))^{-2\gamma } \left(\phi (t)  + \phi (b)\right)^{2\gamma -2i-4 } 
F \left(\frac{1}{2},i+2-\gamma ;\frac{3}{2};\frac{\left(\phi (t)- \phi (b)\right)^2}{\left(\phi (t)  + \phi (b)\right)^2}\right)\Bigg]\,.
\end{eqnarray*} 
For     $z=  \frac{\left(\phi (t)- \phi (b)\right)^2}{\left(\phi (t)  + \phi (b)\right)^2} \in ( 0,1) $ taking into account condition $\ell <4 $  we continue as follows  
 \begin{eqnarray*} 
&  \lesssim  &
\phi (t) (\phi (t)- \phi (b))^2 
\Bigg\{ \left(\phi (t)  + \phi (b)\right)^{ -2    \gamma -4  } \Big[ 1+ (1-z)^{-1-\gamma }\Big] \\
&  &
+ \phi (t)  \phi (b) \left(\phi (t)  + \phi (b)\right)^{ -2    \gamma -6 }  \Big[ 1+ (1-z)^{-2-\gamma }  \Big]\\
&  &
+ \phi (t)^2 \phi (b)^2\left(\phi (t)  + \phi (b)\right)^{ -2    \gamma -8 } \Big[ 1+ (1-z)^{-3-\gamma }    \Big]\\
&  &
  +    
  \phi (t)^{ -2\gamma } \phi (b)^{ -2\gamma } \left(\phi (t)  + \phi (b)\right)^{2\gamma   -4  } \Big[  1+ (1-z)^{-1+\gamma } \Big] \\
  &  &
+ \phi (t)^{1-2\gamma } \phi (b)^{1-2\gamma } \left(\phi (t)  + \phi (b)\right)^{2\gamma   -6 } \Big[  1+ (1-z)^{-2+\gamma }    \Big]\\
&  &
+ \phi (t)^{2-2\gamma } \phi (b)^{2-2\gamma } \left(\phi (t)  + \phi (b)\right)^{2\gamma   -8 } \Big[   1+ (1-z)^{-3+\gamma }  \Big]\Bigg\}\,.
\end{eqnarray*} 
Hence, the integral (\ref{46})   can be estimated from above by $\phi^{1-2\gamma } (t)  \lesssim t$.   Lemma is proved. \hfill $\square$

\end{document}